\documentclass[12pt,a4paper]{amsart}

\usepackage{amsfonts}
\usepackage{amsmath}
\usepackage{amsthm}
\usepackage{amssymb}
\usepackage{mathrsfs}
\usepackage{bm}
\usepackage{here}
\usepackage{comment}

\usepackage{tikz}

\usepackage{marginnote}

\usepackage{mathtools}
\usetikzlibrary{patterns}

\mathtoolsset{showonlyrefs}

\allowdisplaybreaks[1]

\newcommand{\be}{\begin{equation}} 
\newcommand{\ee}{\end{equation}}
\newcommand{\bea}{\begin{eqnarray}} 
\newcommand{\eea}{\end{eqnarray}}
\newcommand{\bean}{\begin{eqnarray*}} 
\newcommand{\eean}{\end{eqnarray*}}

\def\mn{|\!\!|}

\def\mn2{|\!\!|_{M^{d/2}}}

\newcommand{\N}{\nabla}

\newcommand{\x}{\xi}

\def\<{\langle }

\renewcommand{\qed}{\qquad\kern1pt   
   \vbox{\hrule height 0.6pt      
         \hbox{\vrule width 0.6pt 
               \vbox{\vskip 6pt}  
               \hskip 3pt
              \vrule width 1.3pt} 
         \hrule depth 1.3pt}     
   \kern1pt}

\newtheorem{theorem}{Theorem}[section]
\newtheorem{proposition}[theorem]{Proposition}
\newtheorem{lemma}[theorem]{Lemma}

\theoremstyle{definition}
\newtheorem{definition}[theorem]{Definition}
\def\qed{\hfill $\square$}

\theoremstyle{remark}
\newtheorem{remark}[theorem]{Remark}

\numberwithin{equation}{section}
\numberwithin{theorem}{section}

\author[H. Wakui]{Hiroshi Wakui}
\address[H. Wakui]{ 
 Faculty of Engineering, 
 University of Fukui, 
 Fukui-shi,
 Fukui \hbox{910-8507}, Japan}
\email{hwakui@u-fukui.ac.jp}

\author[T. Yamada]{Tetsuya Yamada}
\address[T. Yamada]{Course of General Education, National Institute of Technology, Fukui College, Sabae-shi, Fukui \hbox{916-8507}, Japan}
\email{yamada@fukui-nct.ac.jp}

\title[]{Stability of constant steady states of an
attraction-repulsion chemotaxis system
}

\begin{document}

\begin{abstract}
The Cauchy problem for the attraction-repulsion chemotaxis 
system
in the whole $n$-dimensional space 
has uncountable constant steady states. 
In the attraction chemotaxis 
system, each positive constant steady state is stable if it is in a certain region.
On the other hand, in the repulsion chemotaxis system, every positive constant steady state is stable.
Our main purpose of this paper is to give a suitable condition under which the attraction-repulsion 
chemotaxis system has also stable constant steady states. 
 \end{abstract}

 \keywords{Attraction-repulsion chemotaxis system, Constant steady states, Stability of solutions}

 \subjclass[2020]{35Q92; 35K55; 35B35; 35B40; 35K15; 92C17}

 \date{ \today}
 \maketitle

 \baselineskip=14.5pt


 \section{Introduction}


In this paper, we consider a stability of the following attraction-repulsion chemotaxis 
system:
\begin{equation}\label{eq;at-re}\tag{ARKS}
  \left\{
  \begin{split}
  	\partial_t u - \Delta u \, + \, \, &\nabla \cdot \big(  u \nabla (\beta_1 \psi_1 - \beta_2 \psi_2) \big) = 0,&& \quad t>0,\quad x \in \mathbb{R}^{n},\\
 - &\Delta \psi_1 + \lambda_1 \psi_1 = u,&& \quad  t>0,\quad x \in \mathbb{R}^{n},\\
 - &\Delta \psi_2 + \lambda_2 \psi_2 = u ,&& \quad  t>0,\quad x \in \mathbb{R}^{n},\\
 &u(0,x)  = u_0(x),&& \quad x \in \mathbb{R}^{n}.
 \end{split}
  \right.
\end{equation}
where $\beta_1, \beta_2, \lambda_1, \lambda_2 > 0$ are given constants.
In system \eqref{eq;at-re}, 
{the} function 
$
 u
 =
 u(t,x)
$ 
is the density of cells, 
{the} function 
$
 \psi_1
 =
 \psi_1(t,x)
$
stands for the concentration of an attraction effect of cells, and 
{the} 
function 
$
 \psi_2
 =
 \psi_2(t,x)
$ 
represents the concentration of a repulsive effect of cells.
This system \eqref{eq;at-re} is already investigated in some frameworks on whole space $\mathbb{R}^n$. 
For instance, blow-up or global existence in system \eqref{eq;at-re} has been discussed in \cite{Ho}, \cite{HoOg2}, \cite{NaYa3}, \cite{NaYa1}, \cite{NaYa2}, \cite{ShWa}, 
and some related results is cited on \cite{JiLi}, \cite{NaSeYa2}, \cite{NaSeYa1}, \cite{Ya1}.
In the above results,  
the framework of solutions to system \eqref{eq;at-re} heavily depends on 
integrability of functions at $|x| \to \infty$, and constant functions 
are excluded in such a framework. 
Recently, by \cite{CyKaKrWa} and \cite{Su1},
a mathematical theory concerning local or global solutions to 
system \eqref{eq;at-re} with {$\beta_j=0$ and $\psi_j\equiv 0$ for each $j=1,2$} has been
developed  in the uniformly local Lebesgue spaces 
$
 L^p_{\mathrm{uloc}}(\mathbb{R}^n)
$, 
which contains the constant functions.

For any given constant $A \in \mathbb{R}$, 
the pair of function 
$
 \big( 
  u, \psi_1, \psi_2 
 \big) 
 = 
 \big( 
  A, {A}/{\lambda_1}, {A}/{\lambda_2} 
 \big)
$ 
solves system \eqref{eq;at-re}.
In \cite{CyKaKrWa}, 
Cygan, Karch, Krawczyk and the first author of this paper consider system \eqref{eq;at-re} with $\beta_2=0$ and $\psi_2\equiv 0$:
\begin{equation}\label{eq;at}\tag{AKS}
  \left\{
  \begin{split}
  \partial_t  u - \Delta u  \, + \, \, &\nabla \cdot \big(  u  \nabla (\beta_1 \psi_1) \big) = 0,&& \quad t>0,\quad x \in \mathbb{R}^n,\\
 - &\Delta \psi_1 + \lambda_1 \psi_1 = u ,&& \quad  t>0,\quad x \in \mathbb{R}^n,\\
 &u (0,x)  = u_0(x),&& \quad x \in \mathbb{R}^{n},
 \end{split}
  \right.
\end{equation}
and study stability of constant steady states
$
 (u , \psi_1) 
 =
 \big(  
  A, {A}/{\lambda_1} 
 \big)
$
of system \eqref{eq;at}. More precisely, if
 $
   0 
   < 
   A 
   <
   {\lambda_1}/{\beta_1}
 $,
then each constant steady state 
$
 ( u , \psi_1) 
 = 
 \big( 
  A, {A}/{\lambda_1} 
 \big)
$
is uniform stable 
in a suitable Lebesgue space and is asymptotically stable in higher regular Lebesgue spaces.
On the other hand, 
when $\beta_1 = 0$ and $\psi_1 \equiv 0$ in system \eqref{eq;at-re}, 
it can be reformulated as follows:
\begin{equation}\label{eq;re}\tag{RKS}
  \left\{
  \begin{split}
 \partial_t u - \Delta u \, - \, \, &\nabla \cdot \big( u \nabla (\beta_2 \psi_2) \big) = 0,&& \quad t>0,\quad x \in \mathbb{R}^n,\\
 - &\Delta \psi_2 + \lambda_2 \psi_2 = u,&& \quad  t>0,\quad x \in \mathbb{R}^n,\\
 &u(0,x)  = u_0(x),&& \quad x \in \mathbb{R}^n,
 \end{split}
  \right.
\end{equation}
and if $A>0$, 
then each constant steady state 
$
 (u, \psi_2) 
 = 
 \big( 
 A, A/\lambda_2
 \big)
$ is stable 
by the same argument as in \cite{CyKaKrWa}. 
As observed above,
system \eqref{eq;at-re} 
will have more complicated 
structure in the stability of constant steady states  
in contrast to system \eqref{eq;at} and system \eqref{eq;re}. 
In fact, Liu, Shi and Wang 
\cite{LiShWa} have studied a fully parabolic attraction repulsion chemotaxis system on bounded domains in $\mathbb{R}^n$ ($n\ge 1$) with Neumann boundary condition, 
and shown exponential stability of constant steady states when $\beta_1/\beta_2\le \min\{\lambda_1/\lambda_2, 1/(\lambda_1/\lambda_2)\}$ (see also \cite{LiMuWa}). 
However, as far as we know, it is unclear whether each constant steady state 
$
 \big( 
  u, \psi_1, \psi_2 
 \big) 
 = 
 \big( 
  A, {A}/{\lambda_1}, {A}/{\lambda_2} 
 \big)
$ of system \eqref{eq;at-re} is stable or not. 
The main goal of our paper is to solve this problem. 
In detail, we show that 
a small initial perturbation gives asymptotic stability of constant steady {states}
$
 (u,\psi_1,\psi_2)
  =
 (A,A/\lambda_1,A/\lambda_2)
$
of system \eqref{eq;at-re} under suitable range on $A>0$ (see Theorem \ref{thm;1} below).

Before closing this section, we now 
introduce some notation used in this paper. 

\subsection{Notation} 
Let 
$
 {\mathbb Z}_+
$ 
be the set of nonnegative integers. 
For 
$ 
 x=(x_1,x_2,\ldots,x_n)\in {\mathbb R}^n
$ 
and 
$
 \alpha
 =
 (
  \alpha_1, 
  \alpha_2,
  \ldots,
  \alpha_n
 )
 \in 
 \mathbb{Z}_+^n
$, 
set 
$ 
  x^\alpha
  :=
  x_1^{\alpha_1}
  x_{2}^{\alpha_2}
  \cdots 
  x_{n}^{\alpha_n}
$ 
and 
$
 |\alpha|
 :=
  \alpha_1
  +
  \alpha_2
  +
  \cdots
  +
  \alpha_n
$. 
The symbols 
$
 \frac{\partial}{\partial t}
$ 
and 
$
 \frac{\partial}{\partial x_j}
$ 
are defined by 
$
 \partial_t
$
and 
$ 
 \partial_j
$, 
respectively, 
and define 
$
 \nabla
 =
 {}^t(
  \partial_1, 
  \partial_2,
  \ldots, 
  \partial_n
  )
$ 
and 
$
 \partial_x^\alpha
 =
 \partial_1^{\alpha_1}
 \partial_2^{\alpha_2}
 \cdots
 \partial_n^{\alpha_n}
$ 
for 
$
 \alpha
 =
 (
  \alpha_1,
  \alpha_2,
  \ldots, 
  \alpha_n
 )
 \in 
 \mathbb{Z}_{+}^n
$. 
For 
$
 1
 \le 
 p
 \le 
 \infty
$, 
the norm of the usual Lebesgue space $L^p({\mathbb R}^n)$ is denoted by $\|\cdot\|_p$. 
Let $\mathcal{S}$ be the set of rapidly decreasing functions on ${\mathbb R}^n$. 
For $f\in {\mathcal S}$, the symbol ${\mathcal F}[f]$ stands for the Fourier transform of $f$, that is, 
$
 {\mathcal F}[f](\xi)
 :=
 (2\pi)^{-\frac{n}{2}}
 \int_{{\mathbb R}^n} e^{-ix \cdot \xi}
  f(x)
 \, \mathrm{d}x, 
 \xi
 \in \mathbb R^n
$, 
and the inverse Fourier transform of $f$ is defined by 
$
 {\mathcal F}^{-1}[f](x)
 :=
 (2\pi)^{-\frac{n}{2}}
 \int_{\mathbb R^n} 
  e^{ix \cdot \xi}
   f(\xi)
  \, \mathrm{d}\xi
 , x \in \mathbb{R}^n
$. 
Let $C$ be a positive constant which may change from line to line. In particular, we write $C(\ast,\ldots,\ast)$ for a positive constant depending on the quantities appearing in parentheses.

\section{Main result}
We first solve the second and the third equations of system \eqref{eq;at-re} to obtain
$
 \psi_j
 =
 B_{\lambda_j}
 *u
$
($j=1,2$), 
where $B_{\lambda}=B_{\lambda}(x)$ stands for the Bessel function 
defined later in relation \eqref{eq;Bessel-kernel}.
Then we reduce system \eqref{eq;at-re} to the following problem:
\begin{equation}\label{eq;simplified-ar-ks}\tag{P}
 \left\{
 \begin{split}
  &\partial u - \Delta u 
  + \nabla \cdot \big( u\nabla (\beta_1 B_{\lambda_1}*u - \beta_2 B_{\lambda_2}*u ) \big) = 0,\qquad && t > 0, x \in \mathbb{R}^n,\\
  &u(0,x) = u_0(x),\qquad && x \in \mathbb{R}^n.
 \end{split} 
 \right.
\end{equation}
We now states the result on stability of constant steady states of problem \eqref{eq;simplified-ar-ks} as follows.

\begin{theorem}\label{thm;1}
Let 
  $    
    \beta_1, 
    \beta_2, 
    \lambda_1, 
    \lambda_2 > 0,
    \lambda_1 
    \not = 
    \lambda_2
  $, 
  and 
  $p$ be
  \[
   \left\{
   \begin{split}
    p & \, = 1 ,\qquad && n =1 ,\\
    \frac{n}2 < & \, p \le  n,\qquad && n \ge 2
   \end{split}
   \right.
  \]
  and 
  $ n < q \le 2p $.
 Assume that 
\begin{equation}\label{eq;cond-A}
   0
   <
   A
   <
   A_*,
\end{equation}
where
\begin{equation}\label{eq;cond_A}
 A_*
 :=
\begin{cases}
   \infty, 
   &
    \dfrac{\beta_1}{\beta_2} \le 1,\quad 
    \dfrac{\beta_1}{\beta_2} \le \dfrac{\lambda_1}{\lambda_2}, \vspace{2mm}
    \\
    \dfrac{
     \lambda_1\lambda_2
    }{
     \beta_1 \lambda_2 - \beta_2 \lambda_1
    }, 
    &\dfrac{\beta_1}{\beta_2} > \dfrac{\lambda_1}{\lambda_2},\quad \displaystyle \dfrac{\beta_1}{\beta_2} \ge \bigg( \dfrac{\lambda_1}{\lambda_2} \bigg)^2, \vspace{2mm} \\
 \displaystyle \frac{\lambda_1 - \lambda_2}{(\sqrt{\beta_1}-\sqrt{\beta_2})^2}, &\displaystyle \dfrac{\beta_1}{\beta_2} < \bigg( \dfrac{\lambda_1}{\lambda_2} \bigg)^2,\quad \dfrac{\beta_1}{\beta_2} >1.
 \end{cases}
\end{equation}
Then there exists $\varepsilon_0 > 0$ such that for every $v_0 \in L^p(\mathbb{R}^n)$ with
$
 \|
  v_0
 \|_p
 <
 \varepsilon_0
$,
problem \eqref{eq;simplified-ar-ks} 
corresponding to the initial data 
$
 u_0
 =
 A
 +
 v_0
$
has a unique global-in-time mild solution $u = u(t,x)$ 
satisfying 
$
 u
 -
 A 
 \in
 C([0,\infty);L^p(\mathbb{R}^n)) 
$
with 
\[
 \sup_{t>0}
 \| 
  u(t) - A 
 \|_p
 +
 \sup_{t > 0}
 t^{\frac{n}2(\frac1p-\frac1q)}
 \| 
  u(t) - A 
 \|_q
 \le
 C 
 \|
  u_0
  -
  A
 \|_p
\]
for some $C > 0$.
\end{theorem}

We observe from Theorem \ref{thm;1} that we need a restriction on $A>0$ in spite of the repulsion dominated case $\beta_1/\beta_2 < 1$ with 
$
 \beta_1 / \beta_2
 >
 \lambda_1 / \lambda_2
$ 
even though it is well-known in \cite{NaYa3}, \cite{ShWa} that uniform boundedness of solutions to problem \eqref{eq;simplified-ar-ks} is immediately obtained in the canonical $L^p$-framework in the repulsion dominated case $\beta_1/\beta_2 \le 1$.
Moreover, 
in the suitable $L^p$-framework setting,
large time behavior of solutions to problem \eqref{eq;simplified-ar-ks} definitely depend on the relation 
between $\beta_1, \beta_2$ and $\lambda_1, \lambda_2$.
Indeed,
if $\beta_1/\beta_2 > 1$ and $n = 2$,
then
there exists an enough 
large initial data $u_0$ satisfying 
$
 \|
  u_0
 \|_1
 >
 {8\pi}/(\beta_1-\beta_2)
$
such that
the solution $u$ corresponding to the initial data $u_0$ blows up in finite time, see also \cite{ShWa}.
When $\lambda_1/\lambda_2 < 1$ and $\beta_1 / \beta_2 \ge \lambda_1 / \lambda_2$ for $n \ge 3$,
unboundedness of solutions to problem \eqref{eq;simplified-ar-ks} may occur with
enough large initial data, see \cite{Ho}.
These observation above implies that ratio $\beta_1/\beta_2$, $\lambda_1/\lambda_2$ and relation between $\beta_1/\beta_2$
and $\lambda_1/\lambda_2$ plays an important role to classify large time behavior of solutions to problem \eqref{eq;simplified-ar-ks}.


\section{Preliminaries}

\subsection{Local-in-time solutions to problem \eqref{eq;simplified-ar-ks} in ${\mathcal L}_{\mathrm{uloc}}^p(\mathbb{R}^n)$} 
Let $1 \le p < \infty$. 
Then the uniformly local Lebesgue spaces $L_{\text{uloc}}^p({\mathbb R}^n)$ is defined by 
\[
 L_{\text{uloc}}^p(\mathbb{R}^n)
 \coloneqq 
 \left\{ 
  f \in L_{\text{loc}}^p(\mathbb{R}^n) 
  \,
  \mid 
  \, 
  \| f \|_{p, \mathrm{uloc}}
  \coloneqq 
  \sup_{x\in \mathbb{R}^n}
   \left(
    \int_{|y-x|<1}|f(y)|^{p} \, \mathrm{d}y
   \right)^{1/p}
   <
   \infty
 \right\}
\]
and also 
\[
 \mathcal{L}_{\mathrm{uloc}}^p(\mathbb{R}^n)
 \coloneqq 
 \overline{BUC(\mathbb{R}^n)}^{\| \cdot \|_{p, \text{uloc}}},
\]
where 
$
 L_{\text{loc}}^p({\mathbb R}^n)
$ 
denotes the set of locally $p$-integrable functions on $\mathbb{R}^n$, 
and $BUC(\mathbb{R}^n)$ stands for the set of bounded uniformly continuous functions on $\mathbb{R}^n$. 

For any $v_0 \in L^p(\mathbb{R}^n)$ and any $A \in \mathbb{R}$,
the function 
$
 u_0 
 \coloneqq 
 A
 +
 v_0
$ 
belongs to 
$
 \mathcal{L}_{\mathrm{uloc}}^p(\mathbb{R}^n)
$ 
with 
$
 1 
 \le 
 p 
 < 
 \infty
$ (see \cite[Remark 5.3]{CyKaKrWa}). 
Then it has been shown that the following local well-posedness for problem \eqref{eq;simplified-ar-ks} holds in 
$
 \mathcal{L}_{\mathrm{uloc}}^p(\mathbb{R}^n)
$ 
by the same argument as in \cite[Corollary 2.2]{CyKaKrWa}, \cite[Theorem 1.1]{Su1}.
\begin{proposition}\label{cor;1}
Let $n \ge 1$, 
$
 \max \{1, {n}/2 \} \le p <\infty, 
$
and 
\[
 q_{\ast} 
 :=
 \min
 \left\{ 
  \frac{2np}{2n-p}, 2p
 \right\}
.
\]
Then, 
for all $A \in \mathbb{R}$ and all $v_0\in L^p(\mathbb{R}^n)$, 
there exist a number $T>0$ and a unique mild local-in-time solution $u$ 
to problem \eqref{eq;simplified-ar-ks} with the initial data 
$
 u_0
 =
 A
 +
 v_0 \in \mathcal{L}_{\mathrm{uloc}}^p(\mathbb{R}^n)
$ 
such that 
$
 u 
 \in 
 C([0,T);\mathcal{L}_{\mathrm{uloc}}^p(\mathbb{R}^n)) 
 \cap 
 C((0,T);\mathcal{L}_{\mathrm{uloc}}^q(\mathbb{R}^n))
$ 
and 
$
 u-A\in C([0,T); L^p(\mathbb{R}^n))
$. 
\end{proposition}

\subsection{Useful tools} We give three elementary facts which is often used in this paper.

\begin{lemma}\label{lem;1}
Let $N \in \mathbb{N}$. 
Then, for any $\alpha \in \mathbb{Z}^n_+$ and $\lambda > 0$, 
there exists a constant $C=C(\alpha,\lambda) > 0$ such that 
\begin{equation}\label{eq;lem1-1}
 \left| 
  \partial_\xi^\alpha
  \left(
   \frac1{\lambda + |\xi|^2}
  \right)
 \right|
 \le 
 \frac{C}{(1 + |\xi|)^{2+|\alpha|}}, \quad \xi \in \mathbb{R}^n.
\end{equation}
\end{lemma}

We apply the following identity for smooth functions 
in order to prove Lemma \ref{lem;1}.

\begin{lemma}[Multivariable Fa\`{a} di Bruno's formula]\label{lem;2}
Let 
$ 
 f 
 \in 
 C^\infty( \mathbb{R} \setminus \{ 0 \})
$ 
and 
$
 g 
 \in 
 C^\infty( \mathbb{R}^n \setminus \{ 0 \}) 
 \setminus 
 \{ 0 \}
$.
Then, for any $\alpha \in \mathbb{Z}^n_+$ satisfying $|\alpha| \ge 1$,
it holds that
\begin{equation}\label{eq;lem2-1}
 \partial_\xi^\alpha 
 f\big( g(\xi) \big)
 =
 \sum_{k=1}^{|\alpha|}
 \dfrac{\mathrm{d}^kf}{\mathrm{d}t^k}
 \big(
  g(\xi)
 \big)
 \sum_{\substack{\alpha_1 + \cdots + \alpha_{k} = \alpha \\ |\alpha_{i}| \ge 1}}
 \Gamma_{\alpha_1,\ldots, \alpha_{k}}^k
 \partial_\xi^{\alpha_1}g(\xi) \cdots \partial_\xi^{\alpha_{k}}g(\xi),
\end{equation}
where each $\Gamma_{\alpha_1,\ldots, \alpha_{k}}^k$ is a constant depending only on $k$ and 
$\alpha_1, \cdots, \alpha_{k} \in \mathbb{Z}^{n}_+$.
\end{lemma}

\begin{proof}[Proof of {\bf Lemma \ref{lem;1}}]
Set 
 $
  f(t)
  =
  1/({\lambda+t}) \, (t \ge 0)
 $ 
 and set 
 $
  g(\xi)
  =
  |\xi|^2 \, (\xi \in \mathbb{R}^n)
 $. 
Then, for $k \in \mathbb{N}$, we have
 $
  \frac{\mathrm{d}^{k}f}{\mathrm{d}t^k}(t)
  =
  {(-1)^{k}k!}/{(\lambda+t)^{k+1}}
 $. 
We note that for $l,m=1,2,\ldots,n$,
 $
  \partial_{\ell}|\xi|^2
  =
  2
  \xi_\ell
 $
 and
 $
  \partial_\ell \partial_{m}|\xi|^2
  =
  2
  \delta_{\ell m}
 $, where $\delta_{\ell m}$ is the Kronecker delta. 
Obviously, $\partial_{\xi}^\beta|\xi|^2=0$ for $\beta \in \mathbb{Z}^n_+$ with $|\beta| \ge 3$.
Thus, it holds that 
\[
 |
  \partial_\xi^\beta
  |\xi|^2
 |
 \le 
 2
 |\xi|^{2-|\beta|}
\]
for
$ 
 |\beta|
 \le 
 2
$
and 
$
 \xi
 \in 
 \mathbb{R}^n.
$
Taking $\alpha_1, \ldots, \alpha_{k} \in \mathbb{Z}^n_+$ with 
$
 \sharp \{ \alpha_j \, | \, |\alpha_j| = 1, j = 1, \cdots, k \} = m_{1}
$,
$
 \sharp \{ \alpha_j \, | \, |\alpha_j| = 2, j = 1, \cdots, k \} = m_{2}
$
and
$\alpha = \alpha_1 + \cdots + \alpha_{k}$
,
we have
$m_{1}+m_{2}=k$, $m_{1}+2m_{2}=|\alpha|$.
Then we obtain
\[
 \big|
 \partial_\xi^{\alpha_1}|\xi|^2
 \cdots
 \partial_\xi^{\alpha_{k}}|\xi|^2
 \big|
 \le 
 2^{k}
 |\xi|^{2-|\alpha_1|+2-|\alpha_2|+\cdots+2-|\alpha_{k}|}
 =
 2^{k}
 |\xi|^{2k-|\alpha|}.
\]
Applying formula \eqref{eq;lem2-1} leads to
\begin{align*}
 \left|
  \partial_{\xi}^\alpha
  \left(
  \frac1{\lambda+|\xi|^2}
 \right)
 \right|
 \le \, & 
 \sum_{k=1}^{|\alpha|}
 \bigg|
  \frac{\mathrm{d}^kf}{\mathrm{d}t^k}\big( g(\xi) \big)
 \bigg|
 \sum_{\substack{\alpha_1+\cdots+\alpha_k=\alpha, \\ |\alpha_i| \ge 1}}
 \big|
  \Gamma_{\alpha_1,\ldots, \alpha_k}^k
 \big|
 |
  \partial_{\xi}^{\alpha_1}g(\xi) 
  \cdots 
  \partial_{\xi}^{\alpha_k}g(\xi)
 |\\
 \le \, & 
 \sum_{\frac{|\alpha|}2 \le k \le |\alpha|}
 \sum_{\substack{\alpha_1+\cdots+\alpha_k=\alpha, \\ |\alpha_i|\ge 1}}
 2^k
 k!
 \big|
  \Gamma_{\alpha_1,\ldots, \alpha_k}^k
 \big|
 \frac1{(\lambda+|\xi|^2)^{k+1}}
 |\xi|^{2k-|\alpha|} \\
 = \, & 
 \sum_{\frac{|\alpha|}2 \le k \le |\alpha|}
 \sum_{\substack{\alpha_1+\cdots+\alpha_k=\alpha, \\ |\alpha_i|\ge 1}}
 2^k
 k!
 \big| 
  \Gamma_{\alpha_1,\ldots, \alpha_k}^k
 \big|
 \frac1{(\lambda+|\xi|^2)^{1+\frac{|\alpha|}2}}
 \left(
  \frac{|\xi|^{2}}{\lambda+|\xi|^2}
 \right)^{k-\frac{|\alpha|}2}\\
 \le \, &  
 \frac{C}{(\lambda+|\xi|^2)^{1+\frac{|\alpha|}2}}, \quad \xi\in {\mathbb R}^n
\end{align*}
since $m_1=2k-|\alpha|\ge 0$ and $m_2=|\alpha|-k\ge 0$.
Therefore,
noting that
\[
 \frac1{(\lambda+|\xi|^2)^{1+\frac{|\alpha|}2}}
 \le 
 C
 \begin{cases}
 \dfrac1{|\xi|^{2+|\alpha|}}, &|\xi|\ge 1, \\
 1, &|\xi|\le 1,
\end{cases}
\]
 we obtain estimate \eqref{eq;lem1-1}. 
\end{proof}
Next lemma is needed to show Theorem \ref{thm;1}.
\begin{lemma}\label{lem;ckkw-lem4.5}
Let $N\in {\mathbb N}$ be a positive integer with $N>n/2$. Then, for any $\alpha \in \mathbb{Z}_+^n$ satisfying 
$
 |\alpha|
=N, 
$ 
there exists a constant $C = C(n,N) > 0$ such that  
\begin{equation}\label{eq;es_v1}
 \|
  \varphi
 \|_1
 \le 
 C
 \|
  \mathcal{F}[\varphi]
 \|_2^{1-\frac{n}{2N}}
 \left(\sum_{|\alpha|=N}\|
  \partial_\xi^\alpha\mathcal{F}[\varphi]
 \|_2\right)^{\frac{n}{2N}}, \quad \varphi \in \mathcal{S}.
\end{equation}
\begin{proof}
Let $\varphi \in {\mathcal S}$ and let $N\in {\mathbb N}$ be a positive integer satisfying $N>{n}/{2}$. By the argument similar to that in \cite[Lemma 4.5]{CyKaKrWa}, we have 
\begin{equation}\label{es;l44_1}
 \|
  \varphi
 \|_1
  \le 
  C(n,N)
  \| 
   \varphi
  \|_2^{1-\frac{n}{2N}}
  \|
   |\cdot|^N
   \varphi
  \|_2^{\frac{n}{2N}}.
\end{equation}
Direct calculations also give 
\begin{equation}\label{es;l44_2}
\||\cdot|^N\varphi\|_2\le C(N)\sum_{|\alpha|=N}\|\partial_\xi^\alpha {\mathcal F}[\varphi]\|_2.
\end{equation}
Hence, combining $\|\varphi\|_2=\|{\mathcal F}[\varphi]\|_2$ and \eqref{es;l44_2} with \eqref{es;l44_1}, we obtain the desired estimate \eqref{eq;es_v1}.
\end{proof}
\end{lemma}
For $\lambda>0$, 
the Bessel kernel $B_\lambda(x)$ is defined as 
\begin{equation}\label{eq;Bessel-kernel}
 B_\lambda(x)
 :=
 (4\pi)^{-\frac{n}{2}}
 \int_0^\infty 
  e^{-\lambda\sigma-\frac{|x|^2}{4\sigma}}
  \sigma^{-\frac{n}{2}}
 \, \mathrm{d}\sigma, 
 \quad x \in \mathbb{R}^n.
\end{equation}
Then, 
the following lemma is to obtain $L^p({\mathbb R}^n)$ estimates on $\nabla B_\lambda$, which plays a crucial role in Section 8 (see \cite[Lemma 3.6]{CyKaKrWa}).
\begin{lemma}\label{lem;be}
Let 
\[
1\le p<\frac{n}{n-1} \quad \text{if} \quad n \ge 2
\quad \text{or} \quad 1\le p\le \infty \quad \text{if} \quad n=1.
\]
Then, $\nabla B_\lambda$ belongs to $L^p({\mathbb R}^n)$. Moreover, $\nabla B_\lambda\in L_w^{\frac{n}{n-1}}({\mathbb R}^n)$ if $n\ge 2$, where $L_{w}^p({\mathbb R}^n)$ is the weak Lebesgue spaces on ${\mathbb R}^n$.
\end{lemma}

\section{Linearized problem of equation \eqref{eq;simplified-ar-ks}}
For any $A \in \mathbb{R}$, 
$u_C \equiv A$
is constant steady states of problem \eqref{eq;simplified-ar-ks}. 
Then, 
by introducing the perturbation $v:= u - A$, 
we transform problem \eqref{eq;simplified-ar-ks} into the following problem:
\begin{equation}\label{eq;lin}\tag{Q}
\begin{cases} 
 \partial_t v - \Delta v+ A \Delta K*v 
 =
 -\nabla \cdot (v\nabla K*v), & t>0, \ x \in \mathbb{R}^n, \\
 v(0,x) = v_0(x),  &x \in \mathbb{R}^n,
\end{cases}
\end{equation}
where 
$
 K(x)
 :=
  \beta_1B_{\lambda_1}(x)
  -
  \beta_2B_{\lambda_2}(x)
$ 
and $v_0(x):=u_0(x)-A$. 
In this section, we treat the following {linear} problem \eqref{eq;lin-lin}:
\begin{equation}\label{eq;lin-lin}
\tag{L}
\begin{cases}
\partial_t w-L_Aw=0, &t>0, \ x \in \mathbb{R}^n, \\
w(0,x)=w_0(x), &x \in \mathbb{R}^n,
\end{cases}
\end{equation}
where the linearized operator $L_A$ is given by
\begin{equation}\label{eq;def-li-op}
 L_A
 \coloneqq 
 \Delta 
 -
 A
 \Delta K*
 .
\end{equation}
\subsection{Properties of the linearized operator $L_A$}
In this subsection, we investigate the linearized operator $L_A$ defined by 
\eqref{eq;def-li-op}. First of all, we begin with the following definition:
\begin{definition}\label{def;sy-cl}
Let $m_1 \in \mathbb{R}$ and $0 < m_2 \le 1$.
We define a function space $S_{m_2,0}^{m_1}$ as 
\[
 S_{m_2,0}^{m_1}
 \coloneqq 
 \{
  m \in C^{\infty}(\mathbb{R}^n) 
  \mid
  |\partial_\xi^\alpha m (\xi)|
  \le
  C_\alpha
  (1+|\xi|)^{m_1-m_2|\alpha|}
  \,
  \text{for all}
  \, 
  \xi \in \mathbb{R}^n
  \,
  \text{and each}
  \, 
  \alpha \in \mathbb{Z}_+^n
 \}.
\]
Moreover, for $m \in S_{m_2,0}^{m_1}$, we denote the multiplier $T_m$ on $L^p(\mathbb{R}^n) \ (1 \le p \le \infty)$ by 
\[
 (T_m f)(x)
 =
 \mathcal{F}^{-1}
 \big[
  m
  \mathcal{F}[f]
 \big](x),\quad
 f \in \mathcal{S}(\mathbb{R}^n).
\]
\end{definition}
Next, the operator $L_A$ given by relation \eqref{eq;def-li-op} is closable on $L^p(\mathbb{R}^n)$ for all $1\le p<\infty$, and 
its closure ${\mathcal L}_A$ can be obtained an exact form used the Fourier transform.
\begin{lemma}\label{prop;prop2.2}
Let $A \in \mathbb{R}$. 
Then,
the operator $L_A$ defined by relation \eqref{eq;def-li-op} is closable on $L^p(\mathbb{R}^n)$ for all $1\le p<\infty$. 
Also, the closure $\mathcal{L}_A$ of operator $L_A$ in $L^p(\mathbb{R}^n)$ has the following explicit form 
\begin{equation}\label{eq;def-li-op-cl}
 (\mathcal{L}_Af)(x)
 \coloneqq 
 \mathcal{F}^{-1}
 \big[
  h_A
  \mathcal{F}[f]
 \big](x), \quad f\in {\mathcal S},
\end{equation}
where the function $h_A(\xi)$ is given by 
\begin{equation}\label{eq;def-sy}
h_A(\xi)=-|\xi|^2
 +
 \beta_1A
 \frac{|\xi|^2}{\lambda_1+|\xi|^2}
 -
 \beta_2A
 \frac{|\xi|^2}{\lambda_2+|\xi|^2}.
\end{equation}
\end{lemma}
To show Lemma \ref{prop;prop2.2}, we estimate the function $h_A(\xi)$ given by \eqref{eq;def-sy}. 
\begin{lemma}\label{lem;der-sy}
For any $\alpha \in \mathbb{Z}_+^n$ and any $\xi \in \mathbb{R}^n$,
the function $h_A(\xi)$ defined by relation \eqref{eq;def-sy} satisfies
\begin{equation}\label{eq;der-sy-fo}
 \partial_x^\alpha 
 (
  h_A(\xi)+|\xi|^2
 )
 =
 -
 A\beta_1\lambda_1
 \partial_x^\alpha
 \left(
 \frac1{\lambda_1+|\xi|^2}
 \right)
 +
 A\beta_2\lambda_2
 \partial_\xi^\alpha
 \left(
  \frac1{\lambda_2+|\xi|^2}
 \right).
\end{equation}
Moreover, 
\begin{equation}\label{eq;der-1-2}
 | 
  \partial_{\xi}^\alpha h_A(\xi)
 |
 \le 
 C_1
 |\xi|^{2-|\alpha|}
 \quad
 (\xi \in \mathbb{R}^n)
\end{equation}
holds for any $|\alpha|=1, 2$ and some $C_1=C_1(A,\beta_1,\beta_2,\lambda_1,\lambda_2,\alpha)>0$.
Furthermore, for any $|\alpha| \ge 3$, there exists a constant $C_2=C_2(A,\beta_1,\beta_2,\lambda_1,\lambda_2,\alpha)>0$ such that
the derivative of function $h_A(\xi)$ satisfies
\begin{equation}\label{eq;der-3}
 | 
  \partial_{\xi}^\alpha h_A(\xi)
 |
 \le 
 C_2
 (1+|\xi|)^{2-|\alpha|}
 \quad (\xi \in \mathbb{R}^n).
\end{equation}
\end{lemma}
\begin{proof}
First, we have
\[
 A\beta_1
 \frac{|\xi|^2}{\lambda_1+|\xi|^2}
 -
 A\beta_2
 \frac{|\xi|^2}{\lambda_2+|\xi|^2}
 =
 A(\beta_1-\beta_2)
 -
 A\beta_1\lambda_1
 \frac{1}{\lambda_1+|\xi|^2}
 +
 A\beta_2\lambda_2
 \frac{1}{\lambda_2+|\xi|^2}
\]
for any $A \in \mathbb{R}$, 
$\beta_1, \beta_2, \lambda_1, \lambda_2>0$ and $\xi \in \mathbb{R}^n$, 
which gives relation \eqref{eq;der-sy-fo}. 
Next, we show estimate \eqref{eq;der-1-2}. 
Put 
\begin{equation}\label{fn;ga}
 g_A(\tau)
 \coloneqq
 -
 \tau
 +
 \beta_1
 A
 \frac{\tau}{\lambda_1+\tau}
 -
 \beta_2
 A
 \frac{\tau}{\lambda_2+\tau}, 
 \quad \tau \ge 0.
\end{equation}
Then, elementary calculations together give that there is a  constant $C=C(A, \beta_1, \beta_2, \lambda_1, \lambda_2)>0$ such that  
\[
 |
  g_A'(\tau)
 |
 +
 |
  \tau 
  g_A''(\tau)
 |
 \le 
 C
\]
for all $\tau \ge 0$. 
Therefore, noticing that $h_A(\xi)=g_A(|\xi|^2)$, we observe that 
\[
 |\partial_{\xi}^\alpha h_A(\xi)|=
 |2\xi^\alpha g_A'(|\xi|^2)|
 \le
 C |\xi|
\]
if $|\alpha|=1$, and
\[
 |\partial_{\xi}^\alpha h_A(\xi)|
 \le
 4
 |\xi^{\alpha}g_A''(|\xi|^2)|
 +
 2
 |g_A'(|\xi|^2)|
 \le
 C
\] 
if $|\alpha|=2$. 
Hence these yield estimate \eqref{eq;der-1-2}. 
We finally prove estimate \eqref{eq;der-3}.
It  follows from Lemma \ref{lem;1} that for any $\alpha \in \mathbb{Z}_+^n$ and any $\xi \in \mathbb{R}^n$, 
\begin{equation}\label{eq;lem1-2}
 \left|
  \partial_{\xi}^\alpha
  \left(
  \frac{\beta_1\lambda_1}{\lambda_1+|\xi|^2}
  \right)
 \right|+\left|
  \partial_{\xi}^\alpha
  \left(
  \frac{\beta_2\lambda_2}{\lambda_2+|\xi|^2}
  \right)
 \right|
 \le 
 \frac{C}{(1+|\xi|)^{2+|\alpha|}},
\end{equation}
where $C$ is a positive constant depending only on $\beta_1, \beta_2, \lambda_1,\lambda_2,\alpha$. 
Therefore,
by noting that $\partial_{\xi}^\alpha|\xi|^2=0$ for $\alpha \in \mathbb{Z}_+^n$ with $|\alpha| \ge 3$, the relations \eqref{eq;der-sy-fo} and \eqref{eq;lem1-2}
lead to
\[
 |\partial_\xi^\alpha h_A(\xi)|
 \le 
 C
 (1+|\xi|)^{-4}(1+|\xi|)^{2-|\alpha|}
 \le 
 C(1+|\xi|)^{2-|\alpha|}, \quad \xi \in \mathbb{R}^n,
\]
which implies our conclusion \eqref{eq;der-3}.
\end{proof}

\begin{proof}[Proof of {\bf Lemma \ref{prop;prop2.2}}]
By Lemma \ref{lem;der-sy}, we have 
$
 |
  \partial_\xi^\alpha h_A(\xi)
 |
 \le 
 C
 (1+|\xi|)^{2-|\alpha|}
$
for any $\alpha \in \mathbb{Z}_+^n$ and any $\xi \in \mathbb{R}^n$, 
which yields $h_A \in S_{1,0}^{2}$. 
As a result, Lemma \ref{prop;prop2.2} follows from the argument in \cite[Section 2, p.1390]{Wong}.
\end{proof}

The spectrum of  the operator ${\mathcal L}_A$ on $L^p({\mathbb R}^n)$ $(1<p<\infty)$ is consistent with the range of the function $h_A(\xi)$, and its essential spectrum is also the same.
\begin{lemma}\label{prop;prop2.22}
For each $A \in \mathbb{R}$ and every $1<p<\infty$, 
the operator $\mathcal{L}_A$ on $L^p(\mathbb{R}^n)$ only has essential spectrum $\sigma_{\mathrm{ess}}\big( \mathcal{L}_A \big)$ and satisfies
\begin{equation}\label{es;def-li-op-cl}
 \sigma \big( \mathcal{L}_A \big)
 =
 \sigma_{\mathrm{ess}}\big( \mathcal{L}_A \big)
 =
 \{
  h_A(\xi)
  \mid
  \xi \in \mathbb{R}^n
 \},
\end{equation}
where the function $h_A(\xi)$  is denoted by \eqref{eq;def-sy}.
\end{lemma}
We utilize the following proposition to prove Lemma \ref{prop;prop2.22}: 
\begin{proposition}[{\cite[Theorem 2.1]{Wong}}]\label{lem;l_W}
Let $n \ge 1$, $m_1>0$, $m_2\in (0,1]$ and $m \in S_{m_2,0}^{m_1}$
with
\[
 \frac{1}{m(\xi)}
 =
 O(|\xi|^{-b}), 
 \quad 
 |\xi|
 \to
 \infty
\]
for some $b>0$. Then, if there exists $p\in (1,\infty)$ such that
\[
 \left|
   \frac{1}{p}
   -
   \frac{1}{2}
 \right|
 <
 \frac{b}{n(m_1+1-m_2-b)},
\]
the essential spectrum of the closure $\mathcal{T}_{m,p}$ of operator $T_m$ in $L^{p}(\mathbb{R}^n)$ is given by
\begin{equation}\label{eq;op}
 \sigma(\mathcal{T}_{m,p})
 =
 \sigma_{\mathrm{ess}}(\mathcal{T}_{m,p})
 =
 \{ 
  m(\xi) 
  \mid 
  \xi \in \mathbb{R}^n 
 \}.
\end{equation}
Moreover, if 
\[
 b
 \ge 
 \frac{(m_1+1-m_2)n}{n+2},
\]
then formula \eqref{eq;op} holds for all $1 < p < \infty$. 
\end{proposition}
\begin{proof}[Proof of {\bf Lemma \ref{prop;prop2.22}}]
As mentioned in the proof of Lemma \ref{prop;prop2.2}, We recall $h_A\in \mathcal S_{1,0}^2$. Also, we can rewrite the function $h_A(\xi)$ as
\[
 h_A(\xi) 
 =
 |\xi|^2
 \left(
  -1
  +
  \frac{A\beta_1}{\lambda_1+|\xi|^2}
  -
  \frac{A\beta_2}{\lambda_2+|\xi|^2}
 \right),
\]
and thus
$1/h_A(\xi)=O(|\xi|^{-2})$ as $|\xi| \to \infty$. 
Therefore, making use of Proposition \ref{lem;l_W} and the inequality 
$
 2 
 \ge 
 2n/(n+2)
$ 
for any $n\ge 1$,
we obtain our conclusion \eqref{es;def-li-op-cl},  and the proof of Lemma \ref{prop;prop2.22} is complete.
\end{proof}

\section{$L^p$-$L^q$ estimate of $\mathrm{e}^{t\mathcal{L}_A}$}
For each $A \in \mathbb{R}$, 
it is shown that the operator $\mathcal{L}_A$ defined by relation \eqref{eq;def-li-op-cl} generates an analytic semigroup on 
$
 L^p(\mathbb{R}^n)
$
for all $1 \le p < \infty$ (see
\cite[Lemma 4.2]{CyKaKrWa}). In particular, the analytic semigroup, which is denoted by $\mathrm{e}^{t\mathcal{L}_A}$, is expressed by the following formula:
\[
 \mathrm{e}^{t\mathcal{L}_A}f
 \coloneqq 
 \mathcal{F}^{-1}[\mathrm{e}^{th_A}\mathcal{F}[f]]
 =
 \mathcal{K}_A(t)*{f}, 
 \quad f\in \mathcal{S},
\]
where the function $h_A$ is defined in relation \eqref{eq;def-sy} and 
the function $\mathcal{K}_A$ is given by 
\begin{equation}\label{eq;ker-def}
 \mathcal{K}_A(t,x)
 \coloneqq 
 \mathcal{F}^{-1}[\mathrm{e}^{th_A}](x).
\end{equation}

In this section we are going to give $L^p$-$L^q$ type estimates for the analytic semigroup $\mathrm{e}^{t\mathcal{L}_A}$. 
For this purpose, we begin with the following lemma:

\begin{lemma}\label{lem;es_m}
Let $A \in \mathbb{R}$ and $1 \le p\le \infty$.
Then, 
the linear operator 
\[
 \mathcal{M}_A 
 \coloneqq  
 -
 A\beta_1
 \Delta
 (\lambda_1-\Delta)^{-1}
 +
 A\beta_2 
 \Delta
 (\lambda_2-\Delta)^{-1}
\] 
is bounded on $L^p(\mathbb{R}^n)$, that is, 
we have
\begin{equation}\label{eq;estimate-M_A}
 \| 
  \mathcal{M}_A 
  f
 \|_p 
 \le
 2
 |A|
 \max 
 \{ 
  \beta_1, \beta_2 
 \}
 \| 
  f 
 \|_p
\end{equation}
for all $f \in L^p(\mathbb{R}^n)$. 
\end{lemma}

\begin{proof}[Proof of {\bf Lemma \ref{lem;es_m}}]
Since the proof of \eqref{eq;estimate-M_A} is more or less standard,
we skip to prove it (cf. \cite[Lemma 2. (i), pp.133]{Stein}). 
\end{proof}

We give a comment on the operator $\mathcal{M}_A$.
For all $f\in L^p(\mathbb{R}^n)$ $(1\le p\le \infty)$, 
the operator
$
 \mathrm{e}^{t{\mathcal M}_A}f
 \coloneqq
 \sum_{n=0}^\infty
  {(t{\mathcal M}_A)^nf}/{n!}
$
is the semigroup generated by the bounded operator ${\mathcal M}_A$ on $L^p(\mathbb{R}^n)$, 
and is represented as 
$
 \mathrm{e}^{t{\mathcal M}_A}f
 =
 \mathcal{F}^{-1}
 [
  \mathrm{e}^{t(|\cdot|^2+h_A)}
  \mathcal{F}[f]
 ]
$ 
by use of the Fourier transform. 
Also, the following estimate holds: 
\begin{equation}\label{eq;estimate-semigroup-M_A}
 \|\mathrm{e}^{t\mathcal{M}_A}f\|_p
 \le 
 \mathrm{e}^{2 |A| \max\{\beta_1, \beta_2\}t}\|f\|_p, \quad f\in L^p(\mathbb{R}^n).
\end{equation}

Now we observe that 
the analytic semigroup $\mathrm{e}^{t\mathcal{L}_A}$ admits the $L^p$-$L^q$ type estimates as follows:

\begin{lemma}\label{lem;L^p-L^q-any-A}
Assume that $A \in \mathbb{R}$ and $1 \le q \le p \le \infty$. 
Then there exists a constant $C>0$ such that 
\begin{align}
 &\label{eq;L^p-L^q-any-A}
  \|
   \mathrm{e}^{t\mathcal{L}_A}
   f
  \|_p
  \le 
  C
  t^{-\frac{n}2(\frac1q-\frac1p)}
  \mathrm{e}^{2 |A|\max \{ \beta_1, \beta_2 \}t}
  \|
   f
  \|_q, \\
 &\label{eq;L^p-L^q-derivative-any-A}
  \|
   \nabla 
   \mathrm{e}^{t\mathcal{L}_A}
   f
  \|_p
  \le 
   C
   t^{-\frac{n}2(\frac1q-\frac1p)-\frac12}
   \mathrm{e}^{2|A|\max \{ \beta_1, \beta_2 \}t}
   \|
    f
   \|_q
\end{align}
for all $f \in L^q(\mathbb{R}^n)$ and $t>0$.
\end{lemma}

\begin{proof}
It follows from the property of the operator $\mathcal{M}_A$
that 
\[
\mathrm{e}^{t\mathcal{L}_A}f
=
\mathcal{F}^{-1}
\big[
 \mathrm{e}^{th_A}
 \mathcal{F}[f]
\big]
=
\mathcal{F}^{-1}
\big[
\mathrm{e}^{-t|\cdot|^2}
\mathrm{e}^{t(|\cdot|^2+h_A)}
\mathcal{F}[f]
\big]
=
\mathrm{e}^{t\Delta}(\mathrm{e}^{t\mathcal{M}_A}f).
\]
Thus, the $L^p$-$L^q$ type estimate of the standard heat semigroup $\mathrm{e}^{t\Delta}$ and
estimate \eqref{eq;estimate-semigroup-M_A} lead to
\[
 \| 
  \mathrm{e}^{t\mathcal{L}_A}f 
 \|_p 
 \le 
 C
 t^{-\frac{n}2(\frac1q-\frac1p)}
 \|
  \mathrm{e}^{t\mathcal{M}_A}f
 \|_q
 \le 
 C
 t^{-\frac{n}2(\frac1q-\frac1p)}
 \mathrm{e}^{2|A|\max \{ \beta_1, \beta_2 \}t}
\|f\|_q.
\]
This implies estimate \eqref{eq;L^p-L^q-any-A}. 
Also, we use the argument similar to that in the proof of estimate \eqref{eq;L^p-L^q-any-A} to 
show estimate \eqref{eq;L^p-L^q-derivative-any-A}.
\end{proof}

\section{Decay estimates of $\mathrm{e}^{t\mathcal{L}_A}$}
When the constant $A \in \mathbb{R}$ satisfies condition \eqref{eq;cond-A}, 
we are going to prove that the analytic semigroup $\mathrm{e}^{t\mathcal{L}_A}$ admits the similar $L^p$-$L^q$ type estimate as that of 
the usual heat semigroup $\mathrm{e}^{t\Delta}$. 
As stated in Section 5, we recall that the analytic semigroup $e^{t{\mathcal L}_A}$ is represented as follows: For all 
$
 f \in \mathcal{S}(\mathbb{R}^n),
$
\[
 \mathrm{e}^{t{\mathcal L}_A}
 f
 =
 \mathcal{K}_A(t)*f, 
 \quad 
 \mathcal{K}_A(t)
 \coloneqq
 \mathcal{F}^{-1}[e^{th_A}],
\]
where the function $h_A$ is denoted by relation \eqref{eq;def-sy}. 
Thus, 
the key to show $L^p$-$L^q$ estimates of 
$
 \mathrm{e}^{t{\mathcal L}_A}
$ (see Proposition \ref{prop;lp-lq-li} below) are $L^2({\mathbb R}^n)$ estimates of $\partial_\xi^\alpha {\mathcal F}[{\mathcal K}_A(t)]$ for all $\alpha \in {\mathbb Z}^n$ 
and all $t>0$.
We first begin with the following lemma:
\begin{lemma}\label{lem;81}
Let $A_\ast$ be the same constant as in {relation} \eqref{eq;cond_A}. 
If the constant $A$ satisfies condition $0<A\le A_\ast$, 
then 
\begin{equation}\label{eq;h-up-bd}
 h_A(\xi)
 \le 
 -
 c_*
 |\xi|^2
\end{equation}
holds for all $\xi \in {\mathbb R}^n$. 
Here, the function $h_A(\xi)$ is given by relation \eqref{eq;def-sy}, 
and $c_\ast$ {stand for} a positive constant introduced as
\begin{equation}\label{eq;max--f}
 c_\ast:
 =
  \begin{cases}
    1, 
    &
    \dfrac{\beta_1}{\beta_2} \le 1,\quad \dfrac{\beta_1}{\beta_2} \le \dfrac{\lambda_1}{\lambda_2}, \qquad \\
    \dfrac{\beta_1\lambda_2 - \beta_2 \lambda_1}{\lambda_1 \lambda_2}
    \bigg( 
     \dfrac{\lambda_1 \lambda_2}{\beta_1\lambda_2 - \beta_2 \lambda_1} 
     - 
     A 
    \bigg) ,
    &
    \dfrac{\beta_1}{\beta_2} 
    > 
    \dfrac{\lambda_1}{\lambda_2},\quad 
    \dfrac{\beta_1}{\beta_2} 
    \ge 
    \bigg( 
     \dfrac{\lambda_1}{\lambda_2} 
    \bigg)^2, \qquad \\
    \dfrac{(\sqrt{\beta_1} - \sqrt{\beta_2})^2}{\lambda_1 - \lambda_2}
     \bigg( 
      \dfrac{\lambda_1 - \lambda_2}{(\sqrt{\beta_1} - \sqrt{\beta_2})^2} 
      - 
      A 
     \bigg), \qquad 
    & 
     \dfrac{\beta_1}{\beta_2} 
     < 
     \bigg( 
      \dfrac{\lambda_1}{\lambda_2} 
     \bigg)^2,\quad 
     \dfrac{\beta_1}{\beta_2} >1.
 \end{cases}
\end{equation}
\end{lemma}

\begin{proof}
Rewrite the function $h_A(\xi)=|\xi|^2f_A(|\xi|^2)$, 
where  
\[
 f_A(\tau) 
 \coloneqq  
 -1
 +
 \beta_1
 A 
 \dfrac{1}{\lambda_1 +\tau} 
 - 
 \beta_2
 A
 \dfrac{1}{\lambda_2+\tau}, 
 \quad \tau \ge 0.
\]
Hence, 
all you have to show Lemma \ref{lem;81} is to give the estimate 
\[
 f_A(\tau)
 \le 
 -c_{\ast}, 
 \quad \tau \ge 0.
\]
Let 
$
 \lambda
 \coloneqq
 {\lambda_1}/{\lambda_2}
$, 
$
 \beta
 \coloneqq
 {\beta_1}/{\beta_2}
$ 
and note that the derivative of $f_{A}(\tau)$ on $\tau\ge 0$ is given as
\begin{equation}\label{eq;dif-f}
 f_{A}'(\tau)
 =
 \frac{
   A
   \beta_2
   [
    (
     1
     +
     \sqrt{\beta}
    )
    \tau
    +
    \lambda_2
    (
     \lambda+\sqrt{\beta}
    )
   ]
  }{
   (
    \lambda_1
    +
    \tau
   )^2
   (
    \lambda_2
    +
    \tau
   )^2
   }
   [
    (
     1-\sqrt{\beta}
    )
    \tau
    -
    \lambda_2
    (
     \sqrt{\beta}
     -
     \lambda
    )
   ]
   .
\end{equation}
Now the proof is divided into the following three cases:
\begin{description}
\item[Case 1] $\beta \le 1$, $\beta \le \lambda$;
\item[Case 2] $\beta > \lambda$, $\beta \ge \lambda^2$;
\item[Case 3] $\beta < \lambda^2$, $\beta > 1$.
\end{description}

%

\begin{figure}[H]
\begin{center}
\begin{tikzpicture}[scale=2.0]
\fill[cyan, opacity=0.25] (0,0) --(3,0) --(3,1) --(1,1); 

\fill[red, opacity=0.5] plot[domain=1:sqrt(3), variable=\x, smooth] ({\x},{\x^2}) --(3,3) --(3,1) --(1,1); 

\fill[gray, opacity=0.6] plot[domain=1:sqrt(3), variable=\x, smooth] ({\x},{\x^2}) --(0,3) --(0,0) --(1,1);

 \draw[->,>=stealth,semithick] (-0.25,0)--(3,0) node [below]{{\footnotesize $\lambda$}}; 
 \draw[->,>=stealth,semithick] (0,-0.25)--(0,3) node [above, left]{{\footnotesize $\beta$}}; 
 \draw node[below right]{\textrm{\tiny O}}; 
 
 \draw[densely dotted, variable=\x, domain=0:0.25] plot(-\x,{\x^2}) ; 
 \draw[densely dotted, samples=100, variable=\x, domain=0:1] plot(\x,{\x^2}) ; 
 \draw[semithick, samples=100, domain=1:sqrt(3)+0.02, variable=\x] plot(\x,\x^2) node[above]{{\footnotesize $\beta = \lambda^2$}};

 \draw[densely dotted, domain=-0.25:1, variable=\x] plot(\x,1);
 \draw[semithick, domain=1:3.05, variable=\x] plot(\x,1) node[right]{{\footnotesize $\beta = 1$}};
 
  \draw[densely dotted, variable=\x, domain=0:3.125] plot(1,3-\x) node[below]{{\footnotesize $\lambda = 1$}};
 \draw[densely dotted, domain=1:3, variable=\x] plot(\x,\x) ;
 \draw[semithick, domain=0:1.125, variable=\x] plot(1-\x,1-\x) node[below left]{{\footnotesize $\beta = \lambda$}};

 \draw (0.7,2.0) node[scale=1.0] {\hbox {Case 2}} ; 
 \draw (2.2,1.8) node[scale=1.0] {\hbox {Case 3}} ; 
 \draw (1.8,0.5) node[scale=1.0] {\hbox {Case 1}};
\end{tikzpicture}
\caption{Positivity of $f_A(\tau)$}
\end{center}
\end{figure}
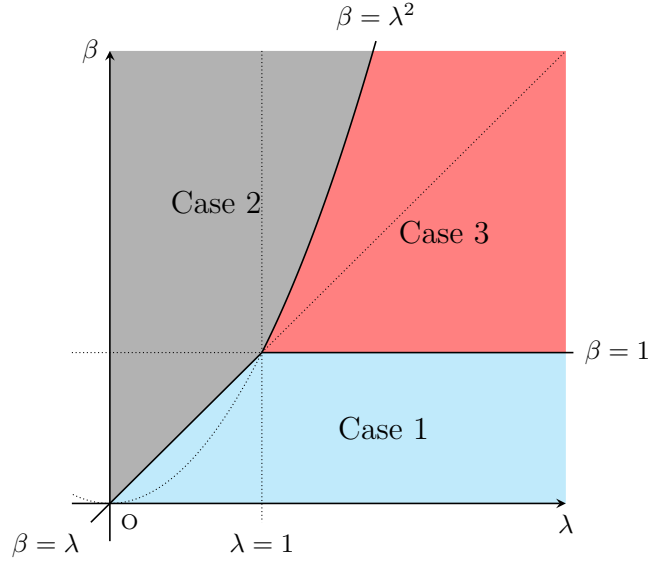

\noindent
\textbf{Case 1:} Remark that 
\begin{equation}\label{con;fai}
\lim_{\tau\to\infty}f_A(\tau)=-1
\end{equation}
for all $(\lambda, \beta)\in {\mathbb R}_{+}\times{\mathbb R}_{+}$, and 
\[
 f_{A}(0)
 =
 -1
 +
 A
 \frac{\beta_1\lambda_2-\beta_2\lambda_1}{\lambda_1\lambda_2}
 =
 -1
 -
 A
 \frac{\beta_2(\lambda-\beta)}{\lambda_1}
 \le 
 -1.
\]
In the case $\beta \le 1$, $\beta \le \lambda^2$, 
the derivative of $f_{A}(\tau)$ given by \eqref{eq;dif-f} is nonnegative for all $\tau \ge 0$. 
On the other hand, 
in the case $\beta \le \lambda$, $\beta > \lambda^2$, 
the function $f_A(\tau)$ on $\tau \ge 0$ attains a negative minimum at 
$
 \tau
 =
 \tau_0
$, where
\begin{equation}\label{no:tau0}
 \tau_0
 \coloneqq
 \frac{
  \lambda_2(\sqrt{\beta}-\lambda)
  }{
   1-\sqrt{\beta}
  }
 >
 0.
\end{equation}
As a result, we see that $f_A(\tau)\le -1$ for all $\tau\ge 0$.

\noindent
\textbf{Case 2:} 
Suppose that 
\begin{equation}\label{con;a3}
0<A<\frac{\lambda_1\lambda_2}{\beta_1\lambda_2-\beta_2\lambda_1}.
\end{equation}
Then, by \eqref{con;a3} we get 
\[
 -1
 <
 f_A(0)
 =
 -1
 +
 A
 \frac{\beta_1\lambda_2-\beta_2\lambda_1}{\lambda_1\lambda_2}
 <
 0.
\]
In the case $\beta\ge \lambda^2$, $\beta\ge 1$, the derivative of $f_A(\tau)$ given by \eqref{eq;dif-f} is not positive for all $\tau\ge 0$,  while in the case $\beta<1$, $\beta>\lambda$, the value $\tau_0$ denoted by \eqref{no:tau0} is positive, and it holds 
that the function $f_A(\tau)$ on $\tau\ge 0$ possesses a negative minimum at $\tau=\tau_0$. Thus it follows that
\[
 \sup_{\tau \ge 0}
 f_A(\tau)
 \le 
 f_A(0)
 =
 -1
 +
 A
 \frac{\beta_1\lambda_2-\beta_2\lambda_1}{\lambda_1\lambda_2}
 <
 0.
\]

\noindent
\textbf{Case 3:}  Assume that 
\begin{equation}\label{con;a2}
0<A<\frac{\lambda_1-\lambda_2}{(\sqrt{\beta_1}-\sqrt{\beta_2})^2}.
\end{equation}
Since the function $f_A(\tau)$ on $\tau\ge 0$ satisfies \eqref{con;fai} and the value $\tau_0$ given by \eqref{no:tau0} is positive, its function has a maximum at $\tau=\tau_0$, that is, 
\begin{equation}\label{max;fat0}
\sup_{\tau\ge 0}f_A(\tau)\le f_A(\tau_0)=\frac{(\sqrt{\beta_1}-\sqrt{\beta_2})^2}{\lambda_1-\lambda_2}\left(A-\frac{\lambda_1-\lambda_2}{(\sqrt{\beta_1}-\sqrt{\beta_2})^2}\right).
\end{equation}
Hence, \eqref{con;a2} and \eqref{max;fat0} imply that 
\[
 \sup_{\tau \ge 0}
 f_A(\tau)
 \le 
 \frac{
   (
    \sqrt{\beta_1}
    -
    \sqrt{\beta_2}
   )^2
  }{
   \lambda_1-\lambda_2
  }
 \left(
  A
  -
  \frac{
   \lambda_1-\lambda_2
   }{
    (
     \sqrt{\beta_1}
     -
     \sqrt{\beta_2}
    )^2
   }
 \right)
<0.
\]
Putting together with Cases 1, 2 and 3 implies our desired conclusion.

\end{proof}
\begin{remark} 
Assume that
\begin{equation}\label{c;ma1}
A=
\begin{cases}
 \dfrac{\lambda_1\lambda_2}{\beta_1 \lambda_2 - \beta_2 \lambda_1}, 
 &\dfrac{\beta_1}{\beta_2}>\dfrac{\lambda_1}{\lambda_2}, \ \dfrac{\beta_1}{\beta_2}\ge \left(\dfrac{\lambda_1}{\lambda_2}\right)^2, \\
 \dfrac{\lambda_1-\lambda_2}{(\sqrt{\beta_1}-\sqrt{\beta_2})^2}, 
 &\dfrac{\beta_1}{\beta_2}<\left(\dfrac{\lambda_1}{\lambda_2}\right)^2,\, \dfrac{\beta_1}{\beta_2}>1.
\end{cases}
\end{equation}
Then, for the function $h_A(\xi)$, 
we see that the asymptotic behavior in \eqref{c;ma1} is 
difference from that in \eqref{eq;cond-A} when $|\xi|\to0$ or $|\xi|\to \infty$. 
Indeed, note that $h_A(\xi)=|\xi|^2f_A(|\xi|^2)$ ($\xi \in {\mathbb R}^n$). 
If 
 $
  \beta_1/\beta_2
  >
  \lambda_1/\lambda_2
 $ 
and 
 $
  \beta_1/\beta_2
  \ge 
  (\lambda_1/\lambda_2)^2
 $, 
then we can write the function $h_A(\xi)$ with 
$
 A
 = 
  \lambda_1\lambda_2
  /
  (
  \beta_1 \lambda_2 
  - 
  \beta_2 \lambda_1
  )
$ 
by 
\[
 h_A(\xi)
 =
 -
  \dfrac{
   |\xi|^4
  }{
   ( 
    \lambda_1
    +
    |\xi|^2
   )
   (
    \lambda_2
    +
    |\xi|^2
   )
   }
  \bigg(
   |\xi|^2
   +
   \dfrac{
    \beta_2\lambda_1^2
    -
    \beta_1\lambda_2^2
    }{
     \beta_1\lambda_2
     -
     \beta_2\lambda_1
    }
  \bigg), 
\]
which behaves as follows:   
\[
 h_A(\xi)
 =
 \begin{cases}
 O(|\xi|^4), &\dfrac{\beta_1}{\beta_2}>\left(\dfrac{\lambda_1}{\lambda_2}\right)^2, \\
 O(|\xi|^6), &\dfrac{\beta_1}{\beta_2}=\left(\dfrac{\lambda_1}{\lambda_2}\right)^2,
 \end{cases}
\]
for $|\xi|\to 0$ and also $h_A(\xi)=O(1)$ for $|\xi|\to \infty$. 
On the other hand, 
if 
 $
  \beta_1/\beta_2
  <
  (\lambda_1/\lambda_2)^2
 $ 
and 
 $
  \beta_1/\beta_2>1
 $, 
the function $h_A(\xi)$
with 
 $
  A
  =
   (
   \lambda_1 
   -
   \lambda_2
   )
   /
   (
    \sqrt{\beta_1}
    -
    \sqrt{\beta_2}
   )^2
 $ 
can be represented as 
\[
 h_A(\xi)
 =
 -
  \dfrac{
   |\xi|^2
   }{
    (
     \lambda_1
     +
     |\xi|^2
    )
    (
     \lambda_2
     +
     |\xi|^2
    )
   }
  \bigg(
   |\xi|^2
   -
   \dfrac{
    \lambda_1\sqrt{\beta_2}
    -
    \lambda_2\sqrt{\beta_1}
    }{
     \sqrt{\beta_1}
     -
     \sqrt{\beta_2}
    }
  \bigg)^2,
\]
and have 
\[
 h_A(\xi)
  =
  \begin{cases}
  O(|\xi|^2), &|\xi|\to 0,\\
  O(|\xi|^2), &|\xi|\to \infty.
 \end{cases}
\]
Moreover it follows that
$
  h_A(\xi)
  =
  0
$
holds for all $\xi\in {\mathbb R}^n$ satisfying 
\[
 |\xi| 
 =
 \frac{\lambda_1\sqrt{\beta_2}-\lambda_2\sqrt{\beta_1}}{\sqrt{\beta_1}-\sqrt{\beta_2}}
 > 
 0
.
 \] 
As a result, for the analytic semigroup $\mathrm{e}^{t{\mathcal L}_A}$, 
it is suggested that the proof of the $L^p$-$L^q$ type estimates in condition \eqref{c;ma1} are extremely difficult in comparison with those in condition \eqref{eq;cond-A}. 
The details will be discussed in a forthcoming paper.
\end{remark}

The following lemma is the key estimate to establish decay estimates of $\mathrm{e}^{{\mathcal L}_A}$. 
\begin{lemma}\label{lem;es_A31}
Let $\alpha \in \mathbb{Z}_+^n$ with $|\alpha|\ge 1$. 
If the constant $A$ satisfies condition \eqref{eq;cond-A}, 
then
there exist constants $C>0$ and $c_*>0$ such that
\begin{equation}\label{eq;es_A31}
 |
  \partial_{\xi}^\alpha 
  \mathcal{F}[{\mathcal K}_A(t)](\xi)
 |
 \le \, 
 C
 \mathrm{e}^{-c_* t|\xi|^2}
  \sum_{\frac{|\alpha|}2 \le k \le |\alpha|}
  t^{k}
  |\xi|^{2k-|\alpha|}
\end{equation}
for all $t > 1$ and all $\xi \in \mathbb{R}^n$,
where the function $\mathcal{K}_A$ is defined by relation \eqref{eq;ker-def}, 
the constant $c_{*} > 0$ is given by \eqref{eq;max--f},
the constant
$C > 0$ depends only on $A$, $\beta_1$, $\beta_2$, 
$\lambda_1$, $\lambda_2$ and $\alpha$.
\end{lemma}

\begin{proof}
Notice that 
$
 \mathcal{F}[{\mathcal K}_A(t)](\xi)
 =
 \mathrm{e}^{th_A(\xi)}
$.
We set 
$
 f(z)
 =
 \mathrm{e}^{z} \ (z \ge 0)
$ 
and  
$
 g_t(\xi)=th_A(\xi)
$
for
$t \ge 0, \xi \in \mathbb{R}^n$.
Applying formula \eqref{eq;lem2-1} to the functions $f$ and $g_t$, 
we have 
\[
 \begin{split}
 \partial_\xi^{\alpha} 
 {{\mathcal F}[{\mathcal K}_A(t)](\xi)}
 = \, &
 \sum_{k=1}^{|\alpha|}
 f^{(k)}(g_t(\xi))
 \sum_{\substack{\alpha_1+\cdots+\alpha_k=\alpha, \\ |\alpha_i|\ge 1}} 
 \Gamma_{\alpha_1,\ldots, \alpha_k}^k
 (\partial_\xi^{\alpha_1}g_t(\xi))\cdots(\partial_\xi^{\alpha_k}g_t(\xi)) \\
 = \, &
 \mathrm{e}^{th_A(\xi)}
 \sum_{k=1}^{|\alpha|}
 \sum_{\substack{\alpha_1+\cdots+\alpha_k=\alpha, \\ |\alpha_i|\ge 1}}
 \Gamma_{\alpha_1,\ldots, \alpha_k}^k
 t^{k}
 (\partial_\xi^{\alpha_1}h_A(\xi))
 \cdots
 (\partial_\xi^{\alpha_k}h_A(\xi))
 \end{split}
\]
for all $\xi \in \mathbb{R}^n$ {and all $\alpha \in \mathbb{Z}_+^n$ with $|\alpha|\ge 1$}.
Combining {estimates} \eqref{eq;der-1-2}{,} \eqref{eq;der-3} {and \eqref{eq;h-up-bd},} 
we find that
\[
\begin{split}
 |\partial_\xi^{\alpha}  {{\mathcal F}[{\mathcal K}_A(t)](\xi)}| 
 \le \, & 
 \mathrm{e}^{th_A(\xi)}
 \sum_{1\le k<\frac{|\alpha|}2}
 \sum_{\substack{\alpha_1+\cdots+\alpha_k=\alpha, \\ |\alpha_i|\ge 1}}
 \left|
  \Gamma_{\alpha_1, \ldots, \alpha_k}^k
 \right|
 t^{k}
 \left| 
  (\partial_\xi^{\alpha_1}h_A(\xi))
  \cdots
  (\partial_\xi^{\alpha_k}h_A(\xi))
 \right| \\
 & \, 
 +
 \mathrm{e}^{th_A(\xi)}
 \sum_{\frac{|\alpha|}2\le k \le |\alpha|}
 \sum_{\substack{\alpha_1+\cdots+\alpha_k=\alpha, \\ |\alpha_i|\ge 1}}
 \left|
  \Gamma_{\alpha_1,\ldots, \alpha_k}^k
 \right|
 t^{k}
 \left|
  (\partial_\xi^{\alpha_1}h_A(\xi))
  \cdots
  (\partial_\xi^{\alpha_k}h_A(\xi))
 \right| \\
 \le \, &
 t^{\frac{|\alpha|}2}
 \mathrm{e}^{th_A(\xi)}
 \sum_{1 \le k < \frac{|\alpha|}2}
 C_k^1
 t^{k-\frac{|\alpha|}2}
 (1+|\xi|)^{2k-|\alpha|}
 +
 \mathrm{e}^{th_A(\xi)}
 \sum_{\frac{|\alpha|}2 \le k \le |\alpha|}
 C_k^2
 t^{k}
 |\xi|^{2k-|\alpha|}\\
 \le \, &
 t^{\frac{|\alpha|}2}
 \mathrm{e}^{-c_* t|\xi|^2}
 \sum_{1 \le k < \frac{|\alpha|}2}
 C_k^1
 +
 \mathrm{e}^{-c_* t|\xi|^2}
 \left(
  \max_{\frac{|\alpha|}2 \le k \le |\alpha|}C_k^2
 \right)
 \sum_{\frac{|\alpha|}2 \le k \le |\alpha|}
 t^{k}
 |\xi|^{2k-|\alpha|}\\
 \le \, &
 C
 \mathrm{e}^{-c_* t|\xi|^2}
  \sum_{\frac{|\alpha|}2 \le k \le |\alpha|}
  t^{k}
  |\xi|^{2k-|\alpha|},
\end{split}
\]
where $t > 1$ and  
\[
\begin{split}
 C
 \coloneqq  \, 
 \max
 \left[
  \sum_{1 \le k < \frac{|\alpha|}2}C_k^1, 
  \max_{\frac{|\alpha|}2 \le k \le |\alpha|}C_k^2
 \right], \quad 
 C_k^j
 \coloneqq  \, 
 \sum_{\substack{\alpha_1 + \cdots + \alpha_k = \alpha, \\ |\alpha_i|\ge 1}} 
 C_{\alpha_1, \alpha_2, \ldots, \alpha_k}^{{j}}
 \left| 
  \Gamma_{\alpha_1, \ldots, \alpha_k}^k
 \right|{,} \quad {j=,1,2}.
\end{split}
\]
Hence we obtain estimate \eqref{eq;es_A31}. 

\end{proof}

By virtue of estimate \eqref{eq;es_A31}, {we obtain the following lemma}.
\begin{lemma}\label{lem;esu_ma}
Let {the constant $A$} satisfy condition \eqref{eq;cond-A}.
For all $\alpha \in \mathbb{Z}_+^n$, 
there exists {a constant} $C=C(n, \beta_1, \beta_2, \lambda_1, \lambda_2, A, \alpha)>0$ such that
\begin{equation}\label{eq;esu_ma}
 \|
  \partial_\xi^\alpha
   {{\mathcal F}[{\mathcal K}_A(t)](\xi)}
 \|_2
 \le 
 C
 t^{-\frac{n}4+\frac{|\alpha|}2}, \quad t {>} 1.
\end{equation}
\end{lemma}

\begin{proof}
When $|\alpha|=0$, inequality \eqref{eq;h-up-bd} shows {that}
\[
 \|
   {{\mathcal F}[{\mathcal K}_A(t)](\xi)}
 \|_2 
 \le 
  \left(
   \int_{\mathbb{R}^n}
    \mathrm{e}^{-2c_{*}t|\xi|^2}
   \,
   \mathrm{d}\xi
  \right)^{\frac12} \\
 = 
  \left(
   \int_{\mathbb{R}^n}
    \mathrm{e}^{-|z|^2}(2c_{*}t)^{-\frac{n}2}
   \, \mathrm{d}z
  \right)^{\frac12}
  = 
  C
  t^{-\frac{n}4}
\]
for all $t>0$. {Moreover,} for $|\alpha| \ge 1$, inequality \eqref{eq;es_A31} {in} Lemma {\ref{lem;es_A31}} gives
\[
\begin{split}
  \|
   \partial_\xi^\alpha
    {{\mathcal F}[{\mathcal K}_A(t)](\xi)}
  \|_2^2
  \le \, &  
   C
   \int_{\mathbb{R}^n} 
   \mathrm{e}^{-2c_* t|\xi|^2}
    \left[
     \sum_{\frac{|\alpha|}2 \le k \le |\alpha|}
     t^{k}
     |\xi|^{2k-|\alpha|}
    \right]^2
   \, \mathrm{d}\xi \\
   = \, &
   C
   t^{|\alpha|}
   \int_{\mathbb{R}^n} 
    \mathrm{e}^{-|\xi\sqrt{2c_* t}|^2}
     \left[
      \sum_{\frac{|\alpha|}2 \le k \le |\alpha|} 
      (2c_*)^{-k+\frac{|\alpha|}2} 
      |
       \xi \sqrt{2c_* t}
      |^{2k-|\alpha|}
     \right]^2
   \, \mathrm{d}\xi \\
   = \, &
   Ct^{-\frac{n}2+|\alpha|}
   \int_{\mathbb{R}^n} 
    \mathrm{e}^{-|z|^2}
     \left[
      \sum_{{\frac{|\alpha|}{2}}\le k\le {|\alpha|}}
      (2c_*)^{-k-\frac{n}4+\frac{|\alpha|}2}
      |z|^{2k-|\alpha|}
    \right]^2
   \, \mathrm{d}z\\
   = \, & 
   Ct^{-\frac{n}2+|\alpha|}.
\end{split}
\]
Therefore, estimate \eqref{eq;esu_ma} holds for all $t {>} 1$ and all $\alpha \in \mathbb{Z}_+^n$, 
and the proof of Lemma \ref{lem;esu_ma} is complete.
\end{proof}

As we mentioned before, Lemma \ref{lem;ckkw-lem4.5} leads to a uniform bound of the norm 
$
 \|
  \mathcal{K}_A(t)
 \|_1{,}
$ 
where the function $\mathcal{K}_A(t)$ is denoted by {relation} \eqref{eq;ker-def}.

\begin{lemma}\label{lem;lem2.8}
Suppose that the constant {$A$} satisfy condition \eqref{eq;cond-A}.
Then there exists {a constant} $C > 0$ depending on $\beta_1, \beta_2, \lambda_1, \lambda_2, n, A$ such that
\begin{equation}\label{es_ma1}
 \|
 { {\mathcal K}_A(t)}
 \|_1
 \le 
 C, \quad 
 t {>} 1.
\end{equation}
\end{lemma}

\begin{proof}
{Let $N$ be a positive integer with $ N > {n}/{2}$. Then, 
for any $\alpha \in \mathbb{Z}_+^n$ with $|\alpha|=N$ and any $t>1$}, {noting \eqref{eq;ker-def} and} combining estimate \eqref{eq;esu_ma} with estimate \eqref{eq;es_v1},
we have
\[
\begin{split}
 \|
 { {\mathcal K}_A(t)}
 \|_1
 = \, &\|
  \mathcal{F}^{-1}
  [
   \mathrm{e}^{th_A}
  ]
 \|_1\\
 \le \, &
 C
 \|
  {e^{th_A}}
 \|_2^{1-\frac{n}{2N}}
 \left(\sum_{|\alpha|=N}\|
  \partial_\xi^\alpha 
  {e^{th_A}}
 \|_2\right)^{\frac{n}{2N}}
 \le 
 C 
 (t^{-\frac{n}4})^{1-\frac{n}{2N}}
 (t^{-\frac{n}4+\frac{N}2})^{\frac{n}{2N}}
 \le 
 C,
\end{split}
\]
which shows estimate \eqref{es_ma1}. 
\end{proof}

Together with Lemma \ref{lem;L^p-L^q-any-A} and Lemma \ref{lem;lem2.8},
we obtain the $L^p$-$L^q$ type estimate of the analytic semigroup $\mathrm{e}^{t\mathcal{L}_A}$ 
when the constant $A>0$ satisfies condition \eqref{eq;cond-A}.

\begin{proposition}\label{prop;lp-lq-li}
Let $1 \le q \le p \le \infty$.
Assume that the constant $A \in \mathbb{R}$ satisfies condition \eqref{eq;cond-A},
then there exists {a constant} $C=C(p,q,n,A,\beta_1,\beta_2,\lambda_1,\lambda_2)>0$ such that 
\begin{align}
 &\label{eq;es_sa} 
 \| 
  \mathrm{e}^{t\mathcal{L}_A}f 
 \|_p
 \le 
 C
 t^{-\frac{n}2(\frac1q-\frac1p)}
 \|
  f
 \|_q, \\
 &\label{eq;es_nsa}
 \|
  \nabla \mathrm{e}^{t\mathcal{L}_A}
  f
 \|_p
 \le 
 C
 t^{-\frac{n}2(\frac1q-\frac1p)-\frac12}
 \|
  f
 \|_q
\end{align}
for all $t>0$ and all $f \in L^q(\mathbb{R}^n)$.
\end{proposition}

\begin{proof}
We first show that estimate \eqref{eq;es_sa} holds. 
For {any} $1 \le {p} \le \infty$, 
estimate \eqref{eq;L^p-L^q-any-A} 
with {$p=q$} lead{s} to
\[
 \|
  \mathrm{e}^{t\mathcal{L}_A}
  f
 \|_{{p}}
 \le
 C
 \|
  f
 \|_{{p}}
 , \quad 0 < t \le 1.
\]
On the other hand, 
the Hausdroff-Young inequality and estimate \eqref{es_ma1} give
\[
 \|
  \mathrm{e}^{t\mathcal{L}_A}
  f
 \|_{{p}}
 =
 \|
  {{\mathcal K}_A(t)}
  *
  f
 \|_{{p}}
 \le 
 \|
  {{\mathcal K}_A(t)}
 \|_1
 \|
  f
 \|_{{p}}
 \le 
 C
 \|
  f
 \|_{{p}}
\]
for all $t > 1$.
Thus, from these estimates we obtain estimate \eqref{eq;es_sa} in the case $p=q$ for all $t > 0$. 
Now we introduce the auxiliary parameters $0 < \varepsilon^\ast \le 1/2$ and $A^\ast > 0$ by
\begin{equation*}
 {\varepsilon^{\ast}}
 =
 \begin{cases}
   \displaystyle \frac{1}{2},\quad & \displaystyle {\frac{\beta_1}{\beta_2}\le 1,\quad \frac{\beta_1}{\beta_2}\le \frac{\lambda_1}{\lambda_2}}, \vspace{2mm} \\ 
   \displaystyle\frac12\left(1  - \frac{A}{A_\ast}\right), &\displaystyle {\frac{\beta_1}{\beta_2} > \frac{\lambda_1}{\lambda_2}, \quad \frac{\beta_1}{\beta_2}\ge \left(\frac{\lambda_1}{\lambda_2}\right)^2} \quad \text{or} \quad {\frac{\beta_1}{\beta_2}<\left(\frac{\lambda_1}{\lambda_2}\right)^2,\quad \frac{\beta_1}{\beta_2}>1 }
 \end{cases}
\end{equation*}
where 
\begin{equation}\label{eq;cond-A-th}
 {A^{\ast}} 
 \coloneqq 
 \left\{
 \begin{split}
  & \, 
  \frac{\lambda_1 \lambda_2}{\beta_1\lambda_2 - \beta_2\lambda_1},\quad 
  && {\frac{\beta_1}{\beta_2} > \frac{\lambda_1}{\lambda_2}, \quad \frac{\beta_1}{\beta_2}\ge \left(\frac{\lambda_1}{\lambda_2}\right)^2}, \\
  & \, \frac{\lambda_1-\lambda_2}{(\sqrt{\beta_1}-\sqrt{\beta_2})^2},\quad 
  && {\frac{\beta_1}{\beta_2}<\left(\frac{\lambda_1}{\lambda_2}\right)^2},\quad   \frac{\beta_1}{\beta_2}>1.
 \end{split}
 \right.
\end{equation}
Combining the relation
\[
 \mathrm{e}^{t\mathcal{L}_A}f
 =
 \mathcal{F}^{-1}
 \Big[
  \mathrm{e}^{-\varepsilon_{\ast}t|{\cdot}|^2}
  \mathrm{e}^{(1-\varepsilon_{\ast})t h_{\frac{A}{1-\varepsilon_{\ast}}}}
  \mathcal{F}[f]
 \Big]
 =
 \mathrm{e}^{\varepsilon_{\ast}t\Delta}
 \Big(
  \mathrm{e}^{(1-\varepsilon_{\ast})t\mathcal{L}_{\frac{A}{1-\varepsilon_{\ast}}}}
  f
 \Big),
\]
the bound 
$
 0
 <
 {A}/({1-\varepsilon_\ast})
 <
 {A^\ast}
$,
the usual $L^p$-$L^q$ type estimate of $\mathrm{e}^{\varepsilon_* t\Delta}$
and {estimate \eqref{eq;es_sa} in the case $p=q$},
we obtain
\[
\begin{split}
 \|
  \mathrm{e}^{t\mathcal{L}_A}
  f
 \|_p
 \le 
 C
 t^{-\frac{n}2(\frac1q-\frac1p)}
 \|
  \mathrm{e}^{(1-\varepsilon_{\ast})t\mathcal{L}_{\frac{A}{1-\varepsilon_{\ast}}}}
  f
 \|_q
 \le 
 C
 t^{-\frac{n}2(\frac1q-\frac1p)}
 \|
  f
 \|_q
\end{split}
\]
for all $1 \le q \le p \le \infty$ and all $t > 0$,
which is our desired conclusion \eqref{eq;es_sa}.  
Furthermore we apply 
the relation 
\[
 \nabla 
 \mathrm{e}^{t\mathcal{L}_A}
 f
 =
 \nabla 
 \mathrm{e}^{\varepsilon_{\ast}t\Delta}
 \Big(
  \mathrm{e}^{(1-\varepsilon_{\ast})t\mathcal{L}_{\frac{A}{1-\varepsilon_{\ast}}}}
  f
 \Big){,}
\]
the usual
$L^p$-$L^q$ type estimate of $\N \mathrm{e}^{\varepsilon_\ast t \Delta}$ and estimate \eqref{eq;es_sa} with $p=q$ in order to prove estimate \eqref{eq;es_nsa}. Hence we complete the proof of Proposition \ref{prop;lp-lq-li}.
\end{proof}

\section{Exponential growth of $\mathrm{e}^{t\mathcal{L}_A}$}
In this section, we {give $L^p$-$L^q$ estimates of the analytic semigroup $\rm{e}^{t{\mathcal L}_A}$}
{when the constant $A>0$ is sufficiently large under the condition}
\[
 \frac{\beta_1}{\beta_2}
 > 
 \frac{\lambda_1}{\lambda_2}
 \quad
 \text{or}
 \quad
 \frac{\beta_1}{\beta_2}
 >
 1.
\]

To this end, it is essential to observe detail properties of the operator $\mathcal{L}_A$. 
According to Lemma \ref{prop;prop2.2}, the operator $\mathcal{L}_A$ has exactly the form 
$
\mathcal{L}_A\varphi
 =\mathcal{F}^{-1}[
 h_A
 \mathcal{F}[\varphi]]
$
for any $\varphi \in \mathcal{S}$. Moreover, it follows from Lemma \ref{prop;prop2.22} that 
the operator $\mathcal{L}_A$ only admits the essential spectrum 
$
 \sigma_{\mathrm{ess}}
 \big(
  \mathcal{L}_A
 \big)
$, 
which lies in $\mathbb{R}$ (see Proposition \ref{prop;ess-sp-li} below).
In particular, 
the maximum value $\max_{\xi \in \mathbb{R}^n} h_A(\xi)$ plays the resonance of 
$
 \sigma_{\mathrm{ess}}
 \big(
  \mathcal{L}_A
 \big)
$. 
\begin{proposition}\label{prop;ess-sp-li}
Let $A > 0$ and $1 < p < \infty$. Then, for the operator $\mathcal{L}_A$ on $L^p(\mathbb{R}^n)$ given by \eqref{eq;def-li-op-cl}, 
the following holds:
\[
 \sigma
 \big(
  \mathcal{L}_A
 \big)
 =
 \sigma_{\mathrm{ess}}
 \big(
  \mathcal{L}_A
 \big)
 =
 \big(
  -\infty,
  M(A)
 \big]{,}
\]
where 
\[
 M(A)
 =
 \max_{\xi \in \mathbb{R}^n}
 h_A(\xi)
 =
 \max_{\tau \ge 0}
 g_A(\tau),
\]
and 
$h_A(\xi)$ {and} $g_A(\tau)$ are defined by relations \eqref{eq;def-sy} {and} \eqref{fn;ga}{, respectively}.
Furthermore, $M(A) = 0$ holds if and only {if} the constant {$A$} satisfies condition {$0<A\le A_\ast$}. {Here the constant $A_\ast$ is defined as \eqref{eq;cond_A}.}
\end{proposition}
\begin{proof}
It follows from \eqref{es;def-li-op-cl} in Lemma \ref{prop;prop2.22} and $h_A(\xi)=g_A(|\xi|^2)$ $(\xi\in{\mathbb R}^n)$ that 
\[
 \sigma (\mathcal{L}_A)
 =
 \sigma_{\text{ess}}(\mathcal{L}_A)
 =
 \{ g_{A}(\tau) \mid \tau\ge 0 \}.
\]
Thus, in order to show Proposition \ref{prop;ess-sp-li}, 
we need to derive the maximum value of $g_A$ since the function $g_A(\tau)$ tends to $-\infty$ as $\tau\to\infty$.

Calculating the derivative of $g_A$ up to second  order with respect to $\tau$ yields that
\[
 g_{A}'(\tau)
 =
  -1
  +
  \frac{
   A\beta_1\lambda_1
   }{
    (
     \lambda_1+\tau
    )^2
    }
  -
  \frac{
   A\beta_2\lambda_2
   }{
    (
     \lambda_2+\tau
    )^2
    }
 \]
 and
 \begin{equation}\label{n;ga2}
 g_{A}''(\tau)
 = 
 -
 \frac{
  2A\beta_2^{\frac{1}{3}}
  \lambda_2^{\frac{1}{3}}
  I(\tau)
   }{
    (
     \lambda_1+\tau
    )^3
    (
     \lambda_2+\tau
    )^3
    }
    [
     (
      \beta^{\frac{1}{3}}
      \lambda^{\frac{1}{3}}
      -1
     )
     \tau
     +
     \lambda_2
     \lambda^{\frac{1}{3}}
     (
      \beta^{\frac{1}{3}}
      -
      \lambda^{\frac{2}{3}}
     )
    ], 
\end{equation} 
respectively, where $\beta \coloneqq \beta_1/\beta_2, \lambda \coloneqq \lambda_1/\lambda_2$ and
\[
 I(\tau)
 \coloneqq
 \left[ 
  \beta_2^{\frac{1}{3}}
  \lambda_2^{\frac{1}{3}}
  (\lambda_1+\tau)
  +
  \frac{\beta_1^{\frac{1}{3}}
  \lambda_1^{\frac{1}{3}}
  (
   \lambda_2+\tau
  )
  }{2}
 \right]^2
 +
 \frac{3}{4}
  \beta_1^{\frac{2}{3}}
  \lambda_1^{\frac{2}{3}}
  (
   \lambda_2+\tau
  )^2
  >
  0.
\]
Then the proof as in Lemma \ref{lem;81} is divided into the following three cases:
\begin{description}
\item[Case 1] the constant $A$ satisfies condition $0<A\le A_\ast$;
\item[Case 2] the constant $A$ satisfies condition
\[
 A
 >
 \frac{
  \lambda_1\lambda_2
  }{
   \beta_1\lambda_2
   -
   \beta_2\lambda_1
  } 
  \quad 
  (
   \beta \ge \lambda^2, \ 
   \beta>\lambda
  );
\]
\item[Case 3] the constant $A>0$ satisfies condition
\[
 A
 >
 \frac{
  \lambda_1
  -
  \lambda_2
  }{
   (
    \sqrt{\beta_1}
    -
    \sqrt{\beta_2}
   )^2
  }
  \quad 
  ( \beta > 1 , \ \beta < \lambda^2 ).
\]
\end{description}

\noindent
{{\bf Case 1.} 
By \eqref{eq;h-up-bd} and $h_A(\xi)=g_A(|\xi|^2)$, we have $g_A(\tau)\le 0$ for any $\tau\ge 0$, which together with $g_{A}(0)=0$ yields $\max_{\tau\ge 0}g_A(\tau)=0$. 
} \vspace{2mm}

\noindent
{{\bf Case 2.} 
It follows from \eqref{n;ga2} that the function $g_A'(\tau)$ on $\tau\ge 0$ is monotone decreasing if $\beta\lambda\ge 1$, and has a local minimum at $\tau=\tau_0$ and $g_A'(\tau_0)<-1$
 if $\beta\lambda<1$, where 
\begin{equation}\label{n;tau0}
\tau_0:=\frac{\lambda_2\lambda^{\frac{1}{3}}(\beta^{\frac{1}{3}}-\lambda^{\frac{2}{3}})}{1-(\beta\lambda)^{\frac{1}{3}}}>0.
\end{equation}
Hence, by these facts, $g_A'(0)>0$ and $\lim_{\tau\to\infty}g_A'(\tau)=-1$, there is a constant $\tau_1>0$ such that the function $g_A(\tau)$ on $\tau\ge 0$ 
has a local maximum at $\tau=\tau_1$. As a consequence, by remarking that $g_A(0)=0$ and $\lim_{\tau\to\infty}g_A(\tau)=-\infty$, it holds that \[\max_{\tau\ge 0}g_A(\tau)=g_A(\tau_1)>0.\]}

\noindent
{{\bf Case 3.} Note that 
\begin{equation}\label{c;a3}
 \frac{
  \lambda_2
  (
   \lambda
   -
   1
  )^2
  }{
   \beta_2
   [
    (\beta \lambda)^{\frac{1}{3}}
    -
    1
   ]^3
   }
 <
 \frac{\lambda_1-\lambda_2}{(\sqrt{\beta_1}-\sqrt{\beta_2})^2}.
\end{equation}
From $\beta\lambda>1$ and \eqref{n;ga2}, we see that the function 
$g_A'(\tau)$ on $\tau\ge 0$ admits a local maximum at $\tau=\tau_0$, where 
$\tau_0$ is denoted by \eqref{n;tau0}.
Hence, by $g_A'(\tau_0)>0$ and $\lim_{\tau\to \infty}g_A'(\tau)=-1$, 
there exists $\tau_1$ with 
$
 0
 <
 \tau_0
 <
 \tau_1
$ 
such that the function $g_A(\tau)$ on $\tau \ge 0$ has a local maximum at $\tau = \tau_1$. 
Here $\tau_1$ is a positive real root of the algebraic equation $g_A'(\tau)=0$. 
From $g_A(0)=0$ and 
$
 \lim_{\tau \to \infty}
 g_A(\tau)
 =
 -\infty
$, we observe that 
\[
 \max_{\tau \ge 0}
 g_A(\tau)
 =
 g_A(\tau_1)
 >
 0.
\]
}  

{A a result, combining Cases 1, 2 and 3, we obtain Proposition \ref{prop;ess-sp-li}.}
\end{proof}
{We give the following two remarks on Proposition \ref{prop;ess-sp-li}. 
For detail, see Appendix A.
\begin{remark}
\begin{enumerate}
\item[] 
\item Ferrari's formula enables us to give the specific form of a point at which the function $g_A$ takes a local maximum.
\item The condition
\begin{equation}\label{eq;ar-ba-co}
 \beta_1\lambda_1 - \beta_2\lambda_2 
 \gtreqqless 
 0
\end{equation}
or equivalently 
\[
 \frac{\beta_1}{\beta_2}\frac{\lambda_1}{\lambda_2} 
 \gtreqqless 
 1
\]
is one of the most important conditions introduced in \cite{HoOg2}, \cite{LiShWa}, \cite{NaYa1},  
and plays a essential role in the monotonicity of function $g_A$ denoted by \eqref{fn;ga} in our problem. 
Indeed, it is shown that the condition $g_A'(\tau) < 0$ holds for all $\tau \ge 0$ if and only if 
the constant $A>0$ satisifes
 \begin{equation}\label{eq;co-A-mo}
\begin{cases}
A> 0, &\dfrac{\beta_1}{\beta_2}{\dfrac{\lambda_1}{\lambda_2}} \le 1,\quad \dfrac{\beta_1}{\beta_2} \le \dfrac{\lambda_1}{\lambda_2}, \vspace{2mm} \\
0 < A <\dfrac{\lambda_1\lambda_2}{\beta_1 \lambda_2 - \beta_2 \lambda_1}, &\dfrac{\beta_1}{\beta_2} > \dfrac{\lambda_1}{\lambda_2},\quad \dfrac{\beta_1}{\beta_2} \ge \left( \dfrac{\lambda_1}{\lambda_2} \right)^2, \vspace{2mm}\\
0 < A <\dfrac{(\lambda_1 - \lambda_2)^2}{\left[(\beta_1\lambda_1)^{\frac13} - (\beta_2\lambda_2)^{\frac13}\right]^3}, &\dfrac{\beta_1}{\beta_2} < \left( \dfrac{\lambda_1}{\lambda_2} \right)^2 ,\quad \dfrac{\beta_1}{\beta_2}{\dfrac{\lambda_1}{\lambda_2}} > 1.
\end{cases}
\end{equation}
\end{enumerate}
\end{remark}}



%
\begin{figure}[h]
\centering
\begin{tikzpicture}[scale=2.0]

\fill[blue, opacity=0.5] plot[domain=1:3, variable=\x, smooth] ({\x},{1/\x}) --(3,0) --(0,0) --(1,1); 
\fill[orange, opacity=0.5] plot[domain=1:sqrt(3), variable=\x, smooth] (\x,{\x^2}) --(3,3) --(3,1); 
\fill[orange, opacity=0.5] plot[domain=1:3, variable=\x, smooth] (\x,{1/\x}) --(3,1); 
\fill[gray, opacity=0.6] plot[domain=1:sqrt(3), variable=\x, smooth] ({\x},{\x^2}) --(0,3) --(0,0) --(1,1); 

 \draw[->,>=stealth,semithick] (-0.25,0)--(3,0) node [below]{{\tiny $\frac{\lambda_1}{\lambda_2}$}}; 
 \draw[->,>=stealth,semithick] (0,-0.25)--(0,3) node [above, left]{{\tiny $\frac{\beta_1}{\beta_2}$}}; 
 \draw node[below right]{\textrm{O}}; 
  
 \draw[densely dotted, variable=\x, domain=0:0.25] plot(-\x,{\x^2}) ;
 \draw[densely dotted, variable=\x, domain=0:1] plot(\x,\x^2) ;

 \draw[semithick, domain=1:sqrt(3)+0.02, variable=\x] plot(\x,\x^2) node[above]{{\tiny $\frac{\beta_{\scalebox{0.7}{$1$}}}{\beta_{\scalebox{0.7}{$2$}}}=\big(\frac{\lambda_{\scalebox{0.7}{$1$}}}{\lambda_{\scalebox{0.7}{$2$}}}\big)^{\scalebox{0.6}{$2$}}$}};
 
 \draw[semithick, domain=1:3.05, variable=\x] plot(\x,{1/\x}) node[right]{{\tiny $\frac{\beta_{\scalebox{0.7}{$1$}}}{\beta_{\scalebox{0.7}{$2$}}}\frac{\lambda_{\scalebox{0.7}{$1$}}}{\lambda_{\scalebox{0.7}{$2$}}}=1$}};
 \draw[densely dotted, variable=\x, domain={1/3}:1] plot(\x,{1/\x})node[right]{};

 \draw[densely dotted, domain=1:3.1, variable=\x] plot(\x,1) node[right]{{\tiny $\frac{\beta_{\scalebox{0.7}{$1$}}}{\beta_{\scalebox{0.7}{$2$}}}=1$}};
 \draw[densely dotted, domain=-0.25:1, variable=\x] plot(\x,1);
 
 \draw[densely dotted, variable=\x, domain=0:3.125] plot(1,3-\x) node[below]{{\tiny $\frac{\lambda_{\scalebox{0.7}{$1$}}}{\lambda_{\scalebox{0.7}{$2$}}}=1$}};
 
 \draw[densely dotted, domain=1:3, variable=\x] plot(\x,\x) ;
 \draw[semithick, domain=0:1.125, variable=\x] plot(1-\x,1-\x) node[below left]{{\tiny $\frac{\beta_{\scalebox{0.7}{$1$}}}{\beta_{\scalebox{0.7}{$2$}}}=\frac{\lambda_{\scalebox{0.7}{$1$}}}{\lambda_{\scalebox{0.7}{$2$}}}$}};
 
 \draw (0.80,2.5) node[scale=0.8] {\hbox {\small $0 < A < \dfrac{\lambda_1\lambda_2}{\beta_1\lambda_2 - \beta_2\lambda_1}$}} ; 
 \draw (2.1,1.2) node[scale=0.7] {\hbox {\small $0 < A < \dfrac{(\lambda_1-\lambda_2)^2}{(\sqrt[3]{\beta_1 \lambda_1} - \sqrt[3]{\beta_2\lambda_2})^3}$}} ; 
   \draw (1.4,0.3) node[scale=1.0] {\hbox {\small $A>0$}} ; 

\end{tikzpicture}
\caption{Monotonicity of $g_A(\tau)$}
\end{figure}
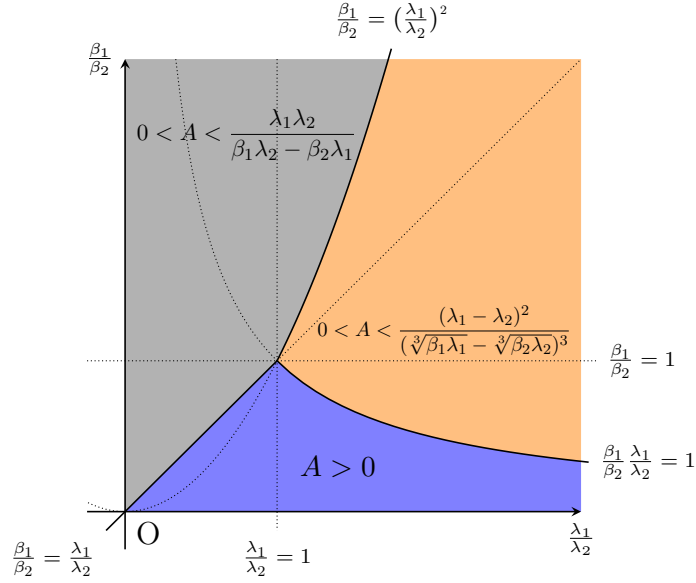
Assume that the constant $A>0$ satisfies 
\[
 A
 >
 A^\ast ,
\]
where the constant {$A^\ast$} defined by relation \eqref{eq;cond-A-th}.
Also, set 
\[
 M 
 \coloneqq  
 M(A)
 =
 \max_{\xi \in \mathbb{R}^n}
 h_A(\xi).
\]
Then we emphasize from Proposition \ref{prop;ess-sp-li} that the constant $M \in \mathbb{R}$ is strictly positive for any $A > A^\ast$,
that is,
\[
 M = M(A) > 0,\quad
 A > A^\ast.
\]
If $A>A^\ast$, 
we can not remove the exponential growth from the $L^p$-$L^q$ type estimate of the analytic semigroup $\mathrm{e}^{t\mathcal{L}_A}$.
\begin{lemma}\label{lem;lplqe}
Let $1 < q \le p < \infty$. 
Then{,} for all $\varepsilon > 0$ and all constants $A>{A^\ast}$, 
there exists a constant 
$
 C
 =
 C(\varepsilon, p,q,n,A) 
 > 
 0
$ 
such that 
\begin{align}
 &\label{eq;lplqe}
  \| 
   \mathrm{e}^{t\mathcal{L}_A} 
   f
  \|_p 
  \le 
  C
  t^{-\frac{n}2(\frac1q-\frac1p)}
  \mathrm{e}^{(M+\varepsilon)t}
  \|
   f
  \|_q, \\
 &\label{nlplqe}
  \|
   \nabla 
   \mathrm{e}^{t\mathcal{L}_A}
   f
  \|_p
  \le 
  C
  t^{-\frac{n}2(\frac1q-\frac1p)-\frac12}
  \mathrm{e}^{(M+\varepsilon)t}
  \|
   f
  \|_q
\end{align}
for all $t>0$ and all $f \in L^{{q}}(\mathbb{R}^n)$, 
where $M>0$ is a constant defined by Proposition \ref{prop;ess-sp-li}.
\end{lemma}

\begin{proof}
We first show that {for all $\varepsilon>0$ and all $1<p<\infty$,} there exists 
$C=C(\varepsilon,p,n) > 0$ depending only on $\varepsilon>0, 1 < p < \infty${, $A>0$} and the dimension $n$ such that 
\begin{equation}\label{eq;lplpe_1}
 \|
  \mathrm{e}^{t\mathcal{L}_A}f
 \|_p
 \le 
 C
 \mathrm{e}^{(M+\varepsilon)t}
 \|
  f
 \|_p, \quad t>0.
\end{equation}
This is {a} direct {consequence} of the well-known fact that 
the analytic semigroup $\mathrm{e}^{t\mathcal{L}_A}$ on $L^p(\mathbb{R}^n)$ has the relation of spectral bound 
(see \textit{e.g.} \cite[Ch.IV, Corollary 3.12]{EN00})
\[
\begin{split}
 & \, \sup
 \big\{
  \mathrm{Re} \, \mu 
  \mid
  \mu \in \sigma(\mathcal{L}_A)
 \big\}\\
 &= 
 \inf 
 \{
  \omega \in \mathbb{R}
  \mid
  \text{there exists } M_\omega \ge 1 \text{ such that } 
  \| \mathrm{e}^{t\mathcal{L}_A} \|_{L^p \to L^p}
  \le
  M_\omega
  \mathrm{e}^{\omega t} 
  \text{ for all } t {>} 0
 \}.
\end{split}
\]
Thus we obtain estimate \eqref{eq;lplpe_1}. 

For all $\varepsilon > 0$ and all $1 < p < \infty$, 
we claim that there exists {a constant} $C=C(\varepsilon,p,n) > 0 $ such that
\begin{equation}\label{eq;nlplpe}
 \|
  \nabla 
  \mathrm{e}^{t\mathcal{L}_A}
  f
 \|_p
 \le 
 C
 t^{-\frac12}
 \mathrm{e}^{(M+\varepsilon)t}
 \|
  f
 \|_p, \quad t>0.
\end{equation}
The proof will be divided into {five} steps.
First, since
$
 \partial_j
 \mathrm{e}^{t\mathcal{L}_A}f
 =
 \mathcal{F}^{-1}
 \big[
  (i\xi_j)\mathrm{e}^{th_A}\mathcal{F}[f]
 \big]
$ 
for
{all} $t>0$ and {each} $j=1,2,\ldots,n$, 
{the relation}
\[
 \begin{split}
   \mathrm{e}^{-(M+\varepsilon)t}
   \partial_j
   \mathrm{e}^{t\mathcal{L}_A}f
   = \, &
   \mathcal{F}^{-1}
   \big[
    (i\xi_j)
    \mathrm{e}^{-t(M+\varepsilon-h_A)}
    \mathcal{F}[f]
   \big]\\ 
  = \, & {(2\pi)^{-\frac{n}{2}}}
  \mathcal{F}^{-1}[i\xi_j \mathrm{e}^{-t(M+\varepsilon-h_A)}]*f
 \end{split}
\]
holds for all {$\varepsilon, t>0$} and {each $j=1,2,\ldots,n$}. 

The second step is to prove 
\begin{equation}\label{eq;es_771}
 \Big\| 
  \mathcal{F}
  \big[
  \mathcal{F}^{-1}
   [
    (i\xi_j)
    \mathrm{e}^{-t(M+\varepsilon-h_A)}
   ]
   \big]
 \Big\|_2
 \le
 Ct^{-\frac{n}4-\frac12}
\end{equation}
for some $C=C(\varepsilon, \beta_1, \beta_2, \lambda_1, \lambda_2, \varepsilon, n, A) > 0$ and {any $t>1$}.
To see this,
from Proposition \ref{prop;ess-sp-li},
we observe $h_A(\xi) \le M$ for all $\xi \in \mathbb{R}^n$.
Thus{,} for all $\varepsilon > 0$ and all $\xi \in \mathbb{R}^n$, 
\[
{{\mathcal H}}_A(\xi)
 \coloneqq 
 M
 +
 \varepsilon
 -
 h_A(\xi)
 \ge 
 \varepsilon
\]
and
\[
 \lim_{|\xi|\to\infty}
 \frac{|\xi|^2}{{{\mathcal H}}_A(\xi)}
 =
 \lim_{|\xi|\to \infty}
 \frac{1}{\displaystyle \frac{M+\varepsilon}{|\xi|^2}+1-\frac{A\beta_1}{\lambda_1+|\xi|^2}+\frac{A\beta_2}{\lambda_2+|\xi|^2}}
 =
 1
\]
hold.
Therefore we obtain
\begin{equation}\label{eq;bouned76}
 0
 \le 
 \max_{\xi\in \mathbb{R}^n}
  \frac{|\xi|^2}{{{\mathcal H}}_A(\xi)}
 <
 \infty.
\end{equation}
According to estimate \eqref{eq;bouned76}, 
we have
\[
\begin{split}
 \big|
  i
  \xi_j 
  \mathrm{e}^{-t{{\mathcal H}}_A(\xi)}
 \big|
 \le \, &
 t^{-\frac{n+1}2}
 \big[
  t{{\mathcal H}}_A(\xi)
 \big]^{\frac{n+1}2}
 \mathrm{e}^{-t{{\mathcal H}}_A(\xi)}
 \left(
  \frac{|\xi|^2}{{{\mathcal H}}_A(\xi)}
 \right)^{\frac{n+1}2}
 |\xi|^{-n}
 \le  
  C
  t^{-\frac{n}4-\frac12}
  |\xi|^{-n}
\end{split}
\]
for all {$t>1$} and all $|\xi| \ge 1$.
On the other hand, {we get}
\[
\begin{split}
 |i\xi_j \mathrm{e}^{-t{{\mathcal H}}_A(\xi)}| 
 \le \, & 
 t^{-\frac{n+1}2}
 |\xi|
 [t{{\mathcal H}}_A(\xi)]^{\frac{n+1}2}
 \mathrm{e}^{-t{{\mathcal H}}_A(\xi)}
 \left(
  \frac{1}{{{\mathcal H}}_A(\xi)}
 \right)^{\frac{n+1}2}
 \le 
 C
 t^{-\frac{n}4-\frac12}
\end{split}
\]
{for all $t >1$}  and {all} $|\xi|\le 1$.
Together with these estimates,
we {find that  for all $t > 1$,} 
\begin{equation}\label{eq;es_77}
 \begin{split}
 \big\|
  \mathcal{F}
  \big[
   \mathcal{F}^{-1}
   [
    (i\xi_j)
    \mathrm{e}^{-t{{\mathcal H}}_A(\xi)}
   ]
  \big]
 \big\|_2^2
 = \, &
  \left(
   \int_{|\xi| \le 1}+\int_{|\xi| \ge 1}
  \right)
  \big| 
   i\xi_j 
   \mathrm{e}^{-t{{\mathcal H}}_A(\xi)}
  \big|^2
  \, 
  \mathrm{d}\xi
  \\
  \le \, & 
  C
  t^{-\frac{n}2-1}
   \int_{|\xi|\le 1}
   \, \mathrm{d}\xi
   +
   \int_{|\xi|\ge 1}
    |\xi|^{-2n}
   \, \mathrm{d}\xi
  \\
  = \, &
  C
  t^{-\frac{n}2-1}
  \left(
   1
   +
   \int_{1}^{\infty}
   r^{-n-1}
   \, \mathrm{d}{r}
  \right)\\
  = \, & 
  C
  t^{-\frac{n}2-1},
 \end{split}
\end{equation}
which yields estimate \eqref{eq;es_771}.

We thirdly show that 
\[
 \big\|
  \partial_\xi^\alpha
  \mathcal{F}
  \big[
   \mathcal{F}^{-1}
   [
    i\xi_j
    \mathrm{e}^{-t{{\mathcal H}}_A}
   ]
  \big]
 \big\|_2
 \le 
 C
 t^{-\frac{n}4-\frac12+\frac{|\alpha|}{2}}, 
 \quad {t > 1}.
\]
For {any} $\alpha \in \mathbb{Z}_+^n$ with 
$
 |\alpha| 
 > 
 {n}/2
$,
Leibniz's rule and the estimate
$
 |\partial_\xi^\gamma\xi_j|\le (1+|\xi|)^{1-|\gamma|} \, (\xi \in \mathbb{R}^n)
$
give 
\begin{equation}\label{eq;es_78}
\begin{split} 
 |
  \partial_\xi^\alpha
  (
   i\xi_j 
   \mathrm{e}^{-t{{\mathcal H}}_A(\xi)}
  )
 | 
 = \, &
 \left|
  i
  \sum_{\gamma\le \alpha}
  \binom{\alpha}{\gamma} 
  \big(
  \partial_\xi^{\alpha-\gamma} 
  \xi_j
  \big)
  \big(
  \partial_\xi^{\gamma}
  \mathrm{e}^{-t{{\mathcal H}}_A(\xi)}
  \big)
 \right| \\
 \le \, & 
 \sum_{\gamma\le \alpha}
 \binom{\alpha}{\gamma}
 \big|
  \partial_\xi^{\alpha-\gamma}\xi_j
 \big|
 \big|
  \partial_\xi^{\gamma}
  \mathrm{e}^{-t{{\mathcal H}}_A(\xi)}
 \big|\\
 \le \, &
 \sum_{\gamma\le \alpha}
 \binom{\alpha}{\gamma}
 (1+|\xi|)^{1-|\alpha|+|\gamma|}
 |
  \partial_\xi^{\gamma}
  \mathrm{e}^{-t{{\mathcal H}}_A(\xi)}
 |.
\end{split}
\end{equation}
Applying formula \eqref{eq;lem2-1} to $f(r)=\mathrm{e}^r$ and $g(\xi)=-t{\mathcal{H}}_A(\xi)$,
we observe from Lemma \ref{lem;der-sy} and $\mathcal{H}_A(\xi)=M+\varepsilon-h_A(\xi)$ that 
the relation
\begin{equation}\label{eq;es_79}
 \begin{split}
 |\partial_\xi^{\gamma}\mathrm{e}^{-t{{\mathcal H}}_A(\xi)}| 
 \le \, & 
 \mathrm{e}^{-t{{\mathcal H}}_A(\xi)}
 \sum_{k=1}^{|\gamma|}
 t^k
 \sum_{\substack{\gamma_1+\cdots+\gamma_k=\gamma, \\ |\gamma_i|\ge 1}}
 \left|
  \Gamma_{\gamma_1, \ldots, \gamma_k}^k
 \right|
 |
  (\partial_\xi^{\gamma_1}{{\mathcal H}}_A(\xi))
  \cdots
  (\partial_\xi^{\gamma_k}{{\mathcal H}}_A(\xi))
 | \\
 = \, & 
 \mathrm{e}^{-t{{\mathcal H}}_A(\xi)}
 \sum_{k=1}^{|\gamma|}
 t^k
 \sum_{\substack{\gamma_1+\cdots+\gamma_k=\gamma, \\ |\gamma_i|\ge 1}}
 \left|
  \Gamma_{\gamma_1, \ldots, \gamma_k}^k
 \right|
 |
  ( 
   \partial_\xi^{\gamma_1}
   h_A(\xi)
   )
  \cdots
  (
   \partial_\xi^{\gamma_k}
   h_A(\xi)
  )
 | \\
 \le \, & 
 \mathrm{e}^{-t{{\mathcal H}}_A(\xi)}
 \sum_{k=1}^{|\gamma|}
  C_k{t^k}
  (1+|\xi|)^{2k-|\gamma|}
 \end{split}
\end{equation}
holds for {all $\gamma\in {\mathbb Z}_+^n$ with} $| \gamma | \ge  1$.
{Remarking estimate \eqref{eq;bouned76} and} combining inequality \eqref{eq;es_78} with inequality \eqref{eq;es_79}, we obtain
\begin{equation}\label{eq;es_710}
 \begin{split}
{t^{\frac{n}{4}+\frac{1}{2}-\frac{|\alpha|}{2}}}&\big|
  \partial_\xi^\alpha
  \big(
   i\xi_j 
   \mathrm{e}^{-t{{\mathcal H}}_A(\xi)}
  \big)
 \big| \\
 \le \, &{Ct^{\frac{n}{4}-\frac{|\alpha|}{2}}|\xi|^{-|\alpha|}\left( 
   \frac{|\xi|^{2}}{{\mathcal H}_A(\xi)}
  \right)^{\frac12}\big[ 
   t{\mathcal H}_A(\xi)
  \big]^{\frac12}\mathrm{e}^{-t{\mathcal H}_A(\xi)}} \\
 &+{t^{\frac{n}{4}-\frac{|\alpha|}{2}}}
 \sum_{\gamma\le \alpha{, \ |\gamma|\ge 1}}
 \sum_{k=1}^{|\gamma|}
  C_k
  \binom{\alpha}{\gamma}
  |\xi|^{-|\alpha|}
  \left( 
   \frac{|\xi|^{2}}{{{\mathcal H}}_A(\xi)}
  \right)^{k+\frac12}
  \big[ 
   t{{\mathcal H}}_A(\xi)
  \big]^{k+\frac12}
  \mathrm{e}^{-t{{\mathcal H}}_A(\xi)}\\
  \le &{C|\xi|^{-|\alpha|}}
\end{split}
\end{equation}
for any {$t>1$} and any $|\xi|\ge 1$.
On the other hand, 
for all {$t >1$} and all $|\xi| \le 1$,
we have
\begin{equation}\label{eq;es_711}
\begin{split}
{t^{\frac{n}{4}+\frac{1}{2}-\frac{|\alpha|}{2}}} & \, |
  \partial_\xi^\alpha
  (
   i\xi_j 
   \mathrm{e}^{-t{{\mathcal H}}_A(\xi)}
  )
 |\\
 \le &\, {t^{\frac{n}{4}-\frac{|\alpha|}{2}}|\xi|^{-|\alpha|}\left( 
   \frac{1}{{\mathcal H}_A(\xi)}
  \right)^{\frac12}\big[ 
   t{\mathcal H}_A(\xi)
  \big]^{\frac12}\mathrm{e}^{-t{\mathcal H}_A(\xi)}} \\
 &+{t^{\frac{n}{4}-\frac{|\alpha|}{2}}} 
 \sum_{\gamma \le \alpha{, \ |\gamma|\ge 1}}
 \sum_{k=1}^{|\gamma|}
 C_k
 \binom{\alpha}{\gamma}{|\xi|^{-|\alpha|}\left(\frac{1}{{\mathcal H}_A(\xi)}\right)^{k+\frac{1}{2}}[t{\mathcal H}_A(\xi)]^{k+\frac{1}{2}}}
 \mathrm{e}^{-t{{\mathcal H}}_A(\xi)}\\
 \le \, & 
{C|\xi|^{-|\alpha|}}.
\end{split}
\end{equation}
Therefore, 
by the similar argument as that of proof for estimate \eqref{eq;es_77},
inequality \eqref{eq;es_710} and inequality \eqref{eq;es_711} lead to 
\begin{equation}\label{eq;es_712}
 \|
  \partial_\xi^\alpha
  \mathcal{F}
  [ 
   \mathcal{F}^{-1}
   [
    i\xi_j
    \mathrm{e}^{-t{{\mathcal H}}_A}
   ]
  ]
 \|_2
 \le 
 C
 t^{-\frac{n}4-\frac12+\frac{|\alpha|}2}
\end{equation}
for any $t > 1$ and any $\alpha\in \mathbb{Z}_+^n$ with 
$ 
 |\alpha| 
 > 
 {n}/2
$.

In {the fourth} step, we {are going to} prove the desired estimate \eqref{eq;nlplpe}. 
Let $N$ be a positive integer with $N > {n}/{2}$. 
Then, 
for all 
$
 \alpha \in {\mathbb Z}_+^n
$ 
with 
$
 |\alpha|
 =
 N
$ 
and all $t>1$, 
applying estimate \eqref{eq;es_v1}, estimate \eqref{eq;es_771} and estimate \eqref{eq;es_712} to 
$
 \mathcal{F}^{-1}
   [
    i\xi_j
    \mathrm{e}^{-t{{\mathcal H}}_A}
   ]
$
gives
\begin{equation}\label{eq;es_713}
\begin{split}
 \|
  &\mathcal{F}^{-1}
   [
    i\xi_j
    \mathrm{e}^{-t{{\mathcal H}}_A}
   ]
 \|_1 \\
 \le \, &
 C
 \| 
  \mathcal{F}
  [
   \mathcal{F}^{-1}
   [
    i\xi_j
    \mathrm{e}^{-t{{\mathcal H}}_A}
   ]
  ]
 \|_2^{1-\frac{n}{2{N}}}
 \left(\sum_{|\alpha|=N}\|
  \partial_\xi^\alpha
  \mathcal{F}
  [
   \mathcal{F}^{-1}
   [
    i\xi_j
    \mathrm{e}^{-t{{\mathcal H}}_A}
   ]
  ]
 \|_2\right)^{\frac{n}{2{N}}}\\
 \le \, &
 C
 (t^{-\frac{n}4-\frac12})^{1-\frac{n}{2{N}}}
 (t^{-\frac{n}4-\frac12+\frac{{N}}2})^{\frac{n}{2{N}}} \\
 = \, &
 C
 t^{-\frac12}.
\end{split}
\end{equation}
Using inequality \eqref{eq;es_713}
together with Hausdorff-Young's inequality,
we obtain
\[
 \mathrm{e}^{-(M+\varepsilon)t}
 \|
  \partial_j
  \mathrm{e}^{t\mathcal{L}_A}
  f
 \|_p
 \le {(2\pi)^{-\frac{n}{2}}}
 \|
  \mathcal{F}^{-1}
   [
    i\xi_j
    \mathrm{e}^{-t{{\mathcal H}}_A}
   ]
 \|_1
 \|
  f
 \|_p
 \le 
 C
 t^{-\frac12}
 \|
  f
 \|_p
\]
for all $1 < p < \infty$ and all {$t>1$}, which implies 
\begin{equation}\label{eq;es_714}
 \mathrm{e}^{-(M+\varepsilon)t}
 \|
  \nabla \mathrm{e}^{t\mathcal{L}_A}
  f
 \|_p
 \le 
 C
 t^{-\frac12}
 \|
  f
 \|_p, \quad {t>1}.
\end{equation}
When $0 < t \le 1$, 
estimate \eqref{eq;L^p-L^q-derivative-any-A}
gives
\begin{equation}\label{eq;es_715}
 \begin{split}
 \|
  \nabla \mathrm{e}^{t\mathcal{L}_A}f\|_p
  \le \, &
  C
   t^{-\frac12}
   \mathrm{e}^{{2A\max\{\beta_1, \beta_2\}}t}
   \|
    f
   \|_p\\
  \le \, &
  C
   \mathrm{e}^{-(M+\varepsilon)t}t^{-\frac12}
   \mathrm{e}^{(M+\varepsilon)t}
   \| 
    f
   \|_p\\
  \le \, &
  C
  t^{-\frac12}
  \mathrm{e}^{(M+\varepsilon)t}
  \|
   f
  \|_p.
 \end{split}
\end{equation}
From inequality \eqref{eq;es_714} and inequality \eqref{eq;es_715},
we conclude {the} desired growth estimate \eqref{eq;nlplpe}.

Fifthly we will prove estimate \eqref{eq;lplqe}. {First of all,}  
for $1 < q < p < \infty$ with 
$
 1/q 
 - 
 1/p 
 < 
 1/n
$,
using the $L^p$-$L^p$ type estimates \eqref{eq;lplpe_1} and \eqref{eq;nlplpe}
together with
the Gagliardo--Nirenberg inequality
leads to
\begin{equation}\label{eq;es_715_2}
\begin{split}
 \|
  \mathrm{e}^{t\mathcal{L}_A}f
 \|_p
 \le \, &
 \|
  \nabla \mathrm{e}^{t\mathcal{L}_A}f
 \|_q^{n(\frac1q-\frac1p)}
 \|
  \mathrm{e}^{t\mathcal{L}_A}f
 \|_q^{1-n(\frac1q-\frac1p)}\\
\le \, &
\{ 
 C
 t^{-\frac12}
 \mathrm{e}^{(M+\varepsilon)t}\|f\|_q
\}^{n(\frac1q-\frac1p)}
\{
 C
 \mathrm{e}^{(M+\varepsilon)t}
 \|
  f
 \|_q
 \}^{1-n(\frac1q-\frac1p)}\\
\le \, &
 C
 t^{-\frac{n}2(\frac1q-\frac1p)}
 \mathrm{e}^{(M+\varepsilon)t}
 \|
  f
 \|_q, 
 \quad t>0.
\end{split}
\end{equation}

Next, for 
$1<q<p<\infty$ with 
$
 1/n
 \le 
 1/q
 -
 1/p
 <
 2/n
$,
taking the index 
$
 q
 <
 q_1
 <
 p
$
with 
$
 0
 <
 1/q
 -
 1/{q_1}
 <
 1/n
$
and
$
 1/{q_1}
 -
 1/p
 <
 1/n,
$
{the inequality}
\begin{equation}\label{eq;es_716}
\begin{split}
 \|
  \mathrm{e}^{t\mathcal{L}_A}
  f
 \|_p
 = \, & 
 \|
  \mathrm{e}^{\frac{t}2\mathcal{L}_A}
  (
   \mathrm{e}^{\frac{t}2\mathcal{L}_A}
   f
  )
 \|_p \\
 \le \, & 
 C
 \Big(
  \frac{t}2
 \Big)^{-\frac{n}2(\frac1{q_1}-\frac1p)}
 \mathrm{e}^{(M+\varepsilon)\frac{t}2}
 \|
  \mathrm{e}^{\frac{t}2\mathcal{L}_A}
  f
 \|_{q_1} \\
 \le \, & 
 C
 \Big(
  \frac{t}2
 \Big)^{-\frac{n}2(\frac1{q_1}-\frac1p)}
 \mathrm{e}^{(M+\varepsilon)\frac{t}2}
 \Big(
  \frac{t}2
 \Big)^{-\frac{n}2(\frac1q-\frac1{q_1})}
 \mathrm{e}^{(M+\varepsilon)\frac{t}2}
 \|
  f
 \|_{q}\\
 \le \, & 
 C
 t^{-\frac{n}2(\frac1q-\frac1p)}
 \mathrm{e}^{(M+\varepsilon)t}
 \|
  f
 \|_q, 
 \quad t>0.
 \end{split}
\end{equation}
{follows from \eqref{eq;es_715_2}.}
 
Now we prove estimate \eqref{eq;lplqe} in the case $1 < q < p < \infty$ satisfying 
$
 (k-1)/n \le {1}/{q}-{1}/{p} < {k}/{n}
$ 
for some $k\in {\mathbb N}$ with $3\le k\le n$. By the argument similar to that in \eqref{eq;es_716}, there are indices $q<q_{k-1}<q_{k-2}<\cdots<q_2<q_1<p$ such that 
\begin{align*}
 &\frac1{q_{1}}
 -
 \frac1{p}
 <
 \frac{1}n, \\
 &\frac{1}{q_{j}}
 -
 \frac{1}{q_{j-1}}
 <
 \frac{1}n, \quad j=2,3,\ldots,k-1, \\
 &0
 <
 \frac{1}{q}
 -
 \frac1{q_{k-1}}
 <
 \frac{1}{n}, 
\end{align*}
and it holds from inequality \eqref{eq;es_715_2} that 
\[
\begin{split}
\|\mathrm{e}^{t\mathcal{L}_A}f\|_p
= \, & 
\|
 \mathrm{e}^{\frac{t}k\mathcal{L}_A}
 (
  \mathrm{e}^{\frac{k-1}{k}t\mathcal{L}_A}f
 )
\|_p\\
 \le \, & 
 C
 \Big(
  \frac{t}k
 \Big)^{-\frac{n}2(\frac1{q_1}-\frac1p)}
 \mathrm{e}^{(M+\varepsilon)\frac{t}{k}}
 \|
  \mathrm{e}^{\frac{k-1}{k}t\mathcal{L}_A}
  f
 \|_{q_1} \\
 \le \, & 
 C
 \bigg(\frac{t}{k}\bigg)^{-\frac{n}2(\frac1{q_1}-\frac1p)-\sum_{j=2}^{k-1}\frac{n}2(\frac1{q_j}-\frac1{q_{j-1}})-\frac{n}2(\frac1{q}-\frac1{q_{k-1}})}
 \mathrm{e}^{(M+\varepsilon)t}
 \|
  f
 \|_{q}\\
 \le \, &
 C
 t^{-\frac{n}{2}(\frac{1}{q}-\frac{1}{p})}
 \mathrm{e}^{(M+\varepsilon)t}
 \|
  f
 \|_q, 
 \quad t>0.
 \end{split}
\]
Thus, gathering the estimates above, we {obtain the desired} estimate \eqref{eq;lplqe} since 
\[
 (p,q) 
 \in 
 \left\{
  (r,s) 
  \, 
  \mid 
  \, 
  \frac{1}{s}-\frac{1}{r}<\frac{1}{n}
 \right\} 
 \cup 
 \bigcup_{k=2}^n 
  \left\{
   (r,s) 
   \, 
   \mid 
   \,
   \frac{k-1}{n}
   \le 
   \frac{1}{s}-\frac{1}{r}
   <
   \frac{k}{n}
  \right\}
\]
if $1<q<p<\infty$.

Finally, we are going to show estimate \eqref{nlplqe}. 
By making use of estimate \eqref{eq;lplqe} and estimate \eqref{eq;nlplpe}, 
there exists a constant $C>0$ independent of $t>0$ such that
\begin{align*}
 \|
  \nabla 
  \mathrm{e}^{t\mathcal{L}_A}
  f
 \|_p
  =& 
  \|
   \nabla 
   \mathrm{e}^{\frac{t}{2}\mathcal{L}_A}
   (
    \mathrm{e}^{\frac{t}{2}\mathcal{L}_A}
    f
   )
  \|_p\\
  \le \, &
  C
  t^{-\frac{1}{2}}
  \mathrm{e}^{(M+\varepsilon)\frac{t}{2}}
  \|
   \mathrm{e}^{\frac{t}{2}\mathcal{L}_A}
   f
  \|_{p}
  \le 
  C
  t^{-\frac{n}2(\frac1q-\frac1p)-\frac12}
  \mathrm{e}^{(M+\varepsilon)t}
  \|
   f
  \|_q
\end{align*}
for $1 < q \le p < \infty$ and $t>0$. 
This completes the proof of Lemma \ref{lem;lplqe}.
\end{proof}

In the symbol $h_{A}(\xi)$ of the operator ${\mathcal L}_A$ defined by relation \eqref{eq;def-li-op-cl},
the constant $M=M(A){>0}${, which is the maximum value of the symbol $h_A(\xi)$}, 
plays {an important} roll as follows:

\begin{lemma}\label{lem;spec-appr}
Let $A>0$ and $1 < p < \infty$.
Then{,} for all $0 < \gamma \le 1$ and all $T>0$, 
there exists $u_{0} \in L^p(\mathbb{R}^n)$ such that 
\begin{align}
\label{75}
 \|
  \mathrm{e}^{t\mathcal{L}_A}
  u_0
  -
  \mathrm{e}^{Mt}
  u_0
 \|_p
 \le  \, &
 \gamma
 \|
  u_0
 \|_p, \\
\label{76}
\|
 \mathrm{e}^{t\mathcal{L}_A}
 u_0
\|_p
\le \, & 
 2
 \mathrm{e}^{Mt}
 \|
  u_0
 \|_p.
\end{align}
for {each} $0 \le t \le T$.
Here $M$ is a constant given in Proposition \ref{prop;ess-sp-li}.
\end{lemma}

\begin{proof}
By Proposition \ref{prop;ess-sp-li},
$\sigma(\mathcal{L}_A)=(-\infty, M]$ holds for $A>0$.
Applying \cite[Lemma 1]{ShSt} to the linearized operator $\mathcal{L}_A$ gives
\[
 \inf_{\substack{u \in D(\mathcal{L}_A) \\ \|u\|_p=1}}\,
 \|
  \mathcal{L}_A
  u
  -
  M
  u
 \|_p
 =
 \inf_{u \in D(\mathcal{L}_A) \setminus \{ 0 \}}
 \frac{\| \mathcal{L}_A u - M u \|_p}{\| u \|_p}=0.
\]
This property shows that there exists $u_\varepsilon \in L^p(\mathbb{R}^n) \setminus \{ 0 \}$ depending on $\varepsilon > 0$ such that
\begin{equation}\label{es77}
 \|
  \mathcal{L}_A
  u_\varepsilon
  -
  M
  u_\varepsilon 
 \|_p
 \le 
 \varepsilon
 \|
  u_\varepsilon
 \|_p.
\end{equation}
Since the semigroup $\{ \mathrm{e}^{t\mathcal{L}_A} \}_{t>0}$ acts an analytic semigroup on $L^p(\mathbb{R}^n)$
together with the $L^p$-$L^q$ estimate \eqref{eq;L^p-L^q-any-A} and inequality \eqref{es77}, 
we obtain
\begin{equation}\label{es78}
 \begin{split}
 \big\|
  \mathrm{e}^{t\mathcal{L}_A}
  u_\varepsilon
  -
  \mathrm{e}^{Mt}
  u_\varepsilon
 \big\|_p
 = \, &
 \left\|
  \int_0^1
   \frac{\mathrm{d}}{\mathrm{d}s}
   \mathrm{e}^{ts\mathcal{L}_A}
   (\mathrm{e}^{Mt(1-s)}u_{\varepsilon})
   \, \mathrm{d}s
 \right\|_p\\
 \le \, &
  t
  \int_0^1
   \big\|{\mathrm{e}^{Mt(1-s)}
    \mathcal{L}_A(
      \mathrm{e}^{ts\mathcal{L}_A}
       u_{\varepsilon})
     -
     M\mathrm{e}^{Mt(1-s)}
     (
      \mathrm{e}^{ts\mathcal{L}_A}
       u_{\varepsilon}
     )}
   \big\|_p
   \, \mathrm{d}s\\
  = \, &
   t
   \int_0^1
    \mathrm{e}^{Mt(1-s)}
    \big\|
     \mathrm{e}^{ts\mathcal{L}_A}
     (
      \mathcal{L}_A 
      u_{\varepsilon}
      -
      M
      u_\varepsilon
     )
    \big\|_p
   \, \mathrm{d}s\\
 \le \, &
  C
   t
    \int_0^1 
    \mathrm{e}^{Mt(1-s)}
    \cdot 
    \mathrm{e}^{{2A\max\{\beta_1, \beta_2\}}ts}
    \big\|
     \mathcal{L}_A 
     u_{\varepsilon}
     -
     M
     u_\varepsilon
    \big\|_p
    \, \mathrm{d}s \\
 \le \, & 
 C
  \varepsilon 
  t
  \mathrm{e}^{Mt}
  \|
   u_\varepsilon
  \|_p
  \int_0^1 
   \mathrm{e}^{{2A\max\{\beta_1, \beta_2\}}ts}
  \, \mathrm{d}s\\
 \le \, &
 C
 \varepsilon 
 T
 \mathrm{e}^{({2A\max\{\beta_1, \beta_2\}}+M)T}
 \|
  u_\varepsilon
 \|_p.
 \end{split}
\end{equation}
for each $0 \le t \le T$.
Set 
\[
 \varepsilon_0
 \coloneqq 
 \frac{\gamma}{CT\mathrm{e}^{({2A\max\{\beta_1, \beta_2\}}+M)T}}
 >
 0
\]
and
$u_0 \coloneqq  u_{\varepsilon_0}$
for given $T>0$ and given $0 < \gamma \le 1$.
Then inequality \eqref{es78} leads to 
\[
 \big\|
  \mathrm{e}^{t\mathcal{L}_A}
  u_0
  -
  \mathrm{e}^{Mt}
  u_0
 \big\|_p
 \le 
 \gamma
 \|
  u_0
 \|_p,
\]
which is {the} desired estimate \eqref{75}. 
Moreover, applying estimate \eqref{75} to the case $\gamma=1$ shows
\[
\begin{split}
 \|
  \mathrm{e}^{t\mathcal{L}_A}
  u_0
 \|_p
 = \, &
 \big\|
  \mathrm{e}^{Mt}u_0
  +
  \mathrm{e}^{t\mathcal{L}_A}
  u_0
  -
  \mathrm{e}^{Mt}
  u_0
 \big\|_p\\
 \le \, &
 \mathrm{e}^{Mt}
 \|
  u_0
 \|_p
 +
 \big\|
  \mathrm{e}^{t\mathcal{L}_A}u_0
  -
  \mathrm{e}^{Mt}u_0
 \big\|_p\\
 \le \, &
 2
 \mathrm{e}^{Mt}
 \|
  u_0
 \|_p{,}
\end{split}
\]
which yields estimate \eqref{76}. 
\end{proof}

\section{Perturbed problem near constant steady states}
Let $u_C \equiv A$ be the constant steady states of problem \eqref{eq;simplified-ar-ks}.
Then, as stated in Section 4, by introducing the new function $v$ as $v \coloneqq  {u} - A$, 
problem \eqref{eq;simplified-ar-ks} is converted to the following perturbation problem \eqref{eq;lin}:
\begin{equation}\label{s;Q1}
\tag{Q}
\begin{cases}
 \partial_t {v} - \Delta {v} + A \Delta K*{v} = -\nabla \cdot ( {v} \nabla K*{v} ), &t>0, \ x \in \mathbb{R}^n, \\
 {v}(0,x) = {v}_0(x), &x\in \mathbb{R}^n,
\end{cases}
\end{equation}
where {$K(x):=(\beta_1B_{\lambda_1}-\beta_2B_{\lambda_2})(x)$ and $v_0(x):=u_0(x) - A$}. 
First of all, we begin with local in time existence of mild solutions to problem \eqref{s;Q1}, which satisfies the following integral equation:
\[
 v(t)
 = 
 \mathrm{e}^{t\mathcal{L}_A}{v}_0+N[{v},{v}](t), \quad {0<t<T,}
\]
where
\begin{equation} 
\label{pp-nuv}
 N[f,g]{(t)}
 = 
 -
 \int_0^t
 \nabla \cdot \mathrm{e}^{(t-s)\mathcal{L}_A} (f \nabla K*g)(s)
 \, \mathrm{d}s.
\end{equation}

\subsection{Existence of local-in-time solutions to problem \eqref{eq;lin}}
We now state existence of local-in-time solutions to problem \eqref{eq;lin}.

\begin{proposition}\label{prop;812}
Assume that $A>0$ and a pair of index $(p,q_*)$ satisfies
\begin{equation}\label{dfn;pq}
 \max\Big\{ 1, \frac{n}2 \Big\}
  \le 
  p
  <
 \infty, 
 \quad 
 q_*
 =
 \min
 \left\{
  \frac{2np}{2n-p}, 2p
 \right\}.
\end{equation}
Then{,} for any ${v}_0 \in L^p(\mathbb{R}^n)$, 
there { exist a number} $T>0$ and a local-in-time mild solution {v} to {problem} \eqref{eq;lin}
such that
$
 {v}
 \in 
 C([0,T]; L^p(\mathbb{R}^n)) \cap C((0,T]; L^{q_*}(\mathbb{R}^n))
$. 
\end{proposition}

\begin{proof}
For fixed $0 < T < \infty$, 
set a Banach space $\mathscr{X}_T$ by
\[
 \mathscr{X}_T 
 \coloneqq 
 \left\{
 {v}  \in X_T 
  \mid
  \sup_{0 \le t \le T}
  \| {v}(t) \|_p \le B,\ 
  \sup_{0< t \le T}
  t^{\frac{n}{2} ( \frac{1}{p} - \frac{1}{q_*} )}
  \| {v}(t) \|_{q_*} \le \eta
 \right\}
\]
with the norm
\[
 \|
  {v}
 \|_{\mathscr{X}_T} 
 \coloneqq 
 \sup_{0 \le t \le T}
 \| {v}(t) \|_p
 +
 \sup_{0 < t \le T}
 t^{\frac{n}{2} ( \frac{1}{p} - \frac{1}{q_*} )}
 \| {v}(t) \|_{q_*},
\]
where
\[
 X_T 
 \coloneqq 
 C([0,T];L^p(\mathbb{R}^n)) \cap C((0,T];L^{q_*}(\mathbb{R}^n))
\]
and the two constants $B > 0$ and $\eta > 0$ are  chosen later.
Also, define
\begin{equation}\label{Ieu}
 \Phi[{v}](t) 
 \coloneqq 
 \mathrm{e}^{t\mathcal{L}_A}
 {v}_0
 +
 N[{v},{v}](t), 
 \quad 0<t<T.
\end{equation}
Then we claim that 
$
 \Phi
$
is {a} contraction mapping from $\mathscr{X}_T$ to $\mathscr{X}_T$ for $0 < T \ll 1$.
To see this,
we have divided the proof into some steps. 

Firstly,
applying 
Lemma \ref{eq;L^p-L^q-any-A}
to $\mathrm{e}^{t\mathcal{L}_A}$ leads to 
\begin{align}
 &
  \label{lcm1}
   \|
    \mathrm{e}^{t\mathcal{L}_A}{v}_0
   \|_p
   \le 
   C_{L1}
   \mathrm{e}^{2A\max\{\beta_1, \beta_2\} t}
   \|
    {v}_0
   \|_p
   \le 
   C_{L}
   \mathrm{e}^{2A\max\{\beta_1, \beta_2\} T}
   \|
    {v}_0
   \|_p, \\
 &
  \label{lcm2}
   t^{\frac{n}{2} ( \frac{1}{p} - \frac{1}{q_*} )}
   \|
    \mathrm{e}^{t\mathcal{L}_A}{v}
    _0
   \|_{q_*}
   \le 
   C_{L2}
  \mathrm{e}^{2A\max\{\beta_1, \beta_2\}t}
   \|
    {v}_0
   \|_p
   \le 
   C_{L}
   \mathrm{e}^{2A\max\{\beta_1, \beta_2\}T}
   \|
    {v}_0
   \|_p{,}
\end{align}
provided $0<t<T$,
where 
\begin{align}
\label{c84}&
C_L
\coloneqq 
\max\{ C_{L1}, C_{L2} \}.
\end{align}

We next prove that there exists $C_{N1} > 0$ independent of $t$, $T$ such that 
\begin{equation}\label{cm1}
 \|
  N[f,g](t)
 \|_{p}
 \le 
  C_{N1}
  \mathrm{e}^{2A\max\{\beta_1, \beta_2\} T}
  D_{n,p}(T)
  [f]_t
  [g]_t
\end{equation}
for $f, g \in \mathscr{X}_T$ and $0<t<T$,
where 
\[
 [f]_t
 \coloneqq 
 \sup_{0<s\le t}
 s^{\frac{n}{2}(\frac{1}{p} - \frac{1}{q_*})}
 \| 
  f(s)
 \|_{q_*}
 \]
 and
 \[
 D_{n,p}(T)
 \coloneqq 
 \begin{cases}
  1, &n \ge 2, \, \dfrac{n}2 \le p \le n \quad \text{or} \quad n=p=1, \\
  T^{\frac12 (1-\frac{n}p)}, &n\ge 2, \ n < p < \infty \quad \text{or} \quad n=1, \, 1<p<\infty.
 \end{cases}
\]
Indeed,
for $q_*', q_1$ with
\[
 \frac{1}{p}
 =
 \frac{1}{q_*}+\frac{1}{q_*'}, \quad q_*'\ge 1, \quad 
 \frac{1}{q_*'}
 =
 \frac{1}{q_1}+\frac{1}{q_*}-1,
\]
the index $q_1$ satisfies
\[
 \frac{1}{q_1}
 =
 \begin{cases}
  \displaystyle 
  1+\frac{1}{n}-\frac{1}{p},\quad &n \ge 2, \, \dfrac{n}2 \le p \le n, \\
  \displaystyle 
  1, &n \ge 2, \ n < p < \infty \quad \text{or} \quad n=1, \ 1 \le p < \infty.
 \end{cases}
\]
Thus {the inequality}
$
 ({n-1})/{n}
 \le 
 {1}/{q_1}
 \le 
 1
$
{follows}.
By $\|\nabla K\|_{L_w^{q_1}}<\infty$ ({see Lemma \ref{lem;be} above}),
H\"{o}lder's inequality and Hausdorff-Young's inequality, 
we have
\[
\begin{split}
 \|
  f(t)(\nabla K*g)(t)
 \|_{p} 
 \le \, &
 \|
  f(t)
 \|_{q_*}
 \|
  \nabla K*g(t)
 \|_{q_*'} \\
 \le \, &
 C
 \| 
  \nabla K
 \|_{L_w^{q_1}}
 \|
  f(t)
 \|_{q_*}
 \|
  g(t)
 \|_{q_*}\\
 \le \, &
 C
 \|
  f(t)
 \|_{q_*}
 \|
  g(t)
 \|_{q_*}{,} 
\end{split}
\]
provided $0<t<T$.
According to estimate \eqref{eq;L^p-L^q-derivative-any-A} and 
\[
\begin{split}
 \|
  f(s)
 \|_{q_*}
 \le \, & 
 s^{-\frac{n}{2}(\frac{1}{p} - \frac{1}{q_*})}
 [f]_t,\\
 \|
  g(s)
 \|_{q_*}
 \le \, &
 s^{-\frac{n}{2}(\frac{1}{p} - \frac{1}{q_*})}
 [g]_t
\end{split}
\]
for  $0 < s < t$, we get
\[
\begin{split}
 \|
  N[f,g](t)
 \|_{p} 
 \le \, &
 C
 \int_0^t
 (t-s)^{-\frac12}
 \mathrm{e}^{2A\max\{\beta_1, \beta_2\} (t-s)}
 \|
  f(s)(\nabla K*g)(s)
 \|_{p}
 \, \mathrm{d}s \\
 \le \, & 
 C
 \mathrm{e}^{2A\max\{\beta_1, \beta_2\} T}
 [f]_t
 [g]_t
 \int_0^t
 (t-s)^{-\frac12}s^{-{n}(\frac{1}{p} - \frac{1}{q_*})}
 \, \mathrm{d}s\\
 \le \, & 
 C_{N1}
 \mathrm{e}^{2A\max\{\beta_1, \beta_2\} T}
 D_{n,p}(T)
 [f]_t
 [g]_t
\end{split}
\]
for $0<t<T$, which shows estimate \eqref{cm1}. {Also, by the same argument as in the proof of estimate \eqref{cm1}, we see} that there exists $C_{N2} > 0$ independent of $t$, $T$ such that
\begin{equation}\label{cm2}
 t^{\frac{n}{2} ( \frac{1}{p} - \frac{1}{q_*} )}
 \|
  N
  [
   f,g
  ](t)
 \|_{q_*}
 \le 
 C_{N2}
 \mathrm{e}^{ 2A\max\{\beta_1, \beta_2\} T }
 D_{n,p}(T)
 [f]_t
 [g]_t
\end{equation}
for $f, g \in \mathscr{X}_T$ and $0<t<T$.

Set 
$
 C_{N} 
 \coloneqq 
 \max\{
  C_{N1}, C_{N2}
 \}
$. 
Then, for any $v_1, v_2\in \mathscr{X}_T$ and any $0<t<T$,
{we have}
\begin{align}
\label{CM1}
\|
 \Phi[{v}_1](t)
 -
 \Phi[{v}_2](t)
\|_{p}
\le \, &  
2
C_N 
\eta 
\mathrm{e}^{ 2 A \max \{ \beta_1, \beta_2 \} T}
D_{n,p}(T)
[
 v_1
 -
 v_2
]_t,
\\
\label{CM2}
 t^{\frac{n}{2} ( \frac{1}{p} - \frac{1}{q_*} )}
 \|
  \Phi[v_1](t)
  -
  \Phi[v_2](t)
 \|_{q_*}
 \le \, &
 2
 C_N 
 \eta 
 \mathrm{e}^{2A\max\{\beta_1, \beta_2\} T}
 D_{n,p}(T)
 [
  v_1
  -
  v_2
 ]_t, \\
\label{CM3}
 \|
  \Phi[v_1]
  -
  \Phi[v_2]
 \|_{\mathscr{X}_T}
 \le \, &
 4
 C_N 
 \eta 
 \mathrm{e}^{ 2 A \max \{ \beta_1, \beta_2 \} T}
 D_{n,p}(T)
 \|
  v_1
  -
  v_2
 \|_{\mathscr{X}_T}.
\end{align}
In fact, since 
$
 \Phi[{v}_1](t)
 -
 \Phi[{v}_2](t)
 =
 N[{v}_1-{v}_2, {v}_1](t)
 +
 N[{v}_2, {v}_1-{v}_2](t)
$,
we observe that 
\begin{equation}\label{es_8.13}
 \|
  \Phi[{v}_1](t)-\Phi[{v}_2](t)
 \|_r
 \le 
 \|
  N[{v}_1-{v}_2, {v}_1](t)
 \|_r
 +
 \|
  N[{v}_2, {v}_1-{v}_2](t)
 \|_r
\end{equation}
for $r = p, q_*$.
By the definition of $\mathscr{X}_T$,
we note that
\begin{equation}\label{ES815}
 [{v_1-v_2}]_T
 \le 
 \eta,
 \quad
 j=1, 2
\end{equation}
for {all} ${v}_1, {v}_2 \in \mathscr{X}_T$.
Combining estimate \eqref{cm1}{,} estimate \eqref{es_8.13} {and estimate \eqref{ES815}},
we obtain
\[
\begin{split}
 \|
  \Phi[v_1](t)-\Phi[v_2](t)
 \|_p
 \le & \, 
 \|
  N[v_1-v_2, v_1](t)
 \|_p
 +
 \|
  N[v_2, v_1-v_2](t)
 \|_p \\
 \le & \, 
  C_{N1}
  \mathrm{e}^{2A\max\{\beta_1, \beta_2\} T}
  D_{n,p}(T)
  (
   [v_1]_t
   +
   [v_2]_t
  )
  [
   v_1
   -
   v_2
  ]_t
  \\
  \le & \, 
  2
  C_N
  \eta
  \mathrm{e}^{ 2 A \max \{ \beta_1, \beta_2 \} T}
  D_{n,p}(T)
  [v_1-v_2]_t, 
  \quad 0<t<T.
\end{split}
\]
This gives the desired estimate \eqref{CM1}. 
In the same way, 
using estimate \eqref{cm2} together with inequality \eqref{es_8.13} and inequality \eqref{ES815} leads
to estimate \eqref{CM2}. 
By the inequality
$
 [ 
  {v}_1
  -
  {v}_2
 ]_{{t}}
 \le 
 \|
  {v}_1
  -
  {v}_2
 \|_{\mathscr{X}_T}
$
{(${v}_1, {v}_2 \in \mathscr{X}_T$), estimate \eqref{CM1} and estimate \eqref{CM2}, we get}
inequality \eqref{CM3}. {Moreover, by remarking that
the relation
\[
 \Phi[v](t)
 =
 \mathrm{e}^{t\mathcal{L}_A}v_0
 +[\Phi[v](t)-\Phi[0](t)],
\]
it holds from estimate \eqref{CM1} and estimate \eqref{CM2}
that 
\begin{align}
\label{cM1} 
 \|
  \Phi[v](t)
 \|_p
 \le \, &
 \|
  \mathrm{e}^{t\mathcal{L}_A}v_0
 \|_p
 +
 2
 C_N
 \eta 
 \mathrm{e}^{2A\max\{\beta_1, \beta_2\}T}
 D_{n,p}(T)
 [v]_t, \\
\label{cM2} 
 t^{\frac{n}2(\frac1p-\frac1{q_*})}
 \|
  \Phi[v](t)
 \|_{q_*} 
 \le \, & 
 t^{\frac{n}{2} (\frac{1}{p} - \frac{1}{q_*} )}
 \|
  \mathrm{e}^{t\mathcal{L}_A}v_0
 \|_{q_*}
 +
 2
 C_N
 \eta 
 \mathrm{e}^{2A\max\{\beta_1, \beta_2\}T}
 D_{n,p}(T)[v]_t
\end{align}
for {any} $v\in \mathscr{X}_T$ and {any}  $0<t<T$.
}

Finally we show that there is a number $T_0 \in (0,1]$ such that $\Phi$ is {a} contraction {mapping on $\mathscr{X}_{T_0}$}.
{Let $n\ge 1$ and} $\max \{ 1, n \} < p < \infty${, and} 
take $\eta, B > 0$ as 
$
 \eta 
 =
 B
 =
 2
 C_{L}
 \mathrm{e}^{2A\max\{\beta_1, \beta_2\}}
 \|
  {v}_0
 \|_p
$.
Then, by vitue of the relation 
$
 D_{n,p}(T)
 =
 T^{\frac12(1-\frac{n}p)}
$, 
estimates \eqref{lcm1}, \eqref{ES815} and \eqref{cM1} give
\begin{equation}\label{Es89}
\sup_{0\le t\le T}
 \|
  \Phi[v](t)
 \|_{p}
 \le 
 C_L
 \mathrm{e}^{2A\max\{\beta_1, \beta_2\}}
 \|
  v_0
 \|_p
 +
 2
 C_{N}
 \mathrm{e}^{2A\max\{\beta_1, \beta_2\}}
 B^2
 T^{\frac12(1-\frac{n}p)}.
\end{equation}
{Here c}hoose {$T_0 \in (0,1]$} as 
$
 2C_{N}\mathrm{e}^{2A\max\{\beta_1, \beta_2\}}
 T^{\frac12{(1-\frac{n}p)}}
 \le
 1/(4B)
$ 
for given $B>0$. 
Then, from inequality \eqref{Es89} we have 
\begin{equation}\label{Es810}
\sup_{0 \le t\le {T_0}}\|\Phi[{v}](t)\|_p
\le 
\frac{B}{2}
+
\frac{B}{{4}}
{<}
B.
\end{equation}
Similarly, it follows from estimate{s} \eqref{lcm2}, \eqref{ES815} and \eqref{cM2} that
\begin{equation}\label{Es813}
\begin{split}
 \sup_{0<t\le {T_0}}
 t^{\frac{n}{2} ( \frac{1}{p} - \frac{1}{q_*} )}
 \|
  \Phi[{v}](t)
 \|_{q_*}
 \le &
 C_L
 \mathrm{e}^{2A\max\{\beta_1, \beta_2\}}
 \| 
  {v}_0
 \|_p
 +
 2
 C_{N}
 \mathrm{e}^{2A\max\{\beta_1, \beta_2\}}
 \eta^2
 {T_0}^{\frac12(1-\frac{n}p)} \\
 \le &
 \frac{\eta}{2}
 +
 \frac{\eta}{{4}}
 {<}
 \eta.
 \end{split}
\end{equation}
Combining 
$
 2
 C_N
 \mathrm{e}^{2A\max\{\beta_1, \beta_2\}}
 {T_0}^{\frac12(1-\frac{n}p)}
 \le 
 1/(4B)
$
with {estimate}
\eqref{CM3} shows {that}
\begin{equation}\label{ES822}
 \|
  \Phi[{v}_1]
  -
  \Phi[{v}_2]
 \|_{\mathscr{X}_{T_*}}
 \le 
 \frac{1}{2}
 \| 
  {v}_1
  -
  {v}_2
 \|_{\mathscr{X}_{T_*}}.
\end{equation}
Therefore, by estimate{s} \eqref{Es810}, \eqref{Es813} and \eqref{ES822},  
there exists {a number $T_0 > 0$} such that 
$\Phi[{v}] \in \mathscr{X}_{{T_0}}$ for all ${v} \in \mathscr{X}_{{T_0}}$ and $\Phi : \mathscr{X}_{{T_0}} \to \mathscr{X}_{{T_0}}$ is {a} contraction {mapping}.
Hence we conclude that there exists a unique {fixed point} ${v}\in \mathscr{X}_{{T_0}}$ such that
\[
 v
 =
 \Phi[v]
\]
by the Banach fixed point {theorem}.

{Let} $n \ge 2$ and $n/2 \le p \le n$ or $n=p=1$. 
Then $D_{n,p}(T)=1$ {holds}. 
{By choosing} $B>0$ as
$
 B
 =
 2
 C_{L}
 \mathrm{e}^{2A\max\{\beta_1, \beta_2\}}
 \|
  {v}_0
 \|_p
${,} inequalit{ies} \eqref{lcm1}, \eqref{ES815} and \eqref{cM1} lead to
\begin{equation}\label{Es891}
 \sup_{0\le t\le T}
  \|
   \Phi[{v}](t)
  \|_{p}
  \le 
  \|
   \mathrm{e}^{t\mathcal{L}_A}{v}_0
  \|_{p}
  +
  2
  C_{N}
  \mathrm{e}^{2A\max\{\beta_1, \beta_2\}}
  \eta^2
  \le
  \frac{B}{2}
  +
  2
  C_{N}
  \mathrm{e}^{2A\max\{\beta_1, \beta_2\}}
  \eta^2.
\end{equation}
{Here, if the constant} $\eta_1 > 0$ is chosen {so} that
$
 2C_{N}
 \mathrm{e}^{2A\max\{\beta_1, \beta_2\}}
 \eta^2
 \le
 B/2
$
for given $B > 0${, then} inequality \eqref{Es891} gives
\[
 \sup_{0 \le t \le T}
  \| 
   \Phi[{v}](t)
  \|_p
  \le 
  \frac{B}{2}
  +
  \frac{B}{2}
  =
  B.
\]
{Moreover,} take
$
 \eta_2 > 0
$
with 
$
 2
 C_N
 \mathrm{e}^{2A\max\{\beta_1, \beta_2\}}
 \eta_2
 \le 
 1/4
$
{and d}efine
$
 \eta
 \coloneqq
 \min\{ \eta_1, \eta_2 \}
$
.
Since  
$
 C_0^\infty(\mathbb{R}^n)
$
is dense in $L^p(\mathbb{R}^n)$, 
there {exists} $\widetilde{{v}}_{0} \in C_0^\infty(\mathbb{R}^n)$ such that
\begin{equation}\label{ap_1}
 \|
  \widetilde{{v}}_{0}
  -
  {v}_0
 \|_p
 \le
 \frac{\eta}{4C_L\mathrm{e}^{2A\max\{\beta_1, \beta_2\}}},
 \end{equation}
where
$C_L > 0$ is given by relation \eqref{c84}. 
{So, picking $0<T_{\ast\ast}\le 1$ as 
 $
 C_L
 \mathrm{e}^{2A\max\{\beta_1, \beta_2\}}
 T_{**}^{\frac14}
 \|
  \widetilde{{v}}_0
 \|_{q_*}
 \le 
 \eta/4
 $, 
 we see from inequalities \eqref{lcm1}, \eqref{lcm2} and \eqref{ap_1} that} 
\begin{equation}\label{les825}
\begin{split}
 t^{\frac{n}{2} ( \frac{1}{p} - \frac{1}{q_*} )}
 \|
  \mathrm{e}^{t\mathcal{L}_A}
  {v}_0
 \|_{q_*}
 \le \, &
 t^{\frac{n}{2} ( \frac{1}{p} - \frac{1}{q_*} )}
 \|
  \mathrm{e}^{t\mathcal{L}_A}
  \widetilde{{v}}_{0}
 \|_{q_*}
 +
 t^{\frac{n}{2} ( \frac{1}{p} - \frac{1}{q_*} )}
 \|
  \mathrm{e}^{t\mathcal{L}_A}
  (
   {v}_0-\widetilde{{v}}_{0}
  )
 \|_{q_*}\\
 \le \, &
 C_L
 \mathrm{e}^{2A\max\{\beta_1, \beta_2\}}
 T_{**}^{\frac14}
 \|
  \widetilde{{v}}_0
 \|_{q_*}
 +
 C_L
 \mathrm{e}^{2A\max\{\beta_1, \beta_2\}}
 \|
  \widetilde{v}_{0}-{v}_0
 \|_p\\
 \le \, &
 \frac{\eta}{4}
 +
 \frac{\eta}{4}\\
 = \, &
 \frac{\eta}{2}
\end{split}
\end{equation}
for all $0 < t < T_{**}$. {Here we have used the relation}  
\[
 \frac{n}{2} 
 \bigg( \frac{1}{p} - \frac{1}{q_*} \bigg)
 =
 \frac{n}{2} 
 \bigg( \frac{1}{p} - \Big( \frac1p - \frac1{2n}\Big) \bigg)
 =
 \frac14
\]
{in the case $n \ge 2$ and ${n}/2 \le p \le n$ or $n=p=1$. Therefore,} for {any} $0 < t < T_{**}$,
inequalit{ies} \eqref{cM2} and \eqref{les825} show {that}
\[
 \sup_{0<t\le T_{**}}
 t^{\frac{n}{2} ( \frac{1}{p} - \frac{1}{q_*} )}
 \|
  \Phi[{v}](t)
 \|_{q_*}
 \le 
 \frac{\eta}{2}
 +
 2
 C_{N}
 \mathrm{e}^{2A\max\{\beta_1, \beta_2\}}
 \eta \cdot \eta 
 \le 
 \frac{\eta}{2}
 +
 \frac{\eta}{4}
 <
 \eta.
\]
According to the definition of $\eta_2$, that is{,}  
$
 2
 C_N
 \mathrm{e}^{2A\max\{\beta_1, \beta_2\}}
 \eta_2
 \le 
 1/4
$,
we obtain 
\[
 \|
  \Phi[{v}_1]
  -
  \Phi[{v}_2]
 \|_{\mathscr{X}_{T_{**}}}
 \le 
 \frac{1}{2}
 \|
  {v}_1
  -
  {v}_2
 \|_{\mathscr{X}_{T_{**}}}
\]
from 
inequality \eqref{CM3}.
As a result, we conclude that 
$
 \Phi : \mathscr{X}_{T_{**}} \to \mathscr{X}_{T_{**}}
$
is {a} contraction {mapping for some $T_{\ast\ast}>0$},
which finally gives the existence of a fixed point $v \in \mathscr{X}_{T_{**}}$ such that $v = \Phi[v]$.
\end{proof}

\subsection{Uniqueness of solution to {problem} \eqref{eq;lin}}
The uniqueness of solution to problem \eqref{eq;lin} will be shown as follows:
\begin{proposition}\label{prop_82} 
Let $n \ge 1$, $A > 0 $, $0 < T \le \infty$ and 
\[
 \begin{cases} 
 p=1 &\text{if $n=1$}, \\
 \dfrac{n}{2}<p<\infty &\text{if $n\ge 2$}.
 \end{cases}
\]
Also, assume that
${v}_1, {v}_2 \in C([0,T);L^p(\mathbb{R}^n))$ are mild solutions to problem \eqref{eq;lin} corresponding to the initial data ${v}_0 \in L^p(\mathbb{R}^n)$ and that 
\begin{equation}
 \label{u10}
 \sup_{0 \le t < T}
 (\| 
  v_1(t)
 \|_p+\| 
  v_2(t)
 \|_p)
 <
 \infty.
\end{equation}
Then ${v}_1={v}_2$.
\end{proposition}

\begin{proof} 
 We give only the proof in the case $n\ge 2$ and ${n}/{2}<p<\infty$ because the one in the case $n=p=1$ is done in a similar fashion.
 Fix $T^\ast\in (0,T)$ and let ${v}_1, {v}_2\in C([0,T^\ast];L^p({\mathbb R}^n))$ be mild solutions to problem \eqref{eq;lin} corresponding to the initial data 
 ${v}_0 \in L^p(\mathbb{R}^n)$ with \eqref{u10}. 
 Then,
 we note that for all $0 < t \le T^\ast$,
\[
 {v}_1(t)
 -
 {v}_2(t)
 =
 N[{v}_1-{v}_2, {v}_1](t)
 +
 N[{v}_2, {v}_1-{v}_2](t).
\]
where the function $N[\cdot, \cdot]$ is defined by relation \eqref{pp-nuv}.
Applying {the} arguments {similar to those in} the proof of Proposition \ref{prop;812} leads to
\begin{equation}\label{es_832}
\begin{split}
 \|
  {v}_1&(t)
  -
  {v}_2(t)
 \|_p\\
 \le \, &
  C 
  \int_0^t 
   (t-s)^{-\frac{n}2(\frac1q-\frac1p)-\frac12}
   \mathrm{e}^{ 2A\max\{\beta_1, \beta_2\}(t-s)}
   \|
    [
     {v}_1(s)
     -
     {v}_2(s)
    ]
    \nabla K*{v}_1(s)
   \|_q
  \, \mathrm{d}s \\
 & \, +
 C
 \int_0^t
  (t-s)^{-\frac{n}2(\frac1q-\frac1p)-\frac12}
  \mathrm{e}^{2A\max\{\beta_1, \beta_2\}(t-s)}
  \|
   {v}_2(s)
   \nabla 
   K*
   [ 
    {v}_1
    -
    {v}_2
   ](s)
  \|_q 
  \, \mathrm{d}s\\
 \le \, &
 C
 \mathrm{e}^{2A \max\{\beta_1, \beta_2\} t}
 \int_0^t
  (t-s)^{-\frac{n}2(\frac1q-\frac1p)-\frac12}
  \|
   {v}_1(s)
   -
   {v}_2(s)
  \|_p
  \|
   \nabla K*{v}_1(s)
  \|_{q^*}
  \, \mathrm{d}s\\
 & \, +
 C
 \mathrm{e}^{2A \max \{ \beta_1, \beta_2 \} t}
 \int_0^t 
  (t-s)^{-\frac{n}2(\frac1q-\frac1p)-\frac12}
  \|
   {v}_2(s)
  \|_p
  \|
   \nabla 
   K*
   [
    {v}_1
    -
    {v}_2
   ](s)
  \|_{q^*}
  \, \mathrm{d}s
\end{split}
\end{equation} 
for all $0<t\le T^\ast$, where 
the exponents $p$, $q$ and $q^\ast$ satisfy
\begin{equation}\label{ex_pq0}
 1 \le q \le p, \quad 
 \frac{1}{2}
 -
 \frac{n}{2}
  \left(
   \frac{1}{q}
   -
   \frac{1}{p}
  \right)
 >0
\end{equation}
and
\[
 \frac{1}{q^{*}}
 =
 \frac{1}{q}
 -
 \frac{1}{p}.
\]
If
\[
 \frac{1}{q_1}
 \coloneqq 
 \frac{1}{q^*}
 -
 \frac{1}{p}
 +
 1
 =
 \frac{1}{q}
 -
 \frac{2}{p}
 +
 1
\]
satisfies 
\[
 1-\frac{1}{n}
 \le 
 \frac{1}{q_1}
 \le 
 1.
\]
Then 
\begin{equation}\label{ex_pq1}
 0
 \le 
 \frac{2}{p}
 -
 \frac{1}{q}
 \le 
 \frac{1}{n}.
\end{equation}
Hence, for 
fixed a $p \in ({n}/2, \infty)$, 
there exists {an exponent} $q$ such that 
inequalit{ies} \eqref{ex_pq0} and \eqref{ex_pq1}
hold{, and then we make use of Lemma \ref{lem;be} and} estimate \eqref{es_832} {to obtain}
\begin{equation}\label{es_833}
 \begin{split}
  & \,
  \|
  {v}_1(t) 
  -
  {v}_2(t)
 \|_p\\
 \le \, & 
 C
 \mathrm{e}^{2A\max\{\beta_1, \beta_2\} t}
 \int_0^t
  (t-s)^{-\frac{n}2(\frac1q-\frac1p)-\frac12}
  \|
   \nabla K
  \|_{L_w^{q_1}}
  \|
   {v}_1(s)
  \|_p
  \|
   {v}_1(s)
   -
   {v}_2(s)
  \|_p
 \, \mathrm{d}s \\
 & \,  +
  C
  \mathrm{e}^{2A\max\{ \beta_1, \beta_2 \} t}
   \int_0^t
    (t-s)^{-\frac{n}2(\frac1q-\frac1p)-\frac12}
     \|
      \nabla K
     \|_{L_w^{q_1}}
     \|
      {v}_2(s)
     \|_p
     \|
      {v}_1(s)
      -
      {v}_2(s)
     \|_p
     \, \mathrm{d}s \\
  \le \, &
  C
   \left(
    \sum_{j=1}^2
    \sup_{0 \le t \le T^\ast}
    \|
     {v}_j(t)
    \|_p
   \right)\\
   &  \qquad \qquad \qquad \times
   \mathrm{e}^{2A\max\{\beta_1, \beta_2\} t}
   \left(
    \sup_{0 \le s \le t}
     \|
      {v}_1(s)
      -
      {v}_2(s)
     \|_p
   \right)
   \int_0^t
    (t-s)^{-\frac{n}2(\frac1q-\frac1p)-\frac12}
   \, \mathrm{d}s \\
  \le \, &
   C
    \left(
     \sum_{j=1}^2
     \sup_{0 \le t \le T^\ast}
      \|
       {v}_j(t)
      \|_p
    \right) \\
    &  \qquad \qquad \qquad \times
    t^{\frac12-\frac{n}2(\frac1q-\frac1p)}
    \mathrm{e}^{2A\max\{\beta_1, \beta_2\} t}
    \sup_{0 \le s \le t}
     \|
      {v}_1(s)
      -
      {v}_2(s)
     \|_p
 \end{split}
\end{equation}
{for all $0<t\le T^\ast$.} Taking $\tau_0 \in (0,T^*)$ so that 
\begin{equation}\label{c_834}
 0
 <
 {C}
 \left(
  \sum_{j=1}^2
  \sup_{0 \le t {\le T^\ast}}
   \|
    {v}_j(t)
   \|_p
 \right)
 {\tau_0}^{\frac12-\frac{n}2(\frac1q-\frac1p)}
 \mathrm{e}^{2A\max\{\beta_1, \beta_2\} \tau_0}
 <
 1,
\end{equation}
we have 
\[
 \sup_{0 \le s \le \tau_0}
 \|
  {v}_1(s)
  -
  {v}_2(s)
 \|_p
 =
 0
\]
{by} inequality \eqref{es_833} and {condition \eqref{c_834}}.
Therefore ${v}_1={v}_2$ in $C([0,\tau_0];L^p(\mathbb{R}^n))$.
If $\tau_0 = T^*$ for a fixed $T^*>0$, 
then we {get} the conclusion of Proposition \ref{prop_82}. 
{On the other hand, if} $\tau_0 < T^*$, {the argument used in the proof of estimate \eqref{es_833}} works {well} after taking $L^p$-norm in relation \eqref{es_832}. Indeed, for {any $\tau_0<t\le T^\ast$}, we have
\begin{equation}\label{es_835}
\begin{split}
 &
  \sup_{0 \le s \le t}
  \|
   \chi_{\{ s > \tau_0 \}}(s)
   [
    {v}_1
    -
    {v}_2
   ]
   (s)
 \|_p
 \\
& \le 
 {C}
 M({T^\ast})
 (t-\tau_0)^{\frac12-\frac{n}2(\frac1q-\frac1p)}
 \mathrm{e}^{2A\max\{\beta_1, \beta_2\}(t-\tau_0)}
 \sup_{0 \le s \le t}
 \|
  \chi_{\{ s > \tau_0 \}}(s)
  [ 
   {v}_1
   -
   {v}_2
  ](s)
 \|_p,
\end{split}
\end{equation}
where
\[
 M({T^\ast})
 \coloneqq 
 \sum_{j=1}^2
 \sup_{0 \le t \le T^\ast}
  \| 
   {v}_j(t)
  \|_p, 
  \quad 
  \chi_{\{ t > \tau_0 \}}(t)
  =
\begin{cases}
 1, & \, \tau_0 < t, \\
 0, & \, 0 < t \le \tau_0.
\end{cases}
\]
{Here set}
$
 \tau_1 
 = 
 2 
 \tau_0
$. 
Then it follows from \eqref{c_834} that
\[
 \begin{split}
 0
 < \, & 
 C
 M({T^\ast})
 (\tau_1-\tau_0)^{\frac12-\frac{n}2(\frac1q-\frac1p)}
 \mathrm{e}^{2A\max\{\beta_1, \beta_2\}(\tau_1-\tau_0)}\\
 = \, &
 C
 M({T^\ast})
 \tau_0^{\frac12-\frac{n}2(\frac1q-\frac1p)}
 \mathrm{e}^{2A\max\{\beta_1, \beta_2\}\tau_0}
 <
 1,
 \end{split}
\]
{which together with inequality \eqref{es_835} implies that}
\[ 
 \sup_{0 \le s \le \tau_1}
 \|
  \chi_{\{ s > \tau_0\}}(s)
  [
   {v}_1
   -
   {v}_2
  ](s)
 \|_{p}
 =
 0.
\] 
{As a result, t}his shows that ${v}_1={v}_2$ in $C([0,\tau_1];L^p(\mathbb{R}^n))$. 

{By r}epeating the previous argument {above, we can choose $m\in {\mathbb N}$ so that} $\tau_\ell = 2\tau_{\ell -1}$ for {each} $\ell = 1, 2, \cdots, m$ and $\tau_{m-1} < T^* < \tau_{m}${, and} conclude that ${v}_1={v}_2$ in $C([0,T^*];L^p(\mathbb{R}^n))$. {Hence, Proposition \ref{prop_82} is proved since $T^\ast>0$ is arbitrary}.
\end{proof}

\subsection{Global-in-time solutions to {problem} \eqref{eq;lin}}

The following theorem gives so called the small data global existence of solutions to problem \eqref{eq;lin} in order to show Theorem \ref{thm;1}.

\begin{theorem}\label{th81}
{Assume that} $p$ satisfies 
\begin{equation}\label{ex;p}
 \begin{cases}
  p=1 &\text{{if} $n=1$}, \\
  p \in \Big( \dfrac{n}2, n \Big] &\text{{if} $n \ge 2$.}
 \end{cases}
\end{equation}
If the constant $A>0$ satisfies \eqref{eq;cond-A}, then
there exists 
$
 \varepsilon_0
 >
 0
$
such that the following statement holds: 
For any ${v}_0 \in L^p(\mathbb{R}^n)$ with 
$
 \|
  {v}_0
 \|_p
 <
 \varepsilon_0
$,
there exists a unique global-in-time {mild} solution 
$
 {v} 
 \in 
 C([0,\infty);L^p(\mathbb{R}^n)) 
 \cap 
 C((0,\infty);L^{2p}(\mathbb{R}^n))
$ 
to {problem \eqref{eq;lin}} such that 
\begin{align}
 & 
 \label{udn}
 \sup_{t>0}
 \| 
  {v}(t)
 \|_p
 +
 \sup_{t>0}
 t^{\frac{n}2(\frac1p-\frac1q)}
 \|
  {v}(t)
 \|_{q}
 \le 
 C_0
 \|
  {v}_0
 \|_p
\end{align}
for some $C_0>0$ and all ${n} < q \le 2p$.
\end{theorem}

\begin{proof}
{Let $N[\cdot,\cdot]$ and $\Phi$ be the same notations as \eqref{pp-nuv} and \eqref{Ieu}, respectively.} By the $L^p$-$L^q$ type estimate of the {analytic} semigroup $\mathrm{e}^{t\mathcal{L}_A}$ {in \eqref{eq;es_sa}}, 
we have 
\[
\begin{split}
 \|
  \mathrm{e}^{t\mathcal{L}_A}
  {v}_0
 \|_{p}
 \le  \, &
 \|
  {v}_0
 \|_p, \\
 \|
  \mathrm{e}^{t\mathcal{L}_A}
  {v}_0
 \|_{{2p}}
 \le \, &
 C
 t^{-{\frac{n}{4p}}}
 \|
  {v}_0
 \|_{p}
\end{split}
\]
for {any $t>0$}.
Thus
\begin{equation}\label{es_ln1}
 \|
  \mathrm{e}^{{\cdot}\mathcal{L}_A}
  {v}_0
 \|_{\mathscr{X}}
 \le 
 C_L
 \|
  {v}_0
 \|_{p}
\end{equation}
holds, where
\begin{align*}
  \mathscr{X}
  \coloneqq \, &
  \{
   {v}
   \in 
   C([0,\infty);L^{p}(\mathbb{R}^n)) \cap C((0,\infty);L^{{2p}}(\mathbb{R}^n)) 
   \mid 
   \|
    {v}
   \|_{\mathscr{X}}
   \le 
   2
   C_L
   \|
    {v}_0
   \|_p
  \}, \\
  \|
   {v} 
  \|_{\mathscr{X}}
  \coloneqq \, & 
  \sup_{t>0}
   \|
    {v}(t)
   \|_{p}
   +
   \sup_{t>0}
    t^{{\frac{n}{4p}}}
    \|
     {v}(t)
    \|_{{2p}}. 
\end{align*}

It follows from the definition of the norm in $\mathscr{X}$ that 
\begin{equation}\label{norm;u}
\sup_{t>0}t^{\frac{n}{2}(\frac{1}{p}-\frac{1}{q})}\|v(t)\|_q\le \left(\sup_{t>0}\|v(t)\|_p\right)^{\frac{2p}{q}-1}\left(\sup_{t>0}t^{\frac{n}{4p}}\|v(t)\|_{2p}\right)^{2-\frac{2p}{q}}\le \|v\|_{\mathscr{X}}
\end{equation}
for all $v\in \mathscr{X}$ and all $p<q<2p$. Therefore, noticing that
the exponent $q_\ast$ denoted by \eqref{dfn;pq} satisfies $
p<q_\ast\le 2p
$
due to \eqref{ex;p},
we use \eqref{norm;u} to get 
\begin{equation}\label{norm;u1}
 \|
  v_1(t)
 \|_{q_*}
 \|
  v_2(t)
 \|_{q_*}
 \le 
 t^{-\frac{1}{2}}
 \|
  v_1
 \|_{\mathscr{X}}
 \|
  v_2
 \|_{\mathscr{X}} 
\end{equation}
for all $v_1, v_2\in {\mathscr X}$ and all $t>0$.

Define $r_\ast$, $q_1$ as 
\[
 \frac{1}{r_\ast}
 =
 \frac{1}{2n}, \quad 
 \frac{1}{q_1}
 =
 -\frac{1}{p}+\frac{1}{n}+1,
\]
respectively. Then 
\[
 \frac{1}{q_1}
 \in 
 \Big(
  \frac{n-1}{n}, 1
 \Big] 
 \quad \text{if $n\ge 2$
} \quad\text{or} \quad \frac{1}{q_1}=1 \quad \text{if $n=1$}
\]
because of \eqref{ex;p} and also
\[
\frac{1}{p}=\frac{1}{q_\ast}+\frac{1}{r_\ast}, \quad \frac{1}{r_\ast}=\frac{1}{q_1}+\frac{1}{q_\ast}-1. 
\] 
Hence, applying Lemma \ref{lem;be}, the $L^p$-$L^q$ type estimate {of the analytic semigroup $\mathrm{e}^{t\mathcal{L}_A}$ in \eqref{eq;es_nsa}}, {estimate \eqref{norm;u1}}, 
H\"{o}lder's inequality and Hausdorff-Young's inequality implies that {for all $v_1$, $v_2\in {\mathscr X}$ and all $t>0$}
\[
\begin{split}
 \|
  N[{v_1,v_2}](t)
 \|_{p}
 \le \, &
 C
  \int_0^t
  (t-s)^{-\frac12}
  \|
   {v}_1(s)
   \nabla K*{v}_2(s)
  \|_{p}
  \, \mathrm{d}s\\
 \le \, &
 C
 \int_0^t
  (t-s)^{-\frac12}
  \|
   {v}_1(s)
  \|_{q_*}
  \|
   \nabla K*{v}_2(s)
  \|_{r_*}
  \, \mathrm{d}s
 \\
 \le \, &
 C
  \int_0^t
   (t-s)^{-\frac12}
   \|
    \nabla K
   \|_{q_1}
   \|
    {v}_1(s)
   \|_{q_*}
   \|
    {v}_2(s)
   \|_{q_*}
  \, \mathrm{d}s
   \\
 \le \, &
 C\|
  {v}_1
 \|_{\mathscr{X}}
 \|
  {v}_2
 \|_{\mathscr{X}}
 \int_0^t
  (t-s)^{-\frac12}
  s^{-{\frac{1}{2}}}
  \, \mathrm{d}s\\
 \le \, & 
 C
  B
  \left(
   \frac{1}{2}, \frac{1}{2}
  \right)
  \|
   {v}_1
  \|_{\mathscr{X}}
  \|
   {v}_2
  \|_{\mathscr{X}}, 
\end{split}
\]
where {$B(\cdot,\cdot)$ is the Beta function}. {This gives}
\begin{equation}\label{es_nn}
 \sup_{t>0}
 \|
  N[{v}_1,{v}_2](t)
 \|_{p}
 \le 
 C_1
 \|
  {v}_1
 \|_{\mathscr{X}}
 \|
  {v}_2
 \|_{\mathscr{X}}
\end{equation}
for all ${v}_1$, ${v}_2 \in \mathscr{X}$ and some $C_1 > 0$ independent of $t > 0$.
{Similarly}, we have
\begin{equation}\label{es_n112n}
 \sup_{t>0}
 t^{{\frac{n}{4p}}}
 \|
  N[{v}_1,{v}_2](t)
 \|_{{2p}}
 \le 
 C_2
 \|
  {v}_1
 \|_{\mathscr{X}}
 \|
  {v}_2
 \|_{\mathscr{X}}
\end{equation}
for all ${v}_1, {v}_2 \in \mathscr{X}$ and some $C_{{2}} > 0$ independent of $t>0$.
{Thus}, combining estimate \eqref{es_nn} with \eqref{es_n112n} leads to 
\begin{equation}\label{es_848}
 \|
  N[{v}_1,{v}_2]
 \|_{\mathscr{X}}
 \le 
 C_{*}
 \|
  {v}_1
 \|_{\mathscr{X}}
 \|
  {v}_2
 \|_{\mathscr{X}}
\end{equation}
for {any $v_1$, $v_2\in {\mathscr X}$ and} some $C_* > 0$.

We will now show that 
\begin{align}
 & \, 
 \label{es_cNN}
 \|
  \Phi[{v}_1]-\Phi[{v}_2]
 \|_{\mathscr{X}}
 \le 
 4C_* C_L
 \|
  {v}_0
 \|_p
 \|
  {v}_1
  -
  {v}_2
 \|_{\mathscr{X}}, \\
 & \, 
 \label{es_cNNN}
 \|
  \Phi[{v}]
 \|_{\mathscr{X}}
 \le 
 {C_L\|{v}_0\|_{p}(1+4C_\ast C_L\|{v}_0\|_p)}
 +
 2
 C_* 
 C_L
 \|
  {v}_0
 \|_{p}
 \|
  {v}
 \|_{\mathscr{X}}
\end{align}
for all ${v}, {v}_1, {v}_2 \in \mathscr{X}$, 
where $C_* > 0$ is given by inequality \eqref{es_848}.
By inequality \eqref{es_848}, we have
\[
\begin{split}
 \|
  \Phi[{v}_1]
  -
  \Phi[{v}_2]
 \|_{\mathscr{X}} 
 \le \, & 
 \|
  N[{v}_1-{v}_2,{v}_1](t)
 \|_{\mathscr{X}}
 +
 \|
  N[{v}_2,{v}_1-{v}_2](t)
 \|_{\mathscr{X}}\\
 \le \, & 
 C_*
 (
  \|
   {v}_1
  \|_{\mathscr{X}}
  +
  \|
   {v}_2
  \|_{\mathscr{X}}
 )
 {\|
   v_1-v_2
  \|_{\mathscr{X}}}\\
 \le \, & 
 4
  C_* 
  C_L
  \|
   {v}_0
  \|_p
  \|
   {v}_1-{v}_2
  \|_{\mathscr{X}}{,}
 \end{split}
\]
which shows inequality \eqref{es_cNN}. 
Moreover, by {inequality \eqref{es_ln1}} and inequality \eqref{es_848}, 
\[
\begin{split}
 \|
  \Phi[{v}]
 \|_{\mathscr{X}}
 \le \, &
 \|
  \mathrm{e}^{t\mathcal{L}_A}
  {v}_0
 \|_{\mathscr{X}}
 +
 \|
  N[{v},{v}]
 \|_{\mathscr{X}}\\
 \le \, &
 {C_L\|v_0\|_p}
 +
 C_*
 \|
  {v}
 \|_{\mathscr{X}}^2\\
 \le \, &
 {C_L\|v_0\|_p(1+4C_\ast C_L\|v_0\|_p)}
\end{split}
\]
holds for {all} ${v}\in \mathscr{X}$.
This is the desired inequality \eqref{es_cNNN}.

Thus{, if ${v}_0 \in L^p(\mathbb{R}^n)$ satisfies} 
$
 4
 C_*
 C_L
 \| 
  {v}_0
 \|_p
 <
 1
$, {then $\Phi[v]\in 
{\mathscr X}$ for all $v\in {\mathscr X}$,}
and {furthermore} $\Phi : \mathscr{X} \to \mathscr{X}$ is {a} contraction {mapping}. {As a result, a}ccording to the Banach fixed point theorem,
there exists {a unique fixed point} $
 {v}
 \in
 \mathscr{X}
$
such that
${v}$ is a mild solution to problem \eqref{eq;lin} {satisfying
$
\|
  v
 \|_{\mathscr{X}}
 \le 
 2
 C_{L}
 \|
  v_0
 \|_{p}
$.} 

We will finally prove {only} inequality \eqref{udn} {in the case $n/2 < p \le n < q \le 2p$} and $n \ge 2$ {since it readily follows from $\|v\|_{{\mathscr X}}\le 2C_L\|v_0\|_1$ in the case $n=p=1$. We employ Proposition \ref{prop;lp-lq-li}, H\"{o}lder's inequality, Young's inequality and \eqref{norm;u1} to obtain} 
\[
\begin{split}
 \|
  {v}(t)
 \|_q
 \le \, &
 \|
  \mathrm{e}^{t\mathcal{L}_A}
  {v}_0
 \|_q
 +
 \int_0^t
  \|
   \nabla \cdot \mathrm{e}^{(t-s)\mathcal{L}_A}
   (
    {v}\nabla K*{v}
   )(s)
  \|_q
  \, \mathrm{d}s \\
 \le \, & 
 C
 t^{-\frac{n}2(\frac1p-\frac1q)}
 \|
  {v}_0
 \|_p
 +
 C
 \int_0^t
  (t-s)^{-\frac{n}2(\frac1p-\frac1q)-\frac12}
  \| 
   {v}(s)
   \nabla K*{v}(s)
  \|_p
  \, \mathrm{d}s\\
 \le \, & 
 C
  t^{-\frac{n}2(\frac1p-\frac1q)}
   \|
    {v}_0
   \|_p
   +
   C
   \int_0^t
    (t-s)^{-\frac{n}2(\frac1p-\frac1q)-\frac12}
    \|
     {v}(s)
    \|_{q_*}^2
    \, \mathrm{d}s\\
 \le \, & 
 C
 t^{-\frac{n}2(\frac1p-\frac1q)}
 \|
  {v}_0
 \|_p
 +
{4CC_L^2\|v_0\|_p^2}
 \int_0^t
  (t-s)^{-\frac{n}2(\frac1p-\frac1q)-\frac12}s^{-\frac12}
 \, \mathrm{d}s\\
 \le \, & 
 C
 t^{-\frac{n}2(\frac1p-\frac1q)}
 \|
  {v}_0
 \|_p
 +
 4
 C
 C_L^2\|
   {v}_0
  \|_{p}^2
 B
  \bigg(
   \frac{1}{2}
   -
   \frac{n}{2}
   \Big(
    \frac{1}{p}
    -
    \frac{1}{q}
   \Big), 
   \frac{1}{2}
  \bigg)
  t^{-\frac{n}2(\frac1p-\frac1q)}
\end{split}
\]
for $v\in \mathscr{X}$, where
\[
 \frac{1}{2}
 -
 \frac{n}{2}
 \bigg(
  \frac{1}{p}-\frac{1}{q}
 \bigg)
 = 
 \frac{1}{2}
 -
 \frac{n}{2p}
 +
 \frac{n}{2q}
 \ge 
 \frac{1}{2}
 -
 \frac{n}{4p}
 >
 0.
\]
Therefore, 
this together with 
$
 \sup_{t>0}
 \|
  v(t)
 \|_p
 \le 
 2
 C_L
 \|
  v_0
 \|_p
 \ 
 ( v \in \mathscr{X})
$ 
yields the desired estimate \eqref{udn} in the case ${n}/{2}<p\le n<q\le 2p$ and $n\ge 2$.
\end{proof}
\subsection{Instability of constant steady states}
We will finally discuss instability of constant steady states when the constant $A>0$ 
satisfies $A>A^\ast$, where the constant $A^* > 0$ is defined by \eqref{eq;cond-A-th}, that is,
\[
A^{\ast} 
 \coloneqq 
 \begin{cases}
\displaystyle  \frac{\lambda_1 \lambda_2}{\beta_1\lambda_2 - \beta_2\lambda_1}, &\displaystyle\frac{\beta_1}{\beta_2} \ge \left(\frac{\lambda_1}{\lambda_2}\right)^2, \quad \frac{\beta_1}{\beta_2}>\frac{\lambda_1}{\lambda_2}, \vspace{2mm}\\
\displaystyle \frac{\lambda_1-\lambda_2}{(\sqrt{\beta_1}-\sqrt{\beta_2})^2}, &\displaystyle \frac{\beta_1}{\beta_2} > 1,\quad \frac{\beta_1}{\beta_2} < \left(\frac{\lambda_1}{\lambda_2}\right)^2.
 \end{cases}
\]

\begin{theorem}\label{t83}
Suppose that $p$ satisfies
\[
 \begin{cases}
  p \in (1,\infty) &\text{{if} $n=1$}, \\
  p \in \Big( \dfrac{n}2, \infty \Big) &\text{{if} $n \ge 2$}. 
 \end{cases}
\]
{and} that 
$
{v}
$
is the global-in-time mild solution to {problem} \eqref{eq;lin}. {If} the constant
$A > 0$ satisfies $A >{A^\ast}${, then} ${v}\equiv 0$ is unstable in the Lyapunov sense.
that is, the constant steady states $u \equiv A$ of problem \eqref{eq;simplified-ar-ks} for $A > A^\ast$ is unstable.
\end{theorem}

\begin{proof}
Fix $\delta \in (0,1)$, and set ${v}_0 \in L^p(\mathbb{R}^n)$, which satisfies $\|{v}_0\|_p=1$ and is chosen later. 
Applying Propositions \ref{prop;812} and \ref{prop_82} yields 
a unique local-in-time mild solution 
$
 {v} 
 \in 
 C([0,T_{\text{max}});L^p(\mathbb{R}^n))
$ 
to problem \eqref{eq;lin} corresponding to 
the initial data $\delta {v}_0$ for any $A > 0$,
where $T_{\text{max}} > 0$ is the maximal existence time of ${v}$ {to problem \eqref{eq;lin}}. 

Let
$
 n
 <
 q
 \le
 \infty
$ 
and assume $T_{\textrm{max}} = \infty$. 
Then we introduce the numbers $T > 0$ and $T_{*} > 0$ by
\begin{align*}
 T 
 \coloneqq \, & 
 \sup 
 \left\{ 
  t 
  \, 
  \mid 
  \,
  \|
   {v}(t)
   -
   \mathrm{e}^{t\mathcal{L}_A}
   (
    \delta {v}_0
   )
  \|_p
  \le 
  \frac{\delta}{2}\mathrm{e}^{M\tau}, 
  \quad 
  0 \le \tau \le t
 \right\}, \\
 T_{*}
 \coloneqq \, &
 \frac{1}{M}
 \log 
  \bigg( 
   \frac{2}{\delta} 
  \bigg),
\end{align*}
where {the constant} $M = M(A) > 0$ is defined in Proposition \ref{prop;ess-sp-li}.
If $T>T_{*}$ or $T=\infty$, {then} 
$v \equiv 0$ is unstable in the Lyapunov sense. 
{Indeed, we apply Lemma \ref{lem;spec-appr} to} choose ${v}_{0,\gamma} \in L^p(\mathbb{R}^n)$ with 
$
 \|
  {v}_{0,\gamma}
 \|_p
 =
 1
$
satisfying 
\begin{equation}\label{es710}
 \|
  \mathrm{e}^{T_*\mathcal{L}_A}(\delta {v}_{0,\gamma})
  -
  \mathrm{e}^{MT_*}(\delta {v}_{0,\gamma})
 \|_p
 \le 
 \gamma
 \|
  \delta
  {v}_{0,\gamma}
 \|_p
 =
 \gamma
 \delta
\end{equation}
for each $0< \gamma \le 1$. 
According to the definition of $T$ and inequality \eqref{es710}, 
we have
\[
\begin{split}
 \|
  {v}(T_*)
 \|_p
 = \, &
 \|
  \mathrm{e}^{MT_*}
   (
    \delta {v}_{0,\gamma}
   )
   + 
   (
    \mathrm{e}^{T_*\mathcal{L}_A}
     (
      \delta {v}_{0,\gamma}
     )
     -
    \mathrm{e}^{MT_*}
     (
      \delta {v}_{0,\gamma}
     )
   )
   +
   (
    {v}(T_*)
    -
    \mathrm{e}^{T_*\mathcal{L}_A}
     (
      \delta {v}_{0,\gamma}
     )
   )
 \|_p \\
 \ge \, &
 \|
  \mathrm{e}^{MT_*}
   (
    \delta {v}_{0,\gamma}
   )
 \|_p
 -
 \|
  \mathrm{e}^{T_*\mathcal{L}_A}
   (
    \delta {v}_{0,\gamma}
   )
  -
  \mathrm{e}^{MT_*}
   (
    \delta {v}_{0,\gamma}
   )
 \|_p
 -
 \|
  {v}(T_*)
  -
  \mathrm{e}^{T_*\mathcal{L}_A}
   (
    \delta {v}_{0,\gamma}
   )
 \|_p\\
 \ge \, &
 \delta 
 \mathrm{e}^{MT_*}
 -
 \gamma
 \delta
 -
 \frac{\delta}{2}
 \mathrm{e}^{MT_*}
 =
 1
 -
 \gamma
 \delta
 >
 1
 -
 \gamma.
\end{split}
\]
In particular, 
for 
$
 \gamma
 =
 1/2
$,
\[
 \|
  {v}(T_*)
 \|_p
 >
 \frac{1}{2}
\]
holds. {On the other hand, s}uppose that 
$
 T
 \le 
 T_*
$. 
Then{, by the argument similar to that in the proof of \eqref{es_833}, there is an exponent $n<q\le \infty$ such that} the mild solution ${v}$ to {problem} \eqref{eq;lin} corresponding to the initial data $\delta {v}_{0,\gamma}$ {is estimated as follows:}
\begin{equation}\label{es811}
\begin{split}
  & \, 
 \|
  {v}(t)
  -
  \mathrm{e}^{t\mathcal{L}_A}
  (
   \delta {v}_{0,\gamma}
  )
 \|_p\\
 \le \, & 
 \int_0^t
  \|
   \nabla 
   \cdot 
   \mathrm{e}^{(t-s)\mathcal{L}_A}
   (
    {v}\nabla K*{v}
   )(s)
  \|_p
 \, \mathrm{d}s\\
 \le \, & 
  C_0
  \int_0^t
   (t-s)^{-\frac{n}{2q}-\frac12}
   \mathrm{e}^{\frac{3}2M(t-s)}
   \|
    {v}(s)
   \|_p^2
   \, \mathrm{d}s \\
 \le \, & 
  C_0
  \int_0^t
   (t-s)^{-\frac{n}{2q}-\frac12}
   \mathrm{e}^{\frac32 M(t-s)}
    (
     \|
      {v}(s)
      -
      \mathrm{e}^{s\mathcal{L}_A}
      (
       \delta {v}_{0,\gamma}
      )
     \|_p^2
     +
     \|
      \mathrm{e}^{s\mathcal{L}_A}
      (
       \delta {v}_{0,\gamma}
      )
     \|_p^2
    )
    \, \mathrm{d}s \\
\le \, &
 C_0
  \int_0^t
   (t-s)^{-\frac{n}{2q}-\frac12}
   \mathrm{e}^{\frac32 M(t-s)}
    \left(
     \frac{\delta^2}{4}
     \mathrm{e}^{2Ms}
     +
     4
     \delta^2
     \mathrm{e}^{2Ms}
     \|
      {v}_{0,{\gamma}}
     \|_p^2
    \right)
   \, \mathrm{d}s  \\
  \le \, & 
  C_0
   \delta^2\mathrm{e}^{\frac32 Mt}
   \int_0^t
    (t-s)^{-\frac{n}{2q}-\frac12}
    \mathrm{e}^{\frac12 Ms}
   \, \mathrm{d}s 
\end{split}
\end{equation}
for some 
$
 C_0
 >
 0
$
independent of $t$, $T$. {Here we have used the definition of $T$ and inequality \eqref{76}.}
{To complete estimate \eqref{es811}, we} will now show that there exists {a constant} $C_1 > 0$ independent of $T > 0$ such that
\begin{equation}\label{es811-2}
 \int_0^t
  (t-s)^{-\frac{n}{2q}-\frac12}
  \mathrm{e}^{\frac12 Ms}
 \, \mathrm{d}s
 \le 
 C_1
 \mathrm{e}^{\frac12 Mt}
\end{equation}
for {all} $0 \le t \le T$.
{For this purpose}, take $\eta > 0$ so that 
$
 0
 <
 \eta
 <
 t
 \le 
 T
$.
Then{,} for $\eta > 0$ and $n < q \le \infty$, we have
\[
\begin{split}
 \int_0^t
  (t-s)^{-\frac{n}{2q}-\frac12}
  \mathrm{e}^{\frac12 Ms}
 \, \mathrm{d}s
 \le \, &
 \int_0^{t-\eta}
  (t-s)^{-\frac{n}{2q}-\frac12}
  \mathrm{e}^{\frac12 Ms}
 \, \mathrm{d}s
 +
 \int_{t-\eta}^t
  (t-s)^{-\frac{n}{2q}-\frac12}
  \mathrm{e}^{\frac12 Ms}
 \, \mathrm{d}s \\
 \le \, & 
 \eta^{-\frac{n}{2q}-\frac12}
 \int_0^{t-\eta}
  \mathrm{e}^{\frac12 Ms}
 \, \mathrm{d}s
 +
 \mathrm{e}^{\frac12 Mt}
 \int_{t-\eta}^t
  (t-s)^{-\frac{n}{2q}-\frac12}
 \, \mathrm{d}s \\
 \le  \, &
 \frac{2}{M}
 \eta^{-\frac{n}{2q}-\frac12}
 \mathrm{e}^{\frac12 Mt}
 +
 \frac{2q}{q-n}
 \eta^{\frac12-\frac{n}{2q}}
 \mathrm{e}^{\frac12 Mt}\\
 \le \, & 
 C_1
 \mathrm{e}^{\frac12 Mt},
\end{split}
\]
where 
\[
 C_1
 \coloneqq
 \min_{\eta > 0}
 \bigg(
 \frac{2}{M}
 \eta^{-\frac{n}{2q}-\frac12}
 +
 \frac{2q}{q-n}
 \eta^{\frac12-\frac{n}{2q}}
 \bigg)
 >
 0
\]
is independent of $t, T$.
Therefore, we observe from estimates \eqref{es811} and \eqref{es811-2} that
\begin{equation}\label{es813}
 \|
  {v}(t)
  -
  \mathrm{e}^{t\mathcal{L}_A}
  (
   \delta {v}_{0,\gamma}
  )
 \|_p
 \le 
 C_0
 C_1
 \delta^2
 \mathrm{e}^{2Mt}
\end{equation}
holds for {any} $0 \le t \le T$ and {some} $C_0, C_1 > 0$ independent of $t, T$.
According to the definition of $T$ and estimate \eqref{es813}, 
we obtain
\[
 \frac{\delta}{2}
 \mathrm{e}^{MT}
 =
 \|
  {v}(T)
  -
  \mathrm{e}^{T\mathcal{L}_A}
  (
   \delta {v}_{0,\gamma}
  )
 \|_p
 \le 
 C_0C_1
 \delta^2
 \mathrm{e}^{2MT}.
\]
Hence, the inequality
$
 1
 \le 
 2C_0C_1
 \delta \mathrm{e}^{MT}
$
holds for fixed $\delta > 0$.
In particular, 
there exists $T_{**} \le T$ depending on $\delta > 0$ such that
$
 1
 =
 2C_0C_1
 \delta 
 \mathrm{e}^{MT_{**}}
$.
Indeed, 
$
 1 
 \le 
 2C_0C_1
 \delta 
 \mathrm{e}^{MT}
$ 
is equivalent to 
$
 \frac1M
 \log
 \big(
  \frac{1}{2C_0C_1\delta}
 \big)
 \le
 T
$,
which claims that
$
 T_{**} = T
$
or
$
 \frac1M
 \log
 \big(
  \frac{1}{2C_0C_1\delta}
 \big)
 =
 T_{**}
 <
 T
$.
Thus, by estimate \eqref{es813}
and
$
 1/(2C_0C_1)
 =
 \delta 
 \mathrm{e}^{MT_{**}}
$, 
we have
\begin{equation}\label{es814}
\begin{split}
 \|
  {v}(T_{**})
 \|_p
 = \, &
 \|
  \mathrm{e}^{T_{**}\mathcal{L}_A}
   (
    \delta {v}_{0,\gamma}
   )
  +
  {v}(T_{**})
  -
  \mathrm{e}^{T_{**}\mathcal{L}_A}
   (
    \delta {v}_{0,\gamma}
   )
 \|_p \\
 \ge \, &
 \|
  \mathrm{e}^{T_{**}\mathcal{L}_A}
   (
    \delta {v}_{0,\gamma}
   )
 \|_p
 -
 \|
  {v}(T_{**})
  -
  \mathrm{e}^{T_{**}
  \mathcal{L}_A}
   (
    \delta {v}_{0,\gamma}
   )
 \|_p \\
 \ge \, &
 \|
  \mathrm{e}^{T_{**}\mathcal{L}_A}
   (
    \delta {v}_{0,\gamma}
   )
 \|_p
 -
 \big(
  \sqrt{
  C_0
  C_1
  }
  \delta 
  \mathrm{e}^{MT_{**}}
 \big)^2 \\
 = \, &
 \|
  \mathrm{e}^{T_{**}\mathcal{L}_A}
   (
    \delta {v}_{0,\gamma}
   )
 \|_p
 -
 \frac{1}{4C_0C_1}.
\end{split}
\end{equation}
{Thus}, by taking $k_0 \gg 1$ so that
\[
 \gamma
 =
 \frac{1}{4k_0C_0C_1}
 \in 
 (0,1),
\]
it follows from estimates \eqref{75} and \eqref{es814} that
\[
\begin{split}
 \|
  {v}(T_{**})
 \|_p
 \ge \, &
 \|
  \mathrm{e}^{MT_{**}}
   (
    \delta 
    {v}_{0,\gamma}
   )
   +
   \mathrm{e}^{T_{**}\mathcal{L}_A}
    (
     \delta {v}_{0,\gamma}
    )
   -
   \mathrm{e}^{MT_{**}}
    (
     \delta {v}_{0,\gamma}
    )
 \|_p
 -
 \frac{1}{4C_0C_1}\\
 \ge \, &
 \|
  \mathrm{e}^{MT_{**}}(\delta {v}_{0,\gamma})
 \|_p
 -
 \|
  \mathrm{e}^{T_{**}\mathcal{L}_A}
   (
    \delta {v}_{0,\gamma}
   )
  -
  \mathrm{e}^{MT_{**}}
   (
    \delta {v}_{0,\gamma}
   )
 \|_p
 -
 \frac{1}{4C_0C_1}\\
 \ge \, &
 \delta 
 \mathrm{e}^{MT_{**}}
 -
 \gamma
 \delta 
 -
 \frac{1}{4C_0C_1} \\
 = \, & 
 \frac{1}{4C_0C_1}
 -
 \gamma \\
 = \, &
 \frac{1}{4C_0C_1}
 \bigg( 
  1-\frac{1}{k_0}
 \bigg)
 >
 0.
\end{split}
\]
As a consequence, the proof of Theorem \ref{t83} is complete. 
\end{proof}

\appendix

\section{Remarks on Proposition \ref{prop;ess-sp-li}}
\subsection{Positive solutions of the algebraic equation $g_A'(\tau)=0$}
We have shown in  Proposition \ref{prop;ess-sp-li} that there is a positive constant $\tau^\ast$ with 
\[
  \tau^\ast
  >
  \tau_0
  =
  \frac{
   \lambda_2\lambda^{\frac{1}{3}}
   (
    \lambda^{\frac{2}{3}}
    -
    \beta^{\frac{1}{3}}
    )
    }{
     (
      \beta
      \lambda)^{\frac{1}{3}}
      -
      1
     }
     >
     0
\] 
such that the function $g_A(\tau)$ on $\tau\ge 0$ given by \eqref{fn;ga} has a positive maximum at $\tau=\tau^\ast$
if $A>A^\ast$, where the constant $A^\ast$ is defined by \eqref{eq;cond-A-th}.
The aim of Appendix A is to give the explicit form of $\tau^\ast$, which is a positive solution to the algebraic equation $g_A'(\tau)=0$, that is 
\begin{equation}\label{eq;quar}
 (\lambda_1+\tau)^2(\lambda_2+\tau)^2
 =
 A[\beta_1\lambda_1(\lambda_2+\tau)^2-\beta_2\lambda_2(\lambda_1+\tau)^2].
\end{equation}
Changing the variable $\tau$ to 
\[
 s 
 \coloneqq 
 \tau
 +
 \dfrac{\lambda_1 + \lambda_2}{2}.
\]
Then equation \eqref{eq;quar} is equivalent to the equation
\begin{equation}\label{eq;tr-qu}
 s^4
 +p
 s^2
 + 
q
 s
 +
 r
 =
 0, 
\end{equation}
where 
\[
\begin{split}
p:=&-\frac{1}{2}[(\lambda_1-\lambda_2)^2+2(\beta_1\lambda_1-\beta_2\lambda_2)A], \\
q:=&A(\lambda_1-\lambda_2)(\beta_1\lambda_1+\beta_2\lambda_2), \\
r:=&\frac{1}{16}\left[(\lambda_1-\lambda_2)^2-4(\beta_1\lambda_1-\beta_2\lambda_2)A\right].
\end{split}
\]
{We apply Ferrari's method  to solve equation \eqref{eq;tr-qu}.
Suppose that $y \in \mathbb{R}$ satisfies 
relation
\begin{equation}\label{eq;fe-1}
 (s^2+y)^2
 =
 (2y-p)s^2
 -qs
 +(y^2-r)
\end{equation}
and ralation
\begin{equation}\label{eq;cu-1}
 q^2-4(y^2-r)(2y-p)
 =
 0.
\end{equation}
Then there exists $a, b \in \mathbb{R}$ depending only on $p, q, r$ such that
\begin{equation}\label{eq;fe-2}
 (2y-p)s^2
 -
 qs
 +
 (y^2-r)
 =
 (as+b)^2.
\end{equation}
Cardano's formula
gives us explicit solutions to 
cubic equation \eqref{eq;cu-1}, that is,
\begin{equation}\label{eq;cu-sol}
y=
  \frac{p}{6}
  +
 \frac{\left(t_1+ 12\sqrt{t_2}\right)^\frac13}{12}+
 \frac{\left(t_1-12\sqrt{t_2}\right)^\frac13}{12},
\end{equation}
where $t_1:=8p^3+108q^2-432pr$, $t_2:=1296p^2r^2-48p^4r-648pq^2r+12p^3r^2+81q^4$, and we also choose the right hand side of relation \eqref{eq;cu-sol} as a real number.
Combining relation \eqref{eq;fe-1}, relation \eqref{eq;fe-2} and relation \eqref{eq;cu-sol},
we have
\begin{equation}\label{eq;qu-fac}
 (s^2+y)^2
 -
 (as+b)^2
 =
 (s^2+as+y+b)
 (s^2-as+y-b).
\end{equation}
Thanks to factorization \eqref{eq;qu-fac}, 
the solutions to quartic equation \eqref{eq;quar} are explicitly given by
\begin{equation}\label{ans;s}
 s
 =
 \frac{-a\pm \sqrt{a^2-4(y+b)}}{2},\quad
 \frac{a\pm \sqrt{a^2-4(y-b)}}{2}, 
\end{equation}
and substituting relation \eqref{eq;cu-sol} into \eqref{ans;s} 
yields the largest real solution $s^*$ to equation \eqref{eq;tr-qu}.
Therefore we conclude  
\[
 \tau^\ast
 =
 s^*
 -
 \dfrac{\lambda_1+\lambda_2}{2}.
\]
\subsection{The monotonicity of $g_A$}
We deduce the necessary and sufficient condition \eqref{eq;co-A-mo} for the monotonicity of $g_A$ under condition \eqref{eq;ar-ba-co}.}

Let $\beta={\beta_1}/{\beta_2}$, $\lambda={\lambda_1}/{\lambda_2}$. First, assume the condition $A>0$ in the case $\beta\lambda\le 1$ and $\beta\le \lambda$. Then it follows that $g_A''(\tau)\ge 0$ for all $\tau\ge 0$
if $\beta\le \lambda^2$, and the function $g_A'(\tau)$ on $\tau\ge0$ admits a local minimum at 
\[
 \tau
 =
 \tau_0
 =
 \frac{
  \lambda_2\lambda^{\frac{1}{3}}(\beta^{\frac{1}{3}}-\lambda^{\frac{2}{3}})
 }{
  1-(\beta\lambda)^{\frac{1}{3}}
  }
  >
  0
\]
 if $\beta>\lambda^2$. 
 Hence this together with $g_A'(0)<-1$ and $\lim_{\tau\to\infty}g_A'(\tau)=-1$ implies that $g_A'(\tau)<0$ for any $\tau\ge 0$.

Next, assume the condition 
\[
  0
  <
  A
  <
  \frac{
   \lambda_1\lambda_2
   }{
   \beta_1\lambda_2
   -
   \beta_2\lambda_1
   }
\] 
in the case $\beta>\lambda$ and $\beta \ge \lambda^2$. 
Then we find that there is a point $\tau_0$ at which the function $g_A$ has a local minimum if $\beta\lambda<1$, and $g_A''(\tau)\le 0$ for any $\tau\ge 0$ if $\beta\lambda\ge 1$. 
Therefore, by these fact, $g_A'(0)<0$ and $\lim_{\tau\to\infty}g_A'(\tau)=-1$, we obtain $g_A'(\tau)<0$ for all $\tau\ge 0$.

Finally, assume the condition
\[
0 < A <\dfrac{(\lambda_1 - \lambda_2)^2}{\left[(\beta_1\lambda_1)^{\frac13} - (\beta_2\lambda_2)^{\frac13}\right]^3}
\]
in the case $\beta<\lambda^2$ and $\beta\lambda>1$. Then it holds that the function $g_A'(\tau)$ on $\tau\ge 0$ has a local maximum at $\tau=\tau_0$, and 
\[
g_A'(\tau_0)=-1+A\dfrac{\left[(\beta_1\lambda_1)^{\frac13} - (\beta_2\lambda_2)^{\frac13}\right]^3}{(\lambda_1 - \lambda_2)^2}<0.
\]
We also see that $g_A'(0)<0$ due to \eqref{c;a3} and $\lim_{\tau\to\infty}g_A'(\tau)=-1$. Thus these gives $g_A'(\tau)<0$ for any $\tau\ge 0$.

 \section*{Acknowledgments}
H.\ Wakui was supported by JSPS KAKENHI Grant Numbers JP23K19005, JP25K07080.
T.\ Yamada was supported by JSPS Grant-in-Aid for Scientific Research (C) Grant number JP24K06806.

\bibliographystyle{siam}
\bibliography{Bait}

 \end{document}